\RequirePackage{ifpdf}
\ifpdf 
\documentclass[pdftex]{sigma}
\else
\documentclass{sigma}
\fi

\usepackage{pict2e}
\usepackage{verbatim}
\usepackage{array}
\usepackage{float}
\usepackage{feynmp}

\def\one{\mbox{\rm 1}\hskip-2.8pt \mbox{\rm l}}
\newcommand{\nco}{\newcommand}
\nco{\Deltanite}{\ensuremath{\,\,\mathrm{l}\!\!\!1}}
\nco{\ZZ}{\mathbb{Z}}
\nco{\CC}{\mathbb{C}}
\nco{\red}{\color{red}}
\nco{\redend}{\normalcolor}

\def\one{\mbox{\rm 1}\hskip-2.8pt \mbox{\rm l}}

\def\sqedges#1#2#3#4{
\begin{picture}(0,0)(0,-.3)
\put(2,-3.2){\makebox(0,0)[b]{\scriptsize \mbox{$#2$}}}
\put(4.4,0){\makebox(0,0)[l]{\scriptsize \mbox{$#4$}}}
\put(2,2.5){\makebox(0,0)[b]{\scriptsize \mbox{$#1$}}}
\put(-.4,0){\makebox(0,0)[r]{\scriptsize \mbox{$#3$}}}
\end{picture}}
\def\sqface#1#2#3#4#5#6#7{\rule[-2.8\unitlength]{0in}{5.6\unitlength}
\begin{picture}(4,4)(-#6,-#7)
\put(0,-2){\vector(1,0){4}}
\put(4,2){\vector(0,-1){4}}
\put(0,2){\vector(1,0){4}}
\put(0,2){\vector(0,-1){4}}
\put(0,-3.2){\makebox(0,0)[b]{\scriptsize \mbox{$#4$}}}
\put(4,-3.2){\makebox(0,0)[b]{\scriptsize \mbox{$#3$}}}
\put(4,2.5){\makebox(0,0)[b]{\scriptsize \mbox{$#2$}}}
\put(0,2.5){\makebox(0,0)[b]{\scriptsize \mbox{$#1$}}}
\put(2,0){\makebox(0,0){\scriptsize \mbox{$#5$}}}
\end{picture}}
\def\punit#1{\hspace{#1\unitlength}}
\newcommand{\racah}[6]
 {\makebox{}
{
\{{}^{ #1\  #2 \ #3}_{ #4 \ #5 \  #6}\}
 }}

\newcommand{\wigner}[6]
{\makebox{}
{
[{}^{ #1\  #2 \   #3}_{ #4 \  #5 \   #6}]
 }}

\newcommand{\TET}[6]
{\makebox{}
{
TET\,({}^{ #1\  #2 \ #3}_{ #4 \ #5 \  #6})
 }}

\newcommand{\frameI}[4]{\unitlength 0.035cm
\parbox{65pt}{\begin{picture}(65,25)
\put(14,13){\circle*{4}}
\put(50,13){\circle*{4}}
\qbezier(14,13)(32,22)(50,13)
\qbezier(14,13)(32,2)(50,13)
\put(35,18){\vector(1,0){0}}
\put(35,8){\vector(1,0){0}}
\put(0,10){$#1$}
\put(55,10){$#2$}
\put(28,21){$#3$}
\put(28,0){$#4$}
\end{picture}}}

\newcommand{\frameIsingle}[2]{\unitlength 0.035cm
\parbox{55pt}{\begin{picture}(55,15)
\put(14,7){\circle*{4}}
\put(40,7){\circle*{4}}
\put(14,7){\line(1,0){28}}
\put(30,7){\vector(1,0){0}}
\put(0,4){$#1$}
\put(45,4){$#2$}
\end{picture}}}

\newcommand{\frameIno}[2]{\unitlength 0.035cm
\parbox{55pt}{\begin{picture}(55,15)
\put(14,7){\circle*{4}}
\put(40,7){\circle*{4}}
\put(14,7){\line(1,0){28}}
\put(0,4){$#1$}
\put(45,4){$#2$}
\end{picture}}}

\newcommand{\Tri}[4]{\unitlength 0.035cm
\parbox{50pt}{\begin{picture}(50,45)
\put(10,5){\circle*{5}}
\put(40,5){\circle*{5}}
\put(25,35){\circle*{5}}
\put(10,5){\line(1,0){30}}
\put(10,5){\line(1,2){15}}
\put(40,5){\line(-1,2){15}}
\put(10,5){\vector(1,0){17}}
\put(40,5){\vector(-1,2){7.6}}
\put(25,35){\vector(-1,-2){9.7}}
\put(-3,0){$#1$}
\put(45,0){$#2$}
\put(23,40){$#3$}
\put(22,15){$#4$}
\end{picture}}}

\newcommand{\Tribis}[7]{\unitlength 0.035cm
\parbox{50pt}{\begin{picture}(50,45)
\put(10,5){\circle*{5}}
\put(40,5){\circle*{5}}
\put(25,35){\circle*{5}}
\put(10,5){\line(1,0){30}}
\put(10,5){\line(1,2){15}}
\put(40,5){\line(-1,2){15}}
\put(10,5){\vector(1,0){17}}
\put(40,5){\vector(-1,2){7.6}}
\put(25,35){\vector(-1,-2){9.7}}
\put(-2,0){$#1$}
\put(44,0){$#2$}
\put(23,40){$#3$}
\put(22,15){$#7$}
\put(22,-3){$#4$}
\put(38,16){$#5$}
\put(2,16){$#6$}
\end{picture}}}

\newcommand{\Tridouble}[7]{\unitlength 0.035cm
\parbox{50pt}{\begin{picture}(50,45)
\put(10,5){\circle*{5}}
\put(40,5){\circle*{5}}
\put(25,35){\circle*{5}}
\qbezier(10,5)(25,7)(40,5)
\qbezier(10,5)(25,-2)(40,5)
\put(28,7){\vector(1,0){0}}
\put(28,1){\vector(1,0){0}}
\put(10,5){\line(1,2){15}}
\put(40,5){\line(-1,2){15}}
\put(40,5){\vector(-1,2){7.6}}
\put(25,35){\vector(-1,-2){9.7}}
\put(0,0){$#1$}
\put(45,0){$#2$}
\put(23,40){$#3$}
\put(22,9){$#7$}
\put(22,-7){$#4$}
\put(38,16){$#5$}
\put(6,16){$#6$}
\end{picture}}}

\newcommand{\basisBig}[8]{\unitlength 0.030cm
\parbox{70pt}{\begin{picture}(70,40)
\put(7.5,10){\begin{picture}(70,25)
\put(0,0){\circle*{6}}
\put(50,0){\circle*{6}}
\put(16,24){\circle*{6}}
\put(66,24){\circle*{6}}
\put(0,0){\line(1,0){50}}
\put(0,0){\line(2,3){16}}
\put(50,0){\line(2,3){16}}
\put(16,24){\line(1,0){50}}
\put(0,0){\vector(1,0){27.5}}
\put(0,0){\vector(2,3){9.5}}
\put(66,24){\vector(-1,0){27.5}}
\put(66,24){\vector(-2,-3){9.5}}
\put(-13,-7){$#1$}
\put(55,-7){$#4$}
\put(10,30){$#2$}
\put(68,29){$#3$}
\put(-6,12){$#5$}
\put(61,10){$#7$}
\put(23,-9){$#6$}
\put(39,27){$#8$}
\end{picture}}
\end{picture}}}

\newcommand{\frameIIIeq}[9]{\unitlength 0.040cm
\parbox{110pt}{\begin{picture}(100,100)
\put(15,10){\begin{picture}(70,25)
\put(0,0){\circle*{5}}
\put(50,0){\circle*{5}}
\put(16,24){\circle*{5}}
\put(66,24){\circle*{5}}
\put(41,74){\circle*{5}}
\put(16,24){\line(1,2){25}}
\put(16,24){\vector(1,2){12}}
\put(66,24){\line(-1,2){25}}
\put(50,0){\line(-9,74){9}}
\put(44,47){\vector(-1,4){0}}
\put(55,45){\vector(1,-4){0}}
\put(12,45){\vector(-1,-4){0}}
\qbezier(0,0)(0,50)(41,74)
\put(0,0){\line(1,0){50}}
\put(0,0){\line(2,3){16}}
\put(50,0){\line(2,3){16}}
\put(16,24){\line(1,0){50}}
\put(0,0){\vector(1,0){27.5}}
\put(0,0){\vector(2,3){9.5}}
\put(66,24){\vector(-1,0){27.5}}
\put(66,24){\vector(-2,-3){9.5}}
\put(-10,-7){$#1$}
\put(55,-7){$#4$}
\put(10,30){$#2$}
\put(68,29){$#3$}
\put(41,80){$#5$}
\put(12,12){$#6$}
\put(61,10){$#8$}
\put(23,-7){$#7$}
\put(35,27){$#9$}
\put(0,45){$\beta_1$}
\put(16,45){$\beta_2$}
\put(35,45){$\beta_4$}
\put(55,45){$\beta_3$}
\end{picture}}
\end{picture}}}

\newcommand{\basisBigLarge}[8]{\unitlength 0.030cm
\parbox{70pt}{\begin{picture}(70,40)
\put(7.5,10){\begin{picture}(70,25)
\put(0,0){\circle*{6}}
\put(50,0){\circle*{6}}
\put(16,24){\circle*{6}}
\put(66,24){\circle*{6}}
\put(0,0){\line(1,0){50}}
\put(0,0){\line(2,3){16}}
\put(50,0){\line(2,3){16}}
\put(16,24){\line(1,0){50}}
\put(0,0){\vector(1,0){27.5}}
\put(0,0){\vector(2,3){9.5}}
\put(66,24){\vector(-1,0){27.5}}
\put(66,24){\vector(-2,-3){9.5}}
\put(-13,-8){$#1$}
\put(55,-7){$#4$}
\put(10,30){$#2$}
\put(68,29){$#3$}
\put(-6,12){$#5$}
\put(61,10){$#7$}
\put(23,-7){$#6$}
\put(39,27){$#8$}
\end{picture}}
\end{picture}}}

\newcommand{\Trimod}[3]{\unitlength 0.030cm
\parbox{52pt}{\begin{picture}(52,50)
\put(-1,5){\line(0,1){35}}
\put(51,40){\line(0,-1){35}}
\put(10,5){\circle*{5}}
\put(40,5){\circle*{5}}
\put(25,35){\circle*{5}}
\put(10,5){\line(1,0){30}}
\put(10,5){\line(1,2){15}}
\put(40,5){\line(-1,2){15}}
\put(52,36){${ }^2$}
\put(1,0){$#1$}
\put(45,0){$#2$}
\put(23,40){$#3$}
\end{picture}}}

\begin{document}

\allowdisplaybreaks

\renewcommand{\PaperNumber}{099}

\FirstPageHeading

\ShortArticleName{$SU(3)$ Triangular Cells}

\ArticleName{Notes on TQFT Wire Models and Coherence\\ Equations for $\boldsymbol{SU(3)}$ Triangular Cells}

\Author{Robert COQUEREAUX~$^\dag$, Esteban ISASI~$^\ddag$ and Gil SCHIEBER~$^\dag$}

\AuthorNameForHeading{R.~Coquereaux, E.~Isasi and G.~Schieber}

\Address{$^\dag$~Centre de Physique Th\'eorique (CPT)  Luminy, Marseille, France}
\EmailD{\href{mailto:Robert.Coquereaux@cpt.univ-mrs.fr}{Robert.Coquereaux@cpt.univ-mrs.fr}, \href{mailto:schieber@cpt.univ-mrs.fr}{schieber@cpt.univ-mrs.fr}}

\Address{$^\ddag$~Departamento de F\'{\i}sica, Universidad Sim\'on Bol\'{\i}var, Caracas, Venezuela}
\EmailD{\href{mailto:eisasi@usb.ve}{eisasi@usb.ve}}

\ArticleDates{Received July 09, 2010, in f\/inal form December 16, 2010;  Published online December 28, 2010}

\Abstract{After a summary of the TQFT wire model formalism we bridge the gap from Kuperberg equations for $SU(3)$ spiders to Ocneanu coherence equations for systems of triangular cells on  fusion graphs that describe modules associated with the fusion category of $SU(3)$ at level $k$. We show how to solve these equations in  a number of examples.}


\Classification{81R50; 81R10; 20C08; 18D10}

\section{Foreword}\label{section1}

Starting with the collection of  irreducible integrable representations (irreps) of $SU(3)$ at some level $k$ (constructed  in the framework of af\/f\/ine algebras or in the framework of quantum groups at roots of unity), the problem is to decide whether a graph encoding the action of these irreps
actually def\/ines a ``healthy fusion graph''  associated with a bona f\/ide $SU(3)$ nimrep, i.e.\ a~module over a particular kind of fusion category.
This is done by associating a complex number (a ``triangular cell'') with every elementary triangle of the given graph in such a way that their collection (a ``self-connection'')  obeys a system of   non trivial quadratic and quartic equations, called Ocneanu coherence equations, that can be themselves derived from  another set of equations (sometimes called Kuperberg equations) describing relations between the intertwiners of the underlying fusion category. One issue is to describe and derive the coherence equations themselves.  Another issue is to use them on the family of examples giving rise to $SU(3)$ nimreps.
Both problems are studied in this article.

Our paper consists of two largely independent parts.
The f\/irst is a set of notes dealing with TQFT graphical models (wire models, spiders, etc.).
Our motivation for this part was to show how the Ocneanu coherence equations for triangular cells of fusions graphs could be deduced from the so-called Kuperberg relations for $SU(3)$,
something that does not seem to be explained in the literature. This section was then enlarged in order to set the discussion in the larger framework of graphical TQFT models
and to discuss some not so well known features that show up when comparing the $SU(2)$ and $SU(3)$ situations.
Despite its size, this f\/irst part should not be considered as a general presentation of the subject starting from f\/irst principles, this would require a book, not an article.
The second part of the paper deals about the coherence equations themselves. In particular we show,  on a selection of examples chosen among quantum subgroups (module-categories) of type $SU(3)$, how to solve them in order to obtain a  self-connection on the corresponding fusion graphs.

\newpage

\section{Framework}\label{section2}

\subsection{Introduction and history}\label{section2.1}

\looseness=1
In recent years the mathematical structure of WZW models has been understood in the framework of fusion categories (also called monoidal), and more generally, in terms of module-categories associated with a given fusion category of type $G$ (a compact Lie group) at level~$k$ (some non-negative integer). The former can be constructed, for instance, in terms of af\/f\/ine Lie algebras at level~$k$, or of quantum groups at a root of unity $q$ that depends both on~$k$ and $G$. The usual notion of action of a group acting on a space is then replaced by a notion describing the action of a fusion category (say ${\mathcal A}_k(G)$) on another category, say ${\mathcal E}$, which is not monoidal in general, and which is called a module-category. To every such action is also associated a~partition function, which is modular invariant, and is characterized by a matrix acting in the vector space spanned by the simple objects of ${\mathcal A}_k(G)$.
If one chooses ${\mathcal E}= {\mathcal A}_k(G)$, this matrix is the identity matrix, and the theory is called diagonal. One obvious problem is to classify the possible~${\mathcal E}$ on which a given ${\mathcal A}_k(G)$ can act. Two dif\/ferent~${\mathcal E}$ can sometimes be associated with the same partition function, so that classifying modular invariants and classifying module-categories of type~$G$ (something that is usually described in terms of appropriate graphs) are two dif\/ferent problems. In the case of $A_1 = su(2)$ this distinction can be forgotten and the classif\/ication, which was shown to coincide with the classif\/ication of ADE Dynkin diagrams, was obtained more than twenty years ago~\cite{CIZ}.
The classif\/ication of $A_2=su(3)$ modular invariant partition functions was obtained by~\cite{Gannon},  and graphs that can now be interpreted as graphs describing the corresponding module-categories
were discovered by~\cite{DiFZuber}, see also~\cite{YellowBook, OujdaRDGH, GilDahmaneHassan, RobertGilSL3Categories}. However, obtaining a graph describing some module action for a~given fusion ring may not be enough to
guarantee the existence of the category ${\mathcal E}$ itself.  Loosely speaking, it may be that the relations existing between the simple objects of ${\mathcal A}_k$  will be correctly represented at the level of ${\mathcal E}$, but there may be also non trivial  identities involving the morphisms of  ${\mathcal A}_k$  and  one should check that such  ``relations between relations''  are also correctly represented in ${\mathcal E}$. For this reason, the list of candidates proposed in  \cite{DiFZuber} was later slightly amended in \cite{Ocneanu:Bariloche}.
 In the case of  the monoidal category ${\mathcal A}_k$
              def\/ined by $SU(3)$ at level $k$, there are two non trivial
              identities for the intertwiners. These identities seem  to belong to the
              folklore and have been rediscovered several times.
              They were obviously used in 1993  by H.~Ewen and O.~Ogievetsky in their classif\/ication of quantum spaces of dimension~$3$ (see~\cite{EwenOgievetsky}, and also~\cite{OgievetskyBariloche}), but they were written down explicitly in 1996, in
              terms of  webs, within the framework of $A_2$ spiders, by G.~Kuperberg~\cite{Kuperberg},  and are sometimes called Kuperberg identities for $SU(3)$.
              We shall later return to them, but here
              it is enough to say that they give rise, at
              the level of the module-categories over ${\mathcal A}_k$,
              to a set of equations between members of a particular
              class of $6j$ symbols, called triangular cells.  These
              equations have been worked out by A.~Ocneanu and
              presented, among other results, in a lecture given at
              Bariloche in January 2000.
              Given a candidate graph,
              they  provide a necessary and
              suf\/f\/icient condition for the existence of the underlying
              $SU(3)$ module-category. They have been used to  disregard one of the
              conjectured examples appearing in the slightly
              overcomplete list previously obtained by~\cite{DiFZuber} and produce the results stated in~\cite{Ocneanu:Bariloche}.
              The details, and the coherence equations themselves, were unfortunately not made available.
              To our knowledge, the coherence equations were f\/irst presented and discussed in Chapter~3.5 of the thesis~\cite{Dahmane:Thesis} (2007), a section that was based, in parts, on a set of notes,
              containing the equations and studying several examples, that was written and
              distributed to students by one of us (R.C., reference~[30] of~\cite{Dahmane:Thesis}).
              A number of other examples was then worked out, about three years ago,  by the authors, originally with no intention to publish them. Some time later,  it became clear that no general presentation of the subject was going to be made available by the author of~\cite{Ocneanu:Bariloche} and we started to work on this manuscript.
              We decided to incorporate a tentatively pedagogical discussion of TQFT graphical models, and stressing the dif\/ferences between the $SU(2)$ and $SU(3)$ cases,
              without relying too much on the not so well-known formalism and terminology of spiders,
              since these graphical models provide a convenient framework to discuss
             Kuperberg equations together with coherence equations for Ocneanu
              triangular cells.
   Our aim was to show
             how to derive these equations and how to use them by selecting a class of examples.
             This explains the genesis of the present work, which, if not for the typing itself, was essentially f\/inished two years ago.
             In the process of comple\-ting our article we became aware of a recent work by D.~Evans et al.~\cite{EvansSU3Cells}:  This publication uses Kuperberg equations and Ocneanu coherence equations as an input, solves them in almost all cases of type $SU(3)$, and shows that the representations of the Hecke algebra (used in generalized RSOS models) that one obtains as a by-product from the values of the associated triangular cells, have the expected properties. This reference is quite complete and contains tables of results for all the examples that we originally planned to discuss.
 Nevertheless the main focus of~\cite{EvansSU3Cells} is about the corresponding representations of the Hecke algebra and associated Boltzmann weights, it is neither on graphical methods nor on the techniques used to solve the coherence equations.
 The present article should therefore be of independent interest,  f\/irst of all because the specif\/ic aspects of $SU(3)$ graphical models, although known by experts, do not seem to be easily available (we use these techniques to manipulate intertwiners and to derive the coherence equations, something that is not done in \cite{EvansSU3Cells}), and also because we provide many details performing calculations  leading to explicit values for the triangular cells (they agree, up to gauge freedom, with those given in~\cite{EvansSU3Cells}). In that respect, our approach and methodology often dif\/fers from the last quoted reference.  Although it is not strictly necessary, some familiarity with the use of graphical models (see for instance the book~\cite{Kauffman:book}, that mostly deals with the $SU(2)$ situation) may help the reader.
We acknowledge conversations with A.~Ocneanu who introduced and studied about ten years ago (but did not make available) most of the material discussed here, often using a dif\/ferent terminology.
 The derivation of the coherence equations that we shall present in the f\/irst part of this paper is ours but it can be traced back to the oral presentation~\cite{Ocneanu:Bariloche}.

\subsection{From module-categories to coherence equations}\label{section2.2}

             Our f\/irst ingredient is the fusion category ${\mathcal
              A}_k$ def\/ined by $SU(3)$ at level $k$. It can be
              constructed in terms of integrable representations of an
              af\/f\/ine algebra, or in terms of particular irreducible representations
              of the quantum groups $SU(3)_q$ at roots of unity ($q = \exp(i \pi/\kappa)$, with $\kappa = k+3$). The
              Grothendieck ring of this monoidal category is
              called the fusion ring,  or the Verlinde algebra.
              For arbitrary values of the non-negative integer $k$,
              it  has two generators  that are conjugate to one another, and
              correspond classically to the two fundamental
              representations of $SU(3)$.
              The table of structure constants describing the
              multiplication of simple objects by a chosen generator is
              encoded by a f\/inite size matrix of dimension $r \times
              r$, with $r=k(k+1)/2$, which can be interpreted as the
              adjacency matrix of a graph (the  Cayley graph of multiplication by
              this generator) called the fundamental fusion
              graph. More generally, multiplication of simple objects $m \times n = \sum_p N_{mn}^p  p$  is encoded by  fusion coef\/f\/icients and fusion matrices~$N_m$.

              The next ingredient is a category ${\mathcal E}$, not
              necessarily monoidal, on which the previous one,~${\mathcal A}_k$, acts.
              In other words, ${\mathcal E}$ is a module-category over
              ${\mathcal A}_k$ and its Grothendieck (abelian) group is
              a module over the fusion ring.
              The structure constants describing this action are
              encoded by matrices with non-negative integer coef\/f\/icients, or
              equivalently, by oriented graphs. In the case of~$sl(3)$, the two generators are conjugated, so that the
              corresponding Cayley graphs, again called fusion (or quantum) graphs, just dif\/fer by orientation of
              edges,  and one graph is enough to fully characterize
              the action of the fusion ring.
            More generally, the (module) action of a~simple object~$m$ of~${\mathcal A}_k$ on a simple object $a$ of ${\mathcal E}$, i.e.\  $m \times a = \sum_b F_{ma}^b   b$  is encoded by ``annular coef\/f\/icients'' and annular matrices $F_m$.

\looseness=1
            A last ingredient to be introduced is the algebra of quantum symmetries \cite{Ocneanu:talks, Ocneanu:paths}:
given a~module-category ${\mathcal E}_k$ over a fusion category  ${\mathcal A}_k(G)$, one may def\/ine its endomorphism ca\-te\-go\-ry  ${\mathcal O(E)} = {\rm End}_{{\mathcal A}} {\mathcal E}$ which is monoidal, like ${\mathcal A}_k(G)$, and  acts on ${\mathcal E}$ in a way compatible with the~${\mathcal A}$ action. In practice one prefers to think in terms of (Grothendieck)  rings and modules, and use the same notations to denote them.
The ring ${\mathcal O({\mathcal E})}$ is naturally a bimo\-dule on ${\mathcal A}_k(G)$,  so that the simple objects $m,n,\ldots$ of the later (in particular the generators) act on ${\mathcal O({\mathcal E})}$, in particular on its own simple objects $x,y,\ldots$ called ``quantum symmetries'',  in two possible ways. One can therefore associate to each generator of ${\mathcal A}_k(G)$  two matrices with non negative integer entries describing left and right multiplication, and therefore two graphs, called left and right chiral graphs whose union is a non connected graph called the Ocneanu graph. The chiral graphs associated with a generator can themselves be non connected.
 More generally, the structure constants describing the bimodule action are encoded by  ``toric matrices'' $W_{xy}$ (write $m   x   n = \sum_y  (W_{xy})_{mn}  y$) that can be physically interpreted as twisted partition functions (presence of defects, see~\cite{PetkovaZuber:Oc}). The particular matrix $Z=W_{00}$ is the usual partition function.  Expressions of toric matrices for all $ADE$ cases can be found in~\cite{Gil:thesis}.

              The category ${\mathcal A}_k$ is modular (action of
              $SL(2,\ZZ)$), and to every
              module-category ${\mathcal E}$ is asso\-ciated, as discussed above, a quantity
              $Z$, called the modular invariant because it commutes with this
              action.
              One problem, of course, is to classify
              module-categories associated with~${\mathcal A}_k$. Very often, ${\mathcal E}$ is not a
              priori known (one does not even know what are its simple
              objects) but in many cases the associated modular
              invariant $Z$ is known, and there exist techniques (like
              the modular splitting method, that we shall not discuss
              in this paper, see \cite{EstebanGil,RobertGilSU4}) that allow one to ``reverse the machine'',
              i.e.\ to obtain the putative Grothendieck group of
              ${\mathcal E}$ and its module structure over the fusion ring, together
              with its fusion graph, from the knowledge of $Z$.
             As already mentioned, obtaining such a module structure over the fusion ring is usually not enough to
              guarantee the existence of the category ${\mathcal E}$ itself, unless we a priori know that it should exist (for instance in those cases of module-categories obtained by conformal embeddings).

\looseness=1
               For a module-category  ${\mathcal E}$ of type $SU(3)$, one can associate a complex-number, called a triangular cell,  to each elementary triangle of its fusion graph.
                This stems from the following simple argument: in the classical theory, like in the quantum one, the tensor cube of the fundamental representation contains the identity representation.
              Therefore, starting from any irreducible representation of $SU(3)$, say $\lambda$,  if we tensor multiply it three times by the fundamental, we obtain a reducible representation that contains $\lambda$ in its reduction.
              The same remark holds if~$\lambda$ is a simple object of a category ${\mathcal E}$ on which ${\mathcal A}_k={\mathcal A}_k(SU(3))$ acts.
              It means that for any elementary triangle (a succession of three oriented edges) of the ${\mathcal A}_k$ fusion graph, or of the graph describing the module structure of ${\mathcal E}$, we have an intertwiner, but this intertwiner should be proportional to the trivial one since $\lambda$ is irreducible.
              In other words, one associates~-- up to some gauge freedom~-- a complex number to each elementary triangle.
              This assignment is sometimes called\footnote{The terminology stems historically from more general constructions describing symmetries of quantum spaces in the framework of operator algebras,
              see \cite{Ocneanu:Galois, Ocneanu:subfactors}.} ``a self-connection'', ``a connection on the system of triangular cells'', or simply ``a cell system''.
                Independently of the question of determining these numbers explicitly, there are non trivial identities between them because of the
    existence of non trivial relations  (Kuperberg identities) between the $SU(3)$ intertwiners.   One obtains two sets of equations,  respectively quadratic and quartic, relating the triangular cells. These
            equations, that we call coherence equations for triangular cells, are sometimes nicknamed ``small and large pocket equations'' because of the shapes of the polygons (frames of the fusion graph) involved in their writing.
             Up to some gauge choice, these equations determine the values of the cells (there may be more than one solution).  Such a solution can be used, in turn, to obtain a representation of the Hecke algebra.

\subsection{Structure of the paper}\label{section2.3}

Our aim is twofold: to
              explain where the coherence equations come from, and to show explicitly, in a
              number of cases, how these equations are used to determine explicitly the values of the triangular cells.
              The  structure of our paper ref\/lects this goal.

              After a general discussion of graphical models that allows
              one to manipulate intertwining
              operators in a  quite ef\/f\/icient manner  (Section~\ref{section3.1}),
              and a  summary of the $SU(2)$ situation (Section~\ref{section3.2}), we
              describe the identities that are specif\/ic to the
              $SU(3)$ cases, namely the quadratic and quartic Kuperberg
              equations between particular webs of the $A_2$ spider
              (Section~\ref{section3.3}). Then, by ``plugging'' these equations into
              graphs describing an $SU(3)$ theory with boundaries, ie
              some chosen module-category ${\mathcal E}$  over ${\mathcal A}_k(SU(3))$ , we
              obtain Ocneanu equations for triangular cells (Section~\ref{section3.4}).

              Section~\ref{section4} is devoted to a discussion of these equations and to
              examples. We
              take, for ${\mathcal E}$, the fusion category ${\mathcal
              A}_k$ itself and analyze three exceptional\footnote{In the sense that they do not belong to series.}
               examples, namely the three ``quantum subgroups'' of type
              $SU(3)$ possessing self-fusion, which are known to exist at levels~$5$,~$9$ and~$21$. Module-categories
              obtained by orbifold techniques from the ${\mathcal A}_k$  (the
              $SU(3)$ analogs\footnote{They have
              self-fusion when $k$ is equal to $0$ modulo~$3$.} of the $D$ graphs), or exceptional
              modules without self-fusion, as well as other examples of ``quantum modules''
              obtained using particular conjugacies, can be
              studied along the same lines.
                  Finally, Subsection~\ref{section4.7} is devoted to  a ``bad example'', i.e.\ a graph that was a candidate for another module-category of $SU(3)$ at level~$9$, but
              for which the coherence equations simply fail.

              The reader who wants to jump directly to the coherence equations and to our description of a few specif\/ic cases may skip the beginning of our paper and starts his reading in Section~\ref{section3.4}, since the
                             derivation of these equations from f\/irst principles requires a fair amount of material that will not be used in the sequel.

\section{Intertwining operators and coherence equations}\label{section3}
\subsection{ Generalities on graphical models (webs or wire models)}\label{section3.1}

From standard category theory we know that objects should be thought of as ``points'' and morphisms as ``arrows'' (connecting points). Although conceptually nice, this kind of visualization appears to be often inappropriate for explicit calculations. Graphical methods like those used in wire models, combinatorial spiders or topological quantum f\/ield theory (TQFT), happen to be more handy.  A detailed presentation of such models goes beyond the scope of this article, but we nevertheless need a few concepts and ideas that we summarize below.

In graphical models, diagrams are drawn on the plane, or more generally on oriented surfaces. Graphical elements of a model include a f\/inite set of distinguished boundary points, possibly dif\/ferent types of wires (oriented or not)  going through, and possibly dif\/ferent types of vertices, where the wires meet (warning: these internal vertices should not be thought of as boundary points).
The general idea underlying graphical models is that boundary points represent simple objects, not arbitrary ones, in contradistinction with category theory, and that
diagrams should be interpreted as morphisms connecting the boundary points (the later are supposed to be tensor multiplied).  Diagrams without boundary points are just numbers.
Diagrams can be formally added,  or multiplied by complex numbers, and it is usually so that arbitrary diagrams can be obtained as linear sums of more basic diagrams. Diagrams can be joined or composed (in~a~way that is very much dependent of the chosen kind of graphical model), and this ref\/lects the composition of morphisms. Diagrams can be read in various directions, which are equivalent, because of the existence of Frobenius isomorphisms.
A given model may also include relations between some particular diagrams (identities between morphisms) and these relations may incorporate formal parameters (indeterminate).
Finally, partial trace operations (tensor contraction or stitch) are def\/ined.
Considering only one type of boundary is not suf\/f\/icient for the description of module-categories, so that more general graphical models also incorporate ``marks'', labelling the dif\/ferent types of  boundaries.

There exist actually two variants of diagrammatical models.  In the f\/irst version (the primitive or elementary  version), boundary points refers to simple objects that are also generators; all other simple objects will appear  when we tensor multiply the fundamental ones (that correspond classically to the fundamental representations of $G$). In other words all other objects are subfactors of tensor products of generators.

 In  the second version (the ``clasp version'', that we could be tempted\footnote{The situation is reminiscent of what happens in the diagrammatic treatment of quantum f\/ield theories (Feynman diagrams) where graphs can incorporate lines, which may be external or internal, describing either fundamental particles or bound states (like clasps), and vertices describing their interactions.} to call the ``bound state version''), f\/inite strings of consecutive fundamental boundary points are grouped together into boxes (clasps) in order to build arbitrary simple objects.
This version may therefore incorporate new kinds of graphical elements, in particular new kinds of boundary points (external clasps representing arbitrary simple objects), new labels for wires, and new kinds of vertices, but one should remember that these new graphical elements should be, in principle, expressible in terms of the primitive ones.

For def\/initeness we suppose, from now on,  that the underlying category is specif\/ied by a~compact simple Lie group~$G$ and a level~$k$ (possibly inf\/inite). When the level is f\/inite, there is a~f\/inite number of simple objects;  they are  highest weight representations, denoted by a weight~$\lambda$, such that $\langle \lambda, \theta \rangle \leq k$, where $\theta$ is the highest root of~$G$ and the scalar product is def\/ined by the fundamental quadratic form. The morphisms are intertwining operators, i.e.\ elements of ${\rm Hom}_{G,k}$ spaces.
The simple objects are often called ``integrable irreducible representations of~$G$ at level~$k$'', a~terminology that makes reference to af\/f\/ine algebras, but one should remember that the same category can be constructed from quantum groups at a root of unity (actually a root of $-1$) $q=\exp(\frac{i \pi}{\kappa})$, where $\kappa = g+k$ is called altitude, and where $g$ is the dual Coxeter number of~$G$. Like in the classical case, simple objects can be labelled by weights or by Young tableaux (reducing to positive integers in the case $G=SU(2)$).
If the level $k$ is too small,  not all fundamental representations of $G$ necessarily appear, however for $G=SU(N)$ all the fundamental representations ``exist'' (or are integrable, in af\/f\/ine algebra parlance), as soon as $k\geq1$.

\begin{figure}[h]
\centering
\includegraphics{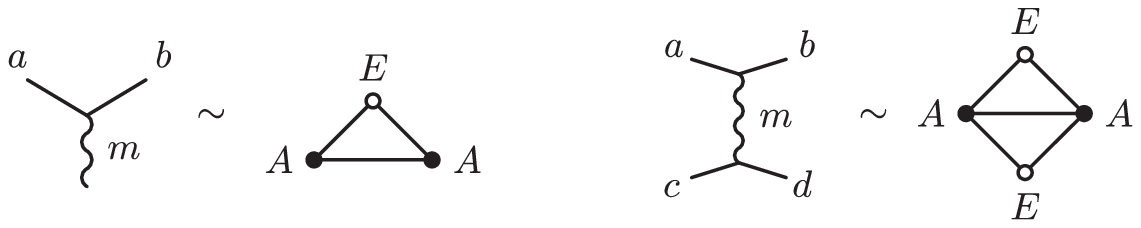}
\caption{Vertex and dif\/fusion graphs for $SU(2)$ coupled to $ADE$ matter.}
\label{fig:vertexsu2}
\end{figure}

The f\/irst graphical models invented to tackle with problems of representation theory go back to the works done, almost a century ago, in mathematics by \cite{RumerTellerWeyl}  and physics or quantum chemistry (recoupling theories for spin,~\cite{Racah}),  but they have been considerably developed in more recent years  by \cite{TemperleyLieb, Jones, Kauffman:book, Lickorish, ReshetikhinTuraev, Kuperberg, Ocneanu:paths, EvansKawa:book} and others.  These models were mostly devoted to the study of $SU(2)$, classical or quantum, in relation with the theory of knots or with the theory of invariants for $3$-manifolds.  They share many features but often dif\/fer by terminology, graphical conventions, signs, and interpretation.
At the end of the 90's this type of  formalism was extended beyond the $SU(2)$ case, \cite{Suciu, Ohtsuki-Yamada,Kuperberg}.  Combinatorial spiders, formally def\/ined in the last reference,  provide a precise terminology and nice graphical models for the study of fusion categories ${\mathcal A}_k(G)$ associated with pairs~$(G,k)$.  In wire models, diagrams are often read from top to bottom (think of this as  time f\/low) but this is just conventional because of the existence of Frobenius maps; spiders diagrams (webs) can also be cut and read in one way or another.
In wire models, two consecutive boundary vertices (or more) located on the same horizontal line (which is not drawn in general) should be tensor multiplied. The same is true, in the case of spiders,  for vertices belonging to a boundary circle.
Moreover the trivial representation is described by an empty drawing (the vacuum).
However, as such, these graphical models are not general enough for our purposes, because although they provide combinatorial descriptions for the morphisms attached to a category  like ${\mathcal A}_k(SU(N))$, they do not describe the collection of inter-relations existing between a monoidal category  ${\mathcal A}_k$, a chosen module ${\mathcal E}$, and the endomorphisms ${\mathcal O({\mathcal E})}$ of the later.
 However, an extension of the recoupling models describing the case of $SU(2)$ ``coupled to matter'' of type $ADE$ (module-categories) was sketched in a section of \cite{Ocneanu:paths}.
The main idea was to introduce dif\/ferent marks (labels called ${\mathcal A}{\mathcal A}$, ${\mathcal A}{\mathcal E}$ or ${\mathcal E}{\mathcal E}$) for boundary points representing simple objects of the three dif\/ferent types, and three kind of  lines (we choose wiggle, continuous, and dashed)  for their morphisms.
Another idea was to choose an alternative way for drawing diagrams, by using elementary star-triangle duality. Let us consider a simple example:
Take $m$ a simple object in ${\mathcal A}$, and $a$ a simple object in ${\mathcal E}$, since the former acts on the later (we have a monoidal functor from ${\mathcal A}$ to the endofunctors of ${\mathcal E}$) we may consider vertices of type~$a$,~$m$,~$b$ where $b$ is a simple  object in  ${\mathcal E}$, or more generally we may consider ``dif\/fusions graphs'' (see Fig.~\ref{fig:vertexsu2}).
Then, by applying a star-triangle duality transformation, we obtain in the f\/irst case a~triangle with vertices marked ${\mathcal E}$, ${\mathcal A}$, ${\mathcal A}$ and edges of type  ${\mathcal A}{\mathcal A}$, ${\mathcal A}{\mathcal E}$, ${\mathcal A}{\mathcal E}$.
 In the second case, we obtain a double triangle (a rhombus) with four edges of type  ${\mathcal A}{\mathcal E}$, and a diagonal of type ${\mathcal A}{\mathcal A}$.
Since the endomorphism category ${\mathcal O}$ of  ${\mathcal E}$ also acts on ${\mathcal E}$, we have also vertices of type~$a$,~$x$,~$b$, where $x$ is a simple object of ${\mathcal O}$, and dually, triangles with edges of type ${\mathcal E}{\mathcal E}$ or ${\mathcal A}{\mathcal E}$, or double triangles with four edges of type  ${\mathcal A}{\mathcal E}$, and a diagonal of type ${\mathcal E}{\mathcal E}$.
One can also introduce vertices of type~$m$,~$n$,~$p$ (for triangles ${\mathcal A}{\mathcal A}{\mathcal A}$) since ${\mathcal A}$ is monoidal, and triangles ${\mathcal E}{\mathcal E}{\mathcal E}$ since ${\mathcal O}$ is monoidal as well. Altogether, there are four kinds of triangles and f\/ive types of  tetrahedra (pairing of double-triangles) describing generalized $6j$~symbols.
This framework seems to be general enough to handle all kinds of situations where a monoidal category  ${\mathcal A}$ is ``coupled'' to a module ${\mathcal E}$. Of course, for any specif\/ic example, we also have specif\/ic relations between the graphical elements of the model. The case of $SU(2)$ and its modules was presented in~\cite{Ocneanu:paths}, but we are not aware of any printed reference presenting graphical models to study $SU(3)$ coupled to its modules.  Fortunately, for the purposes of the present paper, the amount of material that we need from such a model is rather small.
Although we are mainly interested in $SU(3)$, we shall actually start with a description of the $SU(2)$ case since it allows us to present most concepts in a simpler framework.

\subsection[Equations for the $SU(2)$ model]{Equations for the $\boldsymbol{SU(2)}$ model}\label{section3.2}

\subsubsection{The fundamental version}\label{section3.2.1}
\begin{figure}[h]
\centering
\includegraphics{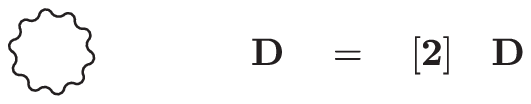}
\caption{The circle for $SU(2)$.}
\label{fig:SU2circle}
\end{figure}

For $(SU(2),k)$ alone, i.e.\ without considering any action of the category ${\mathcal A}_k$ on a module, and calling $\sigma_p$ the simple objects, we have one fundamental object $\sigma=\sigma_1$ corresponding classically to the two-dimensional spin $1/2$ representation.
There is a single non trivial  intertwiner (the determinant map) which sends $\sigma \otimes \sigma$ to the trivial representation. This is ref\/lected in by the fact that $\sigma_0$ appears on the r.h.s.\ of the
equation $\sigma_1 \otimes \sigma_1 = \sigma_0 \oplus \sigma_2$, where $\dim(\sigma_p) = [p+1]$ where
$[n] = \frac{q^n-q^{-n}}{q-q^{-1}}$ with $q=\exp(i\pi/(2+k))$.
The corresponding projector (antisymmetric) is the Jones projector $\sigma \otimes \sigma \mapsto \CC \subset \sigma \otimes \sigma$.
The Frobenius isomorphisms lead to the well known existence of an isomorphism\footnote{Actually two, see the comment in Section~\ref{section3.2.5}.} between $\sigma$ and its conjugate $\overline \sigma$. For this reason, the graphical model uses only one type of wire, which is unoriented.
The fundamental intertwiner $C$ is graphically described by the ``cup diagram'': read from top to bottom, it indeed tells us that $\sigma \otimes \sigma$ contains the trivial representation (remember that the trivial representation is not drawn).
Its adjoint~$C^\dag$ (the same diagram read bottom-up) is the ``cap diagram'', so that the composition of the two gives either an operator  $U=C^\dag C$ from $\sigma \otimes \sigma$ to itself (the Jones projector, up to scale, displayed as a~``cup-cap''), or a number $\beta = C C^\dag$ displayed as a ``cap-cup'', i.e.\ a~circle (closed loop). Obviously $U^2 = \beta U$. It is traditional to introduce the Jones projector $e=\beta^{-1} U$, so that $e^2=e$ and to set $U_1=U$, $U_2= \one \otimes U$, $U_n= \one \otimes \cdots \otimes \one \otimes U$ and $e_n=\beta^{-1} U_n$. A standard calculation shows that  $\beta$ is equal to the $q$-number~$[2]$.
The graphical elements of the $SU(2)$ primitive model are therefore rather simple since we have only one relation (see Fig.~\ref{fig:SU2circle}, where the symbol $D$ refers to any diagram) involving a single parameter, the value of the circle.

\subsubsection{The clasp version (bound states version)}\label{section3.2.2}

The clasp version of the ${\mathcal A}_k(SU(2))$ model is slightly more involved: the representation $\sigma^n = \sigma^{\otimes n}$ is not a simple object and  according to the general philosophy of graphical models, in their primitive version, being a tensor product of $n$ simple objects, it is described by $n$ points along a boundary line. Its endomorphism algebra (operators commuting with the action of~$SU(2)$ or~$SU(2)_q$ on~$\sigma^n$) is the Temperley--Lieb algebra. It is generated by linear combination of diagrams representing crossingless matching of the $2n$ points ($n$ points ``in'' and $n$ points ``out'' in the wire model,  or $2n$ points along a circle, in the spider model), and generated, as a unital algebra, by cup-cap elements $U_i$  (pairs of U-turns in position $(i, i+1)$). The standard recurrence relation for $SU(2)$, namely $\sigma  \otimes \sigma_n = \sigma_{n-1}\oplus \sigma_{n+1}$  (Tchebychev) shows that if $n<k$ there exists a non-trivial intertwiner from $\sigma \otimes \sigma_n$ to $\sigma_{n+1}$ where $\sigma_n$ denotes the irreducible representation of quantum dimension $[n+1]$.
The corresponding projector (the Wenzl projector~$P_n$, which is symmetric) is therefore obtained as the equivariant projection of $\sigma^n$ to its highest weight irreducible summand $\sigma_n$ which is described by a clasp of size $n$ or by a vertical line carrying\footnote{Warning: some authors prefer to denote $\sigma^n$, not $\sigma_n$,  with this drawing.} a label $n$.
 Its expression is given by the Wenzl recurrence formula: $P_1=\one$, $P_2=\one-U_1/\beta$,   $P_{n} = P_{n-1} - \frac{[n-1]}{[n]} P_{n-1}   U_{n-1} \, P_{n-1}$. On the r.h.s.\ of this equation, $P_{n-1}$ is understood as $P_{n-1} \otimes \one$.
The trace of $P_n$,  described by a loop carrying the label $n$, is $ \mu_n = {\rm Tr}(P_n) = [n+1]$.

\begin{figure}[h]
\centering
 \includegraphics{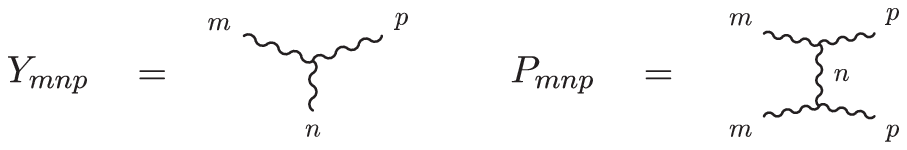}
\caption{The $Y$ and $P$ intertwiners for $SU(2)$.}
\label{fig:YmnpPmnp}
\end{figure}

The non-trivial intertwiners $Y_{mnp}: \sigma_m \otimes \sigma_p \mapsto \sigma_n$ are described by Y-shaped star diagrams (see Fig.~\ref{fig:YmnpPmnp})  with composite wires carrying representation indices ($SU(2)$ weights), or dually, by triangles, where edges are clasps.
It is easy to show that the dimension of the triangle spaces is equal to $0$ or $1$ for all $n$ since all matrix coef\/f\/icients $(N_n)_{mp}$ of fusion matrices obtained from the relation $N_n=N_{n-1} N_1 - N_{n-2}$, where $N_1$ is the adjacency matrix of the fusion graph $A_{k+1}={\mathcal A}_k$ are either equal to $0$ or $1$.  When the value is not $0$, the triangle $(m,n,p)$ is called admissible. It can also be associated with an essential path of length $n$,  from the vertex $\sigma_m$ to the vertex $\sigma_p$, on the fusion graph (see the discussion in Section~\ref{section3.2.4}).

\begin{figure}[h]
\centering
\includegraphics{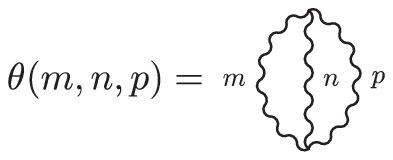}
\caption{The theta symbol for $SU(2)$.}\label{fig:thetamnp}
\end{figure}

 The corresponding endomorphisms $P_{mnp} :  \sigma_m \otimes \sigma_p \mapsto \sigma_n \subset  \sigma_m \otimes \sigma_p$ are displayed as dif\/fusion graphs of a special kind (see Fig.~\ref{fig:YmnpPmnp}), or dually, as particular double triangles.
The trace of~$P_{mnp}$, called theta symbol $\theta(m,n,p)$ because of its shape, is represented in Fig.~\ref{fig:thetamnp}. There are general formulae for the values of these symbols, which, for the pure $SU(2)$ case, are symmetric in~$m$,~$n$,~$p$, see for example~\cite{Kauffman:book};  they are obtained by decomposing $Y_{mnp}$ into elementary tangles or along a web basis.

\begin{figure}[h]
\centering
\includegraphics{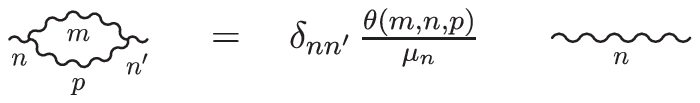}
\caption{Normalizing coef\/f\/icient for a $SU(2)$ loop.}
\label{fig:verticallinewithloop}
\end{figure}

Composing a $Y$ intertwiner and its adjoint in the opposite order gives a ``propagator with a~loop'' which has to be proportional to the identity morphism of  $\sigma_n$ since the later is irreducible, so that, by evaluating the trace on both sides, one f\/inds the identity displayed\footnote{To save space, we display it horizontally.} on Fig.~\ref{fig:verticallinewithloop}.
\noindent
The endomorphisms $\frac{ \mu_n}{\theta(m,n,p)}   P_{mnp}$ are projectors, as it is clear by composing the $P$'s vertically: $P_{mnp} P_{mnp}= \theta(m,n,p)     \mu_n^{-1}   P_{mnp}$.
Notice the particular case $\theta(k-1,1,k) =  \mu_k = [k+1]$ so that we recover  the Wenzl projectors $P_k=P_{k-1,k,1}$, for instance  $P_2 = Y_{121}^\dag Y_{121}=P_{121}$.
As a special case, one recovers  $U=U_1=Y_{101}^\dag Y_{101}=P_{101}$.

\begin{figure}[h]
\centering
\includegraphics{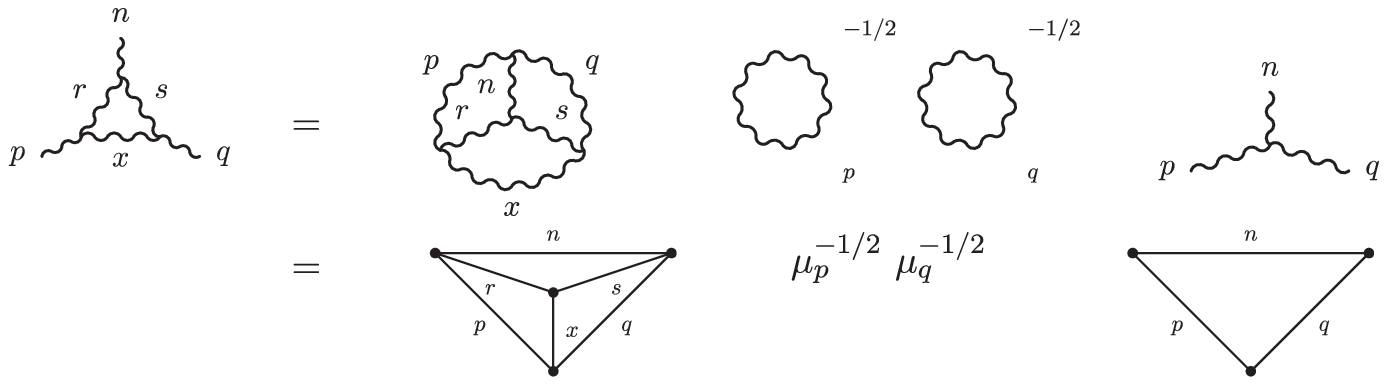}
\caption{The $TET \, ({}^{p \, n \, q}_{s \,  x \,  r} ) $ symbol appears as the f\/irst term on the right-hand side.}
\label{fig:vertexwithaloop}
\end{figure}

Since the triangle spaces are of dimension~$1$, the morphism def\/ined by the left-hand side of Fig.~\ref{fig:vertexwithaloop} should be proportional to the Y-shaped intertwiner $Y_{qpr}^\dag$.
This coef\/f\/icient, that appears on the r.h.s.\ of Fig.~\ref{fig:vertexwithaloop}, is written as the product of quantum dimensions and a new scalar quantity called the  tetrahedral symbol $TET$ symbol (the f\/irst term).
In the ``pure'' $SU(2)$ theory, quantum or not, the general expression of $TET$, which enjoys tetrahedral symmetry as a function of its six arguments, is known~\cite{Kauffman:book}.
The f\/irst three arguments of  $TET$ build an admissible triangle and the last three arguments determine skew lines (in $3$-space) with respect to the f\/irst three.
In the $SU(2)$ theory, this symbol vanishes if $m$, $n$, $p$ all refer to the fundamental representation $\sigma$ because the cube $\sigma^3$ does not contain the trivial representation. This is precisely what is dif\/ferent in the $SU(3)$ theory.

\begin{figure}[h]
\centering
\includegraphics{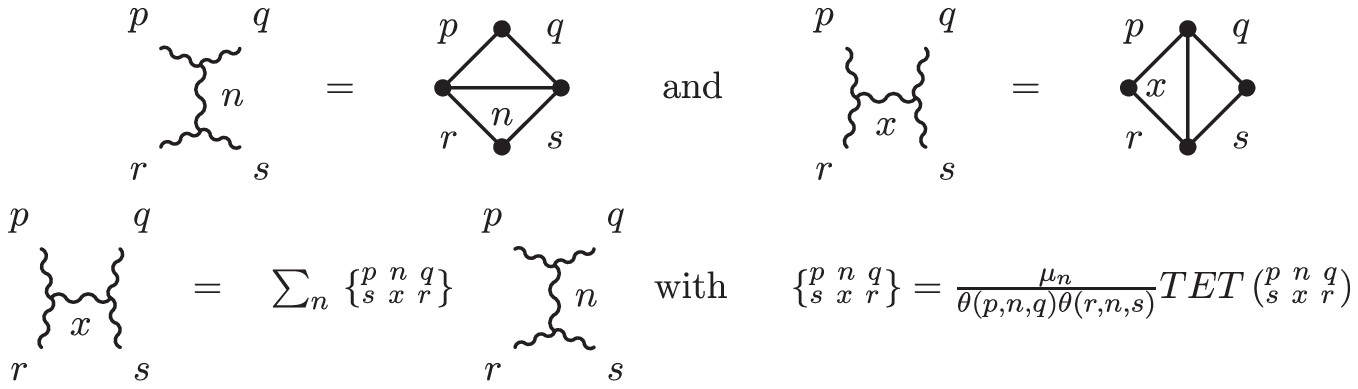}
\caption{Vertical and horizontal dif\/fusion graphs for the pure $SU(2)$ case.}
\label{fig:difgraphs}
\end{figure}

The (quantum) Racah and Wigner $6j$ symbols dif\/fer from the symbol $TET$ by normalizing factors. One should be warned that besides $TET$,  there are at least three types of quantities called~``$6j$'' symbols in the literature,  even in the pure $SU(2)$ case.
The Racah symbols directly enter the structure of the weak bialgebra $\mathcal B$ associated with the data~\cite{Ocneanu:paths}:
For each simple object $n$ one introduces a vector space $H_n$ spanned by a basis whose elements are labelled by admissible triangles\footnote{Or essential paths of length $n$ from $m$ to $p$.} (with f\/ixed edge $n$) representing the intertwiners~$Y_{mnp}$; one def\/ines $\mathcal B  = \bigoplus_n {\mathcal B}_n= \bigoplus_n   {\rm End}(H_n)$. The endomorphism product is depicted by vertical concatenation of vertical dif\/fusion graphs\footnote{Or double triangles sharing a horizontal edge, see Fig.~\ref{fig:difgraphs}.}:
$H_n ({}^{p q}_{r s}) H_{n^\prime} ({}^{r^\prime s^\prime}_{t u}) = \frac{\theta(r n s)}{\mu_n} \delta_{n n^\prime} \delta_{r r^\prime}  \delta_{s s^\prime} H_n ({}^{p q}_{t u})$, where the pre-factor is read from  Fig.~\ref{fig:verticallinewithloop}. These vertical dif\/fusion graphs represent elementary endomorphisms, up to scale (a pre-factor  ${\sqrt{\theta(p,n,q) \theta(r,n,s)}}/{\mu_n})$. Written in terms of matrices, $\mathcal B$ is a f\/inite sum of simple blocks labeled by the simple objects $n$.

\begin{figure}[h]
 \centering
\includegraphics{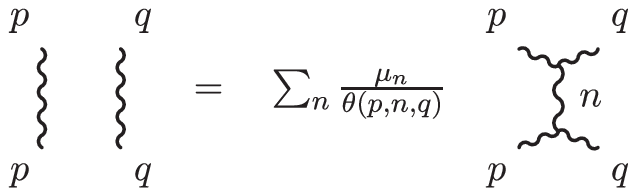}
\caption{Closure relation in the pure $SU(2)$ case.}
\label{fig:closureAA}
\end{figure}

\looseness=-1
One then introduces horizontal dif\/fusion graphs $V_x ({}^{p q}_{r s})$, (equivalently, double triangles sharing a horizontal edge), by the equation displayed on the second line of Fig.~\ref{fig:difgraphs}.
These new dif\/fusion graphs build a new basis for the same algebra  and def\/ine a new grading $\mathcal B = \bigoplus_x {\mathcal B}_x$.
Using this new basis one introduces a multiplication on the dual $\widehat {\mathcal B}$ by concatenation of the horizontal dif\/fusion graphs that can also be used to represent the corresponding dual basis:
$V_x ({}^{p q}_{r s}) V_{x^\prime} ({}^{q^\prime t}_{s^\prime u}) = \frac{\theta(q x s)}{\mu_x} \delta_{x x^\prime} \delta_{q q^\prime}  \delta_{s s^\prime} V_x ({}^{p t}_{r u})$.
Equivalently this def\/ines a coproduct in $\mathcal B$ that can be shown to be compatible with the f\/irst product. The obtained structure is actually a f\/inite dimensional weak Hopf algebra, a quantum groupoid.
The (Racah) symbols $\{ \; \}$ used to def\/ine the pairing do not enjoy tetrahedral symmetry, but only quadrilateral symmetry. If we set $x = 0$ in the duality relation (second line of Fig.~\ref{fig:difgraphs}), the tetrahedron degenerates to a triangle and $\TET {p}{n}{q}{q}{0}{p} = \theta(p,n,q)$ that cancels one of the denominators that enters the def\/inition of the Racah symbol, so that one obtains the closure relation Fig.~\ref{fig:closureAA} that will be used later.
The proportionality coef\/f\/icient entering this relation can be simply checked by writing $p \otimes q = \sum n$, taking the trace on both sides,  the fact that the quantum dimension is a homomorphism, and using the def\/inition of~$\theta$.

One often def\/ines (Wigner) $6j$ symbols by the equation
\[
\wigner {p}{n}{q}{s}{x}{r}  = {\TET  {p}{n}{q}{s}{x}{r} }/{\sqrt{\theta(p,n,q) \theta(p,x,r) \theta(r,n,s) \theta(s,x,q)}}.
 \]
 The normalizing factor in the denominator involves the product of theta functions for all triangles of the given tetrahedra, so that the $6j$ symbols $[\; ]$ also enjoy tetrahedral symmetry. The terminology is unfortunately not standard: The classical limit of the $[ \; ]$ symbols are called $6j$ by most physicists (and Racah $6j$ by Mathematica) but the $\{ \; \}$ are called $6j$ by~\cite{Kauffman:book} and \cite{Ocneanu:paths};  one can even f\/ind, for instance in \cite{Carter:book} and in several references studying the geometry of $3$-manifolds, a def\/inition of  ``unnormalized $6j$ symbols'' that dif\/fer from $TET$, 
  and $[\;]$ by other normalizing factors.
The reader may look at \cite{Coquereaux:6j} for relations between these quantities\footnote{Warning: In this reference the dif\/fusion graphs are rescaled, and  the loop has value $1$, not $[2]$.}  and their use in studying the quantum groupoid structure.

\subsubsection[Coupling of $SU(2)$ to $ADE$ matter]{Coupling of $\boldsymbol{SU(2)}$ to $\boldsymbol{ADE}$ matter}\label{section3.2.3}

\begin{figure}[h]
\centering
\includegraphics{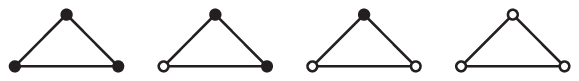}
\caption{The four types of triangles for $SU(2)$ coupled to $ADE$ matter.}
\label{fig:typesoftriangles}
\end{figure}

If we  now couple $SU(2)$ to $ADE$ matter\footnote{Actually to $DE$ matter since the $SU(2)$ at level $k$ is the $A$ case itself.}, the story becomes more complicated (see our general discussion),  since we have in general four types of triangles (see Fig.~\ref{fig:typesoftriangles}) and f\/ive types of genera\-li\-zed $TET$, Racah, or Wigner $6j$ symbols coupled by  pentagonal identities nicknamed ``the Big Pentagon'' equation (this was commented in \cite{Ocneanu:StFrancois},  see~\cite{BohmSzlachanyi}).

\begin{figure}[h]
\centering
 \includegraphics{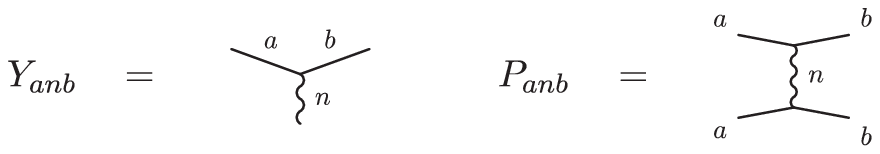}
\caption{The $Y$ and $P$ intertwiners for $SU(2)$ coupled to $ADE$ matter.}
\label{fig:YanbPanb}
\end{figure}

\begin{figure}[h]
\centering
\includegraphics{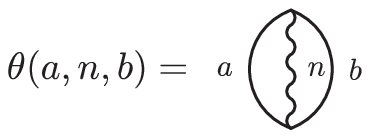}
\caption{The theta symbol for $SU(2)$ coupled to $ADE$ matter.}
\label{fig:thetaanb}
\end{figure}

Since ${\mathcal E}$ is a module-category over  ${\mathcal A}_k(SU(2))$ we obtain endofunctors of ${\mathcal E}$ labelled by the simple objects $\sigma_n$ of  ${\mathcal A}$, sending  objects of ${\mathcal E}$ to objects of ${\mathcal E}$. Again, one represents diagrammatically such morphisms from $a$ to $b$  by $Y$-shaped diagrams called $Y_{anb}$ (see the f\/irst  diagram on Fig.~\ref{fig:YanbPanb}) where $a$ and $b$ refer to simple objects of ${\mathcal E}$.
By composition with the upside-down diagrams, we get the second diagram of Fig.~\ref{fig:YanbPanb} that can be traced, as in the pure $SU(2)$ case,  to get  numbers $\theta_{anb}$, see  Fig.~\ref{fig:thetaanb} that replaces Fig.~\ref{fig:thetamnp}.

\begin{figure}[h!]
\centering
\includegraphics{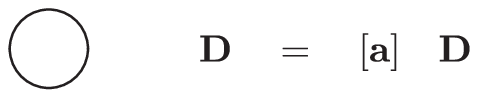}
\caption{The circle for an irreducible object $a$ of type  $\mathcal{E}$.}
\label{fig:SU2circleE}
\end{figure}

\begin{figure}[h!]
\centering
\includegraphics{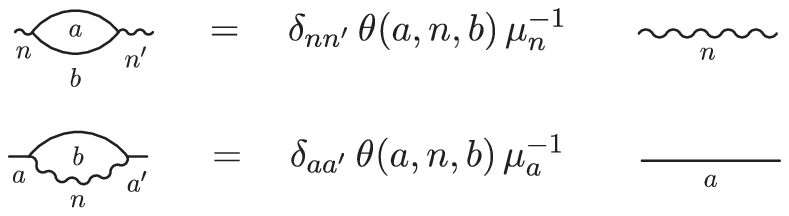}
\caption{Normalizing coef\/f\/icients for $SU(2)$ coupled to $ADE$ matter.}
\label{fig:Elinewithloop}
\end{figure}

\noindent
 In particular, when $n$ is trivial, the value of the circle associated with an irreducible object $a$ of type $\mathcal{E}$ is the corresponding quantum number\footnote{In general we call $ \mu_a = [a]$ the $q$-dimension of $a$, but for $SU(2)$ and $a = \sigma_n$ we have $ \mu_a= [\sigma_n]=[n+1]$.} $ \mu_a$ (see Fig.~\ref{fig:SU2circleE}).
Composing the morphisms in other directions gives the equations depicted in Fig.~\ref{fig:Elinewithloop}.

\begin{figure}[h]
\centering
 \includegraphics{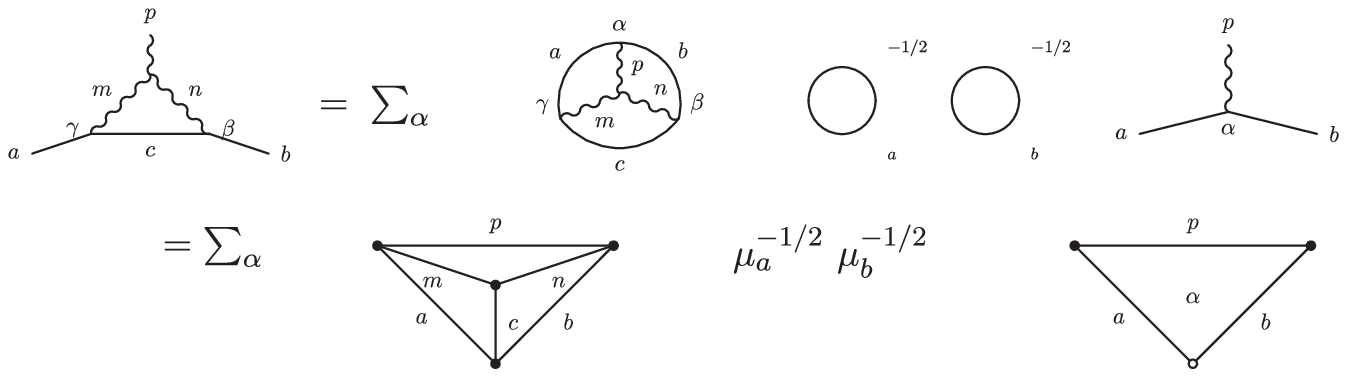}
\caption{The $TET \, ({}^{a \, p \, b}_{n\,  c\,  m} ) $ symbol on the r.h.s.\ may depend on $\alpha$, $\beta$, $\gamma$.}
\label{fig:AEEvertexwithaloop}
\end{figure}

One has to be cautious here because the dimensionality of triangle spaces like $anb$ can be bigger\footnote{The matrix coef\/f\/icients of the annular matrices obtained from the relation $F_n=F_{n-1} F_1 - F_{n-2}$, where $F_1$ is the adjacency matrix of the fusion graph of ${\mathcal E}$, can be bigger than $1$.} than $1$;  in other words the space of essential paths from $a$ to $b$ and of length $n$ may have a~dimension bigger than $1$.  One should therefore sometimes introduce an extra label~$\alpha$ for vertices such as~$anb$.
 Comparing the $TET$ coef\/f\/icients introduced in Figs.~\ref{fig:vertexwithaloop} and \ref{fig:AEEvertexwithaloop} one notices that~$m$,~$n$,~$p$ still label simple objects of ${\mathcal A}_k(SU(2))$, but  $a$, $b$, $c$, replacing $p$, $q$, $x$, refer to simple objects of~${\mathcal E}$.
The other dif\/ference is that the right hand side of Fig.~\ref{fig:AEEvertexwithaloop} now incorporates a~summation over the extra label $\alpha$, in those cases where the triangle space $a$, $p$, $b$ has a dimension bigger than~$1$.

\begin{figure}[h]
 \centering
 \includegraphics{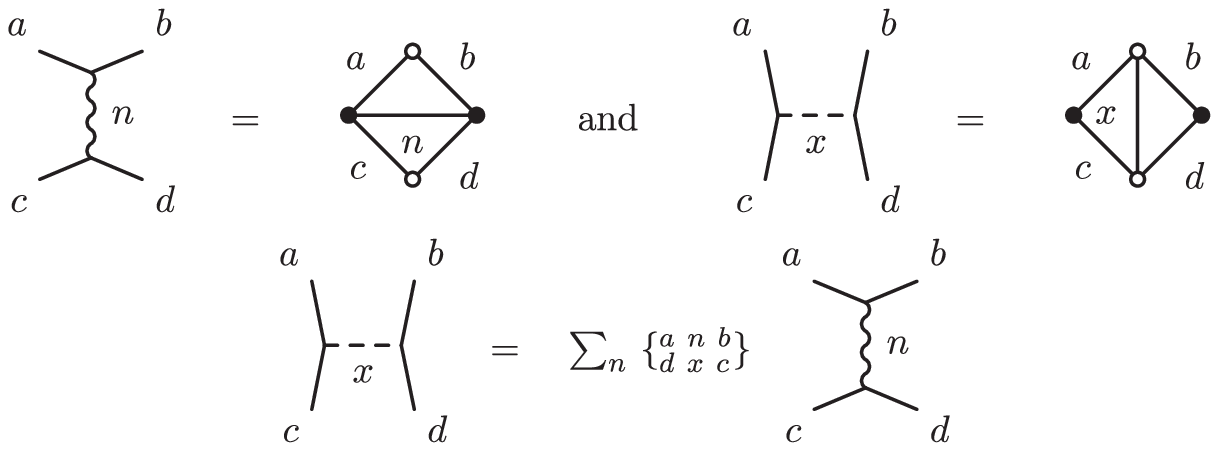}
\caption{Vertical and horizontal dif\/fusion graphs for $SU(2)$ coupled to $ADE$ matter.}
\label{fig:difgraphsAE}
\end{figure}

\begin{figure}[h]
\centering
\includegraphics{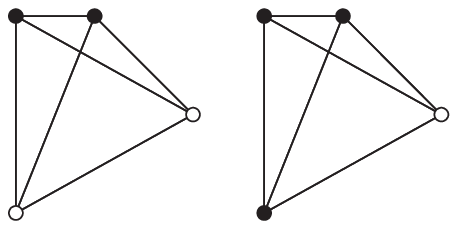}
\caption{Two types of tetrahedra: $AAEE$ and $AAAE$.}
\label{fig:twotypesoftetrahedra}
\end{figure}

\begin{figure}[h]
 \centering
\includegraphics{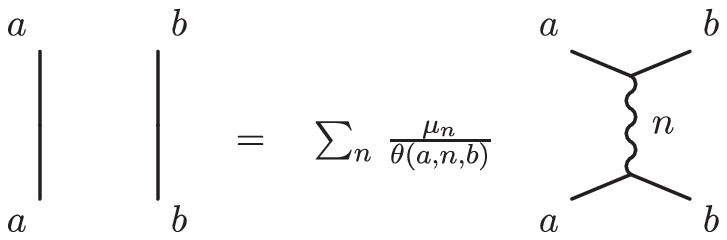}
\caption{Closure relation for $SU(2)$ coupled to matter.}
\label{fig:closureAE}
\end{figure}

Like in the pure $SU(2)$ case one can introduce a vector space $H = \bigoplus_n H_n$  spanned by a basis whose elements are labelled by admissible triangles (with f\/ixed edges $n$ standing for simple objects of ${\mathcal A}_k(SU(2)$) representing intertwiners $Y_{anb}$. These elements can be interpreted in terms of essential paths from $a$ to $b$, of length $n$, on the fusion graph of ${\mathcal E}$ (see the discussion in Section~\ref{section3.2.4}).
In this $SU(2)$ situation,  it is known that $H$ can also be constructed as the vector space underlying the Gelfand--Ponomarev preprojective algebra \cite{GP} associated with the corresponding $ADE$ unoriented quiver.
With $k \neq \infty$, identif\/ication stems from the fact that dimensions of the f\/inite dimensional  vector spaces $H_n$, calculated according to both def\/initions, are equal (compare for instance~\cite{RobertGilSL3Categories}  and~\cite{MOV}).
In the case of $SU(3)$ considered in the next section, the grading label of the horizontal  space $H_n$ refers to a pair of integers $n=(n_1,n_2)$, seen as a Young tableau, and can no longer be interpreted as a length.
When ${\mathcal E} \neq {\mathcal A}_k$ the ring~${\mathcal O({\mathcal E})}$ is usually not isomorphic with  ${\mathcal A}_k$ but the construction of the bialgebra $\mathcal B$ proceeds like in the pure $SU(2)$ case, by using vertical and horizontal dif\/fusion graphs together with the $6j$ symbols appearing in Fig.~\ref{fig:difgraphsAE}.
We shall not discuss here the structure of $\mathcal B$, the interested reader may look at \cite{Ocneanu:paths, BohmSzlachanyi, PetkovaZuber:Oc, Coquereaux:6j, TrincheroDTA},
but notice that these $6j$ symbols have four labels for the simple objects of type $\mathcal E$ (the module), one label for the type $\mathcal A$ and one label for the type $\mathcal O$ (the quantum symmetries), whereas the $6j$ symbols appearing in Fig.~\ref{fig:AEEvertexwithaloop} have three indices  of type $\mathcal E$ and three indices of type $\mathcal A$. Using star-triangle duality,  vertices become edges, $AA$, $A\mathcal{E}$ or $\mathcal{E}\mathcal{E}$, and the two kinds of $6j$ symbols just considered belong to the tetrahedral types described in Figs.~\ref{fig:twotypesoftetrahedra},
where black and white vertices refer to $A$ and $\mathcal{E}$ marks.
If we set $x=0$ in the duality relation (second line of Fig.~\ref{fig:difgraphsAE}), the tetrahedron degenerates and becomes a triangle,  like in the pure $SU(2)$ situation, and the Racah symbol $\racah {a}{n}{b}{d}{x}{c} =  \frac {\mu_n}{\theta(a,n,b) \theta(c,n,d)} \TET  {a}{n}{b}{d}{x}{c}$ becomes simply $\racah {a}{n}{b}{b}{0}{a} =  \frac {\mu_n}{\theta(a,n,b)}$.
We shall need later the $SU(3)$ analog of this closure relation, which is displayed on Fig.~\ref{fig:closureAE}.

\subsubsection[$SU(2)$ braiding, Hecke and Jones-Temperley-Lieb algebras]{$\boldsymbol{SU(2)}$ braiding, Hecke and Jones--Temperley--Lieb algebras}\label{section3.2.4}

\begin{figure}[h!]
\centering
\includegraphics{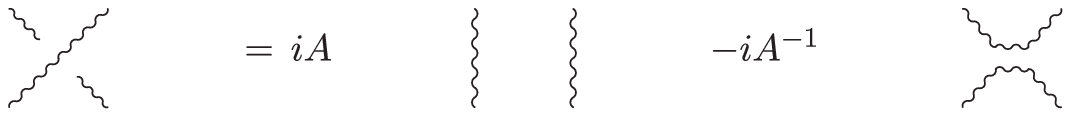}
\caption{The braid relation for $SU(2)$.}
\label{fig:braidingSU2}
\end{figure}

In the $SU(2)$ theory, one also chooses a square root $A = \exp(\frac{i \pi}{2 \kappa})$ of $q$ and introduces a braiding\footnote{The factor $i$, absent in \cite{Kauffman}, insures that the classical situation is obtained by letting $A \mapsto 1$.} def\/ined by  Fig.~\ref{fig:braidingSU2}.
It reads $X = i A  \one - i A^{-1}  U = \overline {\epsilon}  \one + \epsilon  U$ where $\overline{\epsilon} = i A$, so that $\epsilon \overline{\epsilon} =1$, and $U$ is the already def\/ined cup-cap generator.
Notice that  $\overline {\epsilon} X = -1/q + U$. In physics, the braiding can be interpreted as a Boltzman weight at criticality in the context of RSOS height models of type $SU(2)$.
The algebra of the braid group with $s+1$ strands is def\/ined by generators
$g_1,\ldots, g_s$ obeying relations $g_n g_{n+1} g_n = g_{n+1} g_n g_{n+1}$ and $g_n g_m = g_m g_n $ when $\vert m-n  \vert \geq 2$. The Hecke algebra with parameter $\beta = [2]= q+1/q = 2 \cos(\pi/\kappa)$ is obtained as its quotient by the ideal generated by the quadratic relation $(g_n+1)(g_n-q^2)=0$.  Introducing the change of generators def\/ined by\footnote{Another favorite possibility is to take the r.h.s.\ equal to $q U_n - \one$ which amounts to replace the projector $U_n/\beta$ by the projector  $\one - U_n/\beta$.}  $g_n = q(q-U_n)$,  the above three relations def\/ining the Hecke algebra can be written
\begin{gather*}
U_n U_m = U_m U_n \qquad {\rm when} \quad \vert m-n  \vert \geq 2, \\
U_n^2 = \beta  U_n, \qquad
U_n U_{n+1} U_n -  U_n = U_{n+1} U_{n} U_{n+1} -  U_{n+1}.
\end{gather*}
In the $SU(2)$ theory, the operator $U$ is realized as the cup-cap operator appearing on the r.h.s.\ of Fig.~\ref{fig:braidingSU2}. A simple graphical calculation shows that $U_n$ actually obeys  the  stronger equality $U_n U_{n\pm 1} U_n =  U_n$, that replaces the third relation and also ensures the vanishing of a~quantum Young projector def\/ined in the algebra of the braid group with three strands (see for instance~\cite{DF-meanders}).
The obtained algebra is the Jones--Temperley--Lieb algebra $C_{\beta}^s(SU(2))$.

In order to better understand what is modif\/ied when we trade $SU(2)$ for $SU(3)$, we now remind the reader how to build  representations of  $C_\beta(SU(2))$ on vector spaces obtained from paths def\/ined on the fusion graph of any module-category ${\mathcal E}$ associated with ${\mathcal A}_{k}(SU(2))$, i.e.\ on $ADE$ diagrams, or on the fusion graphs of ${\mathcal A}_{\infty}(SU(2))$, i.e.\ $SU(2)$ itself and its subgroups (which are af\/f\/ine $ADE$ diagrams because of the McKay correspondence).

An elementary path is a f\/inite
sequence of consecutive  edges on a fusion graph of the chosen module ${\mathcal E}$:
$\xi(1) = \xi_{a_{1}a_{2}} $,
$\xi(2) = \xi_{a_{2}a_{3}} $, etc.
Vertices are paths of length $0$.
The length of the (possibly backtracking) path $( \xi(1)\xi(2)\cdots
\xi(p) )$ is $p$.
We call $r(\xi_{cd})={d}$,
and $s(\xi_{cd})={c}$, the range and source of  $\xi_{cd}$.
For all edges $\xi(n+1) = \xi_{cd}$ that appear in an elementary path,
we set  ${\xi(n+1)}^{-1} \doteq \xi_{dc}$, and there is no ambiguity since, in the $SU(2)$ theory, all fusion graphs are simply laced.
We call ${\rm Path}^p_{a,b}$ the vector space spanned by all elementary paths of length $p$ from $a$ to $b$.
For every integer $n >0$, the annihilation operator $C_{n}$,
acting on elementary paths of length $p$ is def\/ined
as follows: If $p \leq n$, $C_{n} = 0$, whereas if $ p \geq  n+1$ then:
\begin{gather*}
C_{n} (\xi(1)\cdots\xi(n)\xi(n+1)\cdots \xi(p)) =
\sqrt\frac{ \mu_{r(\xi(n))}}{ \mu_{s(\xi(n))}}
\delta_{\xi(n),{\xi(n+1)}^{-1}}
 (\xi(1)\cdots{\hat\xi(n)}{\hat\xi(n+1)}\cdots  \xi(p)).
\end{gather*}
The symbol ``hat'', like  in $\hat \xi$, denotes omission.
The result is therefore either $0$ or a linear combination of elementary paths of length $p-2$.
 This def\/inition is extended by linearity to arbitrary elements of  ${\rm Path}^p_{a,b}$, called ``paths''.
Intuitively, $C_{n}$ chops the round trip that possibly appears
at position~$n$.
A path, not elementary in general, is called essential if it belongs to
the intersection of the kernels
of the annihilators $C_{n}$'s (or of Jones projectors, see below).
One introduces a scalar product in the vector space ${\rm Path}^p_{a,b}$ by declaring that elementary paths are orthonormal.
Acting on elementary paths of length $p$, the creating operators $C^{\dag}_{n}$ act as follows.
If $n > p+1$, $C^{\dag}_{n}=0$, whereas if $n \leq p+1$ then,
setting $c = r(\xi(n-1))$,
$
C^{\dag}_{n} (\xi(1)\cdots\xi(n-1)\cdots) = \sum_{\vert d-c\vert=1}
\sqrt {\frac{ \mu_{d}}{ \mu_{c}}}  (\xi(1)\cdots\xi(n-1)\xi_{cd}\xi_{dc}\cdots).
$
The previous sum is taken over the neighbors ${d}$ of ${c}$ on the fusion graph.
Intuitively, this operator inserts all possible small round trip(s)
in position~$n$.
The result is therefore either $0$ or a linear combination of paths of
length $p+2$. In particular, on paths of length zero (i.e.\ vertices), one obtains:
$
C^{\dag}_{1} ({a}) = \sum_{\vert b-a\vert=1}
\sqrt{\frac{ \mu_{b}}{ \mu_{a}}}  \xi_{ab}\xi_{ba}.$
The Jones' projectors $e_{n}$ are obtained as endomorphisms of
${\rm Path}^p$ by
$
e_{n} \doteq \frac{1}{\beta} C^{\dag}_{n} C_{n}.
$
All expected relations between the $e_n$, or between the $U_n = \beta   e_n$, are satisf\/ied, and we obtain in this way a path realization of the Jones--Temperley--Lieb algebra, hence a representation of the Hecke algebra with parameter $\beta$.
In particular,  for every pair $a$, $b$ of neighboring vertices:
\[
U(\xi_{ab} \xi_{ba}) = \sum_{\{c : \; \vert c-a \vert = 1\}}     \frac{ \sqrt{ \mu_b    \mu_c}} { \mu_a}
 \, \xi_{ac} \xi_{ca}
\]
so that for every vertex $a$ of the fusion graph of ${\mathcal E}$, i.e.\ a~simple object of the later, one obtains a~representation on the space ${\rm Path}_{a,a}^2$ spanned by elementary paths of length $2$ from $a$ to $a$ (round trips).  In the $SU(2)$ theory there are only single edges, they are un-oriented, and vertices have $1$, $2$ or $3$ neighbors. The matrix representative of $U$ is  a direct sum of blocks of sizes $(1,1)$, $(2,2)$ or $(3,3)$.
Once $U$ is represented as above, the braiding $X$ can be considered as one of the two fundamental connections (say $1_L$) attached to the cell system def\/ined by ${\mathcal E}$.  Its value on a~basic\footnote{The edges of the quadrilateral are of length~$1$, so that consecutive vertices are neighbors on the fusion graph.} cell is:
\[
\begin{array}{@{}c}
\setlength{\unitlength}{.1in}
\sqedges {}{}{}{}
\sqface {a}{b}{d}{c}{X}{0}{.5}
\punit1
=
\punit1
\end{array}
 \overline {\epsilon}   \delta_{bc} + \epsilon   \delta_{ad}    R
\quad
\text{with}
\quad
R = \sqrt{\frac{ \mu_b  \mu_c}{ \mu_a  \mu_d}}.
\]

The other fundamental connection ($1_{R}$) is obtained by replacing  $\epsilon$ by its conjugate $\overline \epsilon$ in the above  expression. Biunitarity of the connection is ensured by the relation $\epsilon^2 + \overline{\epsilon}^2 + \beta = 0$.
The quantity $R$ describes a kind of generalized parallel transport from top to bottom horizontal edges  along the vertical edges.
 The above basic cell is  also a fundamental Racah symbol  $\racah {a}{n}{b}{d}{x}{c} $, for $n=1$ and $x=1_L$. It gives the pairing of matrix units belonging to the block of the bialgebra ${\mathcal B}$ labeled by $n=1$, the fundamental representation of ${\mathcal A}_k(SU(2))$, and the block of $\widehat{\mathcal B}$ labelled by $x=1_L$, one of the two chiral generators of ${\mathcal O} = {\rm End}_{{\mathcal A}_k}{\mathcal E} $. It can be pictured as a tetrahedron pairing a double triangle $(anb)(cnd)$, i.e.\ a~matrix unit of ${\mathcal B}$,  and a double triangle $(axc)(bxd)$, i.e.\ a matrix unit of $\widehat{{\mathcal B}}$.

\subsubsection[Comment about the $SU(2)$ epsilon tensor]{Comment about the $\boldsymbol{SU(2)}$ epsilon tensor}\label{section3.2.5}

Before embarking on $SU(3)$ let us return  to the determinant map (or epsilon tensor  $\epsilon$) of~$SU(2)$.
Call $\{\alpha_1, \alpha_2\}$ an orthonormal basis of the representation space  associated with the fundamental~$\sigma$.
The intertwiner $\sigma \otimes \sigma \mapsto \CC$ is the epsilon tensor $\epsilon$, def\/ined by $\epsilon(\alpha_i \otimes \alpha_j) = \epsilon_{ij}$, the signature of the permutation~$(ij)$, which is $0$ if $i= j$ and $\pm 1$ otherwise.
Take $v,w$ two vectors in $\sigma$, then $\epsilon(v\otimes w) = \epsilon_{ij} v^i w^j = v^1w^2 - v^2 w^1$, which is  the determinant $\det(v,w)$.  One often says that this map $\sigma \otimes \sigma \mapsto \CC$  induces, by duality, an isomorphism from $\sigma$ to $\overline \sigma$ (a particular case of Frobenius isomorphism), but there is a subtlety here:  since $\det$ is antisymmetric under a circular permutation of the columns (a transposition in this case) one does not obtain a~single isomorphism from $\sigma$ to $\overline \sigma$, but two, with zero sum, namely $(\alpha_1, \alpha_2) \mapsto ( - \overline \alpha_2,\overline \alpha_1) $ and  $(\alpha_1, \alpha_2) \mapsto ( \overline \alpha_2,- \overline \alpha_1)$. In other words, in any wire model of $SU(2)$ the wires will be unoriented (since~$\sigma$ is equivalent to its conjugate), but one should, in principle, use two distinct wires joining boundary points because there are two fundamental intertwiners. Since they have zero sum, it does not harm to make a~choice once and for all, so that in practice only one kind of wire is used. The discussion is not modif\/ied when we move from  $SU(2)$ to $SU(2)_q$.
 The only place where this subtlety should be remembered is when one tries to def\/ine an action of the category ${\mathcal A}_k(SU(2))$ on a module whose fusion graph would contain tadpoles:
the presence of a tadpole~-- loop~-- at the vertex~$a$ (a~simple object) would mean that, besides the identity, we could have another intertwiner from $a$ to $a$, but they should be proportional since~$a$ is simple, however the sign of this proportionality factor $\theta(a,1,a)$ will change if we select the opposite convention for the isomorphism between~$\sigma$ and $\overline \sigma$.   The only possible conclusion is that the proportionality factor is zero, which means that such a graph cannot be associated with any $SU(2)$ module. This argument was presented in~\cite{Ocneanu:Bariloche}. Another proof eliminating tadpole graphs was given in~\cite{Ostrik:classificationSU2}.

\subsection[Equations for the $SU(3)$ model]{Equations for the $\boldsymbol{SU(3)}$ model}\label{section3.3}

The $SU(3)$ case shares many features with the $SU(2)$ case. This was one of the reasons to discuss the simpler later case in the previous section, so that we can now concentrate on novel features.
Besides the fact that one has to take into account orientation of edges, complex conjugation (non trivial Frobenius isomorphisms), and the fact that  have usually a dimension bigger than~$1$, the general discussion about the structure of the weak bialgebra $\mathcal B$ and the dif\/ferent types of~$6j$ symbols stays essentially the same. In particular  there are again f\/ive types of tetrahedra and therefore f\/ive types of cells\footnote{$6j$ symbols for $G$ at level $k$ coupled to a module $\mathcal E$ are often called ``Ocneanu cells''.}, but those that we need to consider here are of type  $AAA\mathcal{E}$.  As already stated, $6j$ of that type vanish for $SU(2)$ when the three edges $AA$ (simple objects of $\mathcal A$) all coincide with the fundamental representation~$\sigma$, this is precisely what is not true for~$SU(3)$.

\subsubsection{The elementary version} \label{sec:su3graphicalelements}\label{section3.3.1}

\begin{figure}[h]
\centering
\includegraphics{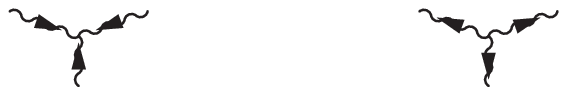}
\caption{The triple vertex for $SU(3)$.}
\label{triplevertex}
\end{figure}

\looseness=-1
For $SU(3)$, we call  $\{\alpha_1, \alpha_2, \alpha_3\}$ the orthonormal basis of $\sigma$. Again,
the  intertwiner (determinant map) $\det : \sigma \otimes \sigma \otimes \sigma \mapsto \CC$ is def\/ined by the epsilon tensor $\epsilon(\alpha_i \otimes \alpha_j \otimes \alpha_k) = \epsilon_{ijk} \alpha_k$, with~$\epsilon_{ijk}$ the signature of the permutation~$(ijk)$.
Using Frobenius isomorphisms on the identity intertwiner~$\sigma \mapsto \sigma$
gives an intertwiner $\sigma \otimes \overline \sigma \mapsto \CC$ represented as a cup diagram.
Here~$\sigma$ is not equivalent to~$\overline \sigma$ and the basic wire (representing an intertwiner from $\sigma$ to itself, or from $\overline \sigma$ to itself) is oriented: a downward oriented wire, read from top to bottom, denotes the identity morphism from~$\sigma$ to~$\sigma$ and an upward oriented wire, also read from top to bottom, denotes the identity morphism from~$\overline \sigma$ to~$\overline  \sigma$. The determinant map is symmetric under circular permutation of the columns and the problem discussed for $SU(2)$ at the end of the previous section does not arise. In particular, tadpoles in graphs describing $SU(3)$ modules may exist.
Using Frobenius isomorphisms on the determinant map $\sigma^3 \mapsto \CC$ leads to an intertwining operator $\sigma \otimes \sigma \mapsto \overline \sigma$ graphically represented\footnote{$\{e_i\}$ being a basis of $\CC^3$, this is the map $e_i \otimes e_j \mapsto \epsilon_{ijk} {\overline e_k}$.} by the  triple vertex (Fig.~\ref{triplevertex}) with all wires oriented in, and an  intertwining operator $\overline \sigma \otimes \overline \sigma \mapsto  \sigma$
represented by the  triple vertex (Fig.~\ref{triplevertex}) with all wires oriented out.

A graphical TQFT model for $SU(3)$ is therefore given by an orientable surface with a f\/inite set of distinguished boundary points, and every such point should have a direction (inner or outer) associated with it. A diagram in the model (primitive version) is a wire diagram made of oriented wires that meet at triple points, or end at boundary points, in a way compatible with the predef\/ined exits  (directions associated with boundary points).  Notice that triple vertices are not boundary points.
Existence of the determinant map, an intertwiner from $\sigma^3$  to $\CC$, the trivial representation,  can be interpreted in terms or the well known representation theory of~$SU(3)$, or of  quantum $SU(3)_q$.
Indeed  $\sigma^3$, which is not irreducible, decomposes as follows into a direct sum (we use highest weights notations for  irreducible representations):
\[
((1,0) \otimes (1,0)) \otimes (1,0) =( (0,1) \oplus (2,0)) \otimes (1,0) =((0,0) \oplus (1,1)) \oplus ((1,1) \oplus (3,0)).
\]
 In terms of classical notations: $(3)^3=(\overline{3} + 6)\,  3 = (\overline{3}\,   3)  +  (6 \, 3) =  (1 + 8) + (8 +  10)$, and  in terms of quantum dimensions, if the level $k$ is large enough,  $[3]^3 = (\overline{[3]} + [3] [4]/[2])\times [3] = (\overline{[3]} \times  [3])  +  ( [3] [4]/[2]   \times[3]) =  ([1] + [2][4]) + ([2][4] +  [4][5]/[2]).$ These elementary calculations show that indeed, the trivial representation appears on the r.h.s.

\begin{figure}[h]
\centering
\includegraphics{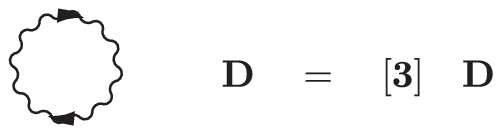}
\caption{The circle for $SU(3)$.}
\label{fig:SU3circle}
\end{figure}

We now turn to ``relations between relations'' (i.e.\ to relations between intertwiners),  there are three of them. The f\/irst, reminiscent of the $SU(2)$ theory, says that the circle has value $[3]$. Indeed, in the $SU(3)$ theory, we have the same cup-cap morphisms as before, but they are now  oriented. The cap graph, read from top to bottom, says that there is an intertwiner from the trivial representation to  $3 \otimes  \overline 3$, and the cup graph, read from top to bottom, says that  $3 \otimes \overline 3$ contains the trivial representation. Composition of the two, from top to bottom, gives a number, i.e.\  the rule displayed on Fig.~\ref{fig:SU3circle}.  In the $SU(2)$ theory there are no further relations, but in the $SU(3)$ theory, there are two.

\begin{figure}[h]
\centering
\includegraphics{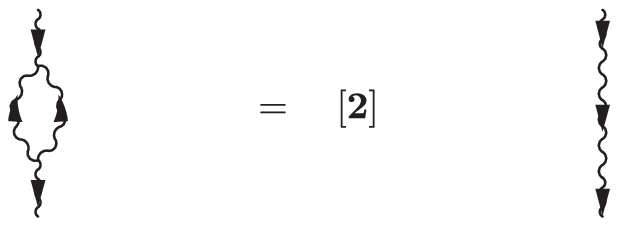}
\caption{The quadratic relation for $SU(3)$.}
\label{fig:quadraticrelation}
\end{figure}

{\it First non-trivial relation}. Consider the morphism $\sigma \mapsto \overline \sigma \otimes \overline \sigma \mapsto \sigma$.  Read  from top to bottom it is represented by the left hand side of Fig.~\ref{fig:quadraticrelation}.
  This morphism goes from $\sigma $ to $\sigma $, however $\sigma $ is irreducible, so this morphism is proportional to the  identity morphism. The equality is depicted by Fig.~\ref{fig:quadraticrelation}.
   Conventionally  the proportionality constant is equal to $\beta = q+1/q = [2]$.
 Tracing this equality gives immediately the value of the theta symbol:  $\theta(\sigma, \sigma, \sigma) = [2] \, [3]$.

\begin{figure}[h!]
\centering
\includegraphics{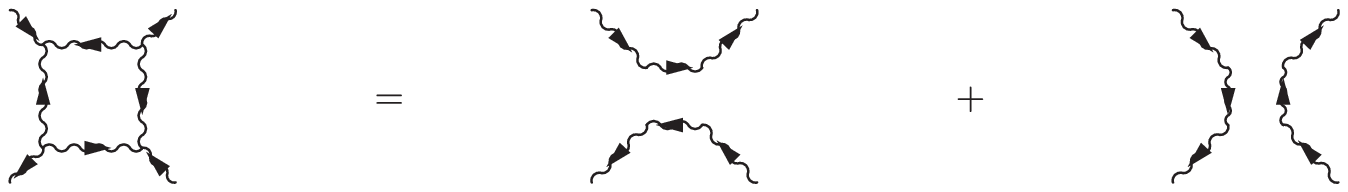}
\caption{The quartic relation for $SU(3)$.}
\label{fig:cuarticrelation}
\end{figure}

 {\it Second non-trivial relation}. Let us consider the possible morphisms from  $\sigma \otimes \overline \sigma$ to $\sigma \otimes \overline \sigma$. There are three of them (we suppress the tensor product signs in this discussion):
 \begin{enumerate}\itemsep=0pt
 \item  We do nothing, i.e.\ $\sigma$ goes to $\sigma$ and $ \overline \sigma$ goes to  $\overline \sigma$. It is  an identity morphism depicted by the last term on the r.h.s.\ of Fig.~\ref{fig:cuarticrelation}.
 \item  We compose the cup and cap morphisms: $\sigma \otimes \overline \sigma \to \CC \to  \sigma \otimes \overline \sigma$. This is depicted by the f\/irst term on the r.h.s.\ of  Fig.~\ref{fig:cuarticrelation}.
 \item The last possibility is to consider  the chain depicted by the l.h.s.\ of Fig.~\ref{fig:cuarticrelation}:
\[
\sigma \otimes \overline \sigma \to ( \overline \sigma  \otimes \overline \sigma) \otimes \overline \sigma =  \overline \sigma  \otimes ( \overline \sigma \otimes \overline \sigma )  \to    \overline \sigma   \otimes \sigma \to  (\sigma  \otimes  \sigma)  \otimes \sigma =  \sigma  \otimes  (\sigma  \otimes \sigma) =  \sigma  \otimes \overline \sigma.
\]
\end{enumerate}

 However, as we know, $\sigma \otimes \overline \sigma = \sigma_{(1,0)} \otimes \sigma_{(0,1)} = \sigma_{(0,0)} \oplus \sigma_{(1,1)}$, in terms of highest weights labels (or, in terms of quantum dimensions:  $[3] \times [3] = [1] + [2][4]$, or, classically, $3 \times \overline 3 = 1 + 8$).
So, from  $\sigma \otimes \overline \sigma$ we have two intertwiners, respectively to $\CC$ and to $\sigma_{(1,1)}$.  By composing them with their adjoints, i.e.\ reversing the arrows, we have two morphisms, and only two,  up to linear sums,  from $\sigma \otimes \overline \sigma$ to itself. This means that the {\it three} morphisms found previously cannot be independent: there is a linear relation between them. The third (for instance) can be written in terms of the f\/irst two.
Actually, as displayed on Fig.~\ref{fig:cuarticrelation}, the third is just {\it the sum} of the f\/irst two, this was shown by \cite{Kuperberg}  where it is also shown that all other relations between morphisms of $SU(3)$ are consequences of Figs.~\ref{fig:quadraticrelation} and~\ref{fig:cuarticrelation}.
Although belonging to the folklore (we already mentioned the fact that they were explicitly used by \cite{EwenOgievetsky}) these two relations are often called ``quadratic and quartic equations for the Kuperberg $SU(3)$ spider''.
The f\/irst is called  quadratic because it involves {\it two} triple points and the other is quartic  because it involves four triple points.
Since these relations hold for $SU(3)$, classical or quantum, they should also hold for any (classical or quantum) module, i.e.\ they should be true for $SU(3)$ itself, its f\/inite subgroups, for the fusion categories ${\mathcal A}_k(SU(3))$ and their module~-- categories $\mathcal E$.

\subsubsection{The clasp version (bound states version)}\label{section3.3.2}

There is of course a clasp version of the wire models, with new types of wires for every irreducible representation existing at  level $k$.  As usual, those ``bound states wires'' are labelled by highest weights $(p,q)$ or  Young tableaux. By looking at the expressions of fusion matrices for $SU(3)$ at level $k$,  one notices that the dimension of triangle spaces can be bigger than $1$,  even in the pure $SU(3)$ case, in contrast with the $SU(2)$ theory.
The endomorphism algebra of the object $\sigma^m \otimes \overline \sigma^n$, i.e.\ ${\rm Hom}_{SU(3)}(\sigma^m \otimes \overline \sigma^n,\sigma^m \otimes \overline \sigma^n)$ is the $SU(3)$ version of the Temperley--Lieb algebra.
It can be  generated by  linear combination of diagrams generalizing those discussed in the $SU(2)$ section; one can also def\/ine  generalizations $P_{mn}$ of the Wenzl projectors, with image the irreducible representation of $SU(3)$ indexed by the highest weight $(m,n)$, in the endomorphism algebra of $\sigma^m \otimes \overline \sigma^n$ (see \cite{Ohtsuki-Yamada,Kuperberg, Suciu}),  obeying a three terms recurrence relation. One can even def\/ine a~natural basis in arbitrary intertwiner spaces, see~\cite{Suciu}, however we shall not need these explicit constructions in full generality. The def\/inition of theta symbols and $6j$ symbols then follows the general rules described in the $SU(2)$ section, but this has to be done with much more care because edges are oriented, and also because the dimension $d_n$ of triangle spaces $H_n$, where $n=\{n_1,n_2\}$, is usually bigger than $1$, so that the notations, for example, have to incorporates new labels.
The recurrence relations for $SU(3)$ fusion matrices $N_n = (N_{n,p}^q)  = ({N_{(n_1,n_2),(p_1,p_2)}}^{(q_1,q_2)})$,  using the seed $N_{(1,0)}$, the adjacency matrix of the fusion graph, are:
\begin{gather*}
N_{(\lambda,\mu)}  =  N_{(1,0)}   N_{(\lambda-1,\mu)} - N_{(\lambda-1,\mu-1)} -
N_{(\lambda-2,\mu+1)} \qquad   \textrm{if}\quad  \mu \not= 0, \nonumber \\
N_{(\lambda,0)}  =  N_{(1,0)}   N_{(\lambda-1,0)} - N_{(\lambda-2,1)},\qquad
N_{(0,\lambda)}  =  (N_{(\lambda,0)})^{\rm  tr}. \nonumber
\end{gather*}
As usual,  $d_n = \sum_{p,q}   N_{n,p}^q $.

\subsection[Coupling to $SU(3)$ matter: Ocneanu coherence equations for $SU(3)$ triangular cells (self-connections)]{Coupling to $\boldsymbol{SU(3)}$ matter: Ocneanu coherence equations\\ for $\boldsymbol{SU(3)}$ triangular cells (self-connections)}
\label{section3.4}

\begin{figure}[h]
\centering
 \includegraphics{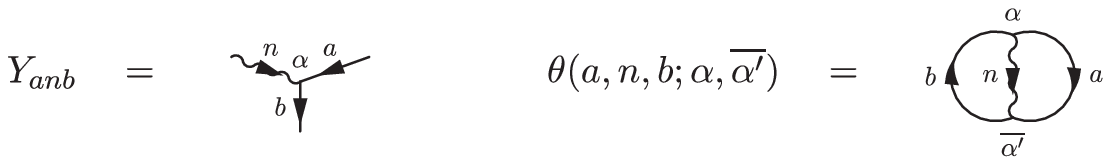}
\caption{The $Y$ intertwiner for $SU(3)$ coupled to  $SU(3)$  matter and its theta symbol.}
\label{fig:thetaanbSU3}
\end{figure}

Before deriving these equations we need to consider an extension of the $SU(3)$ recoupling system where ${\mathcal A}_k(SU(3))$ is allowed to act on a module ${\mathcal E}$.
 In a wire model, the solid wire labelled $a$ refers to a simple object of  ${\mathcal E}$ and the corresponding circle, labelled $a$ denotes its Perron--Frobenius dimension $ \mu_a$. A morphism from $a$ to $b$, induced from the action of the simple object\footnote{Here $n$ should be thought of as a multi-index since it refers to an irreducible representation of $SU(3)$.}  $n$ of ${\mathcal A}_k(SU(3))$ is represented, in the extended wire model, as a vertex $anb$, see Fig.~\ref{fig:thetaanbSU3} where the extra index $\alpha$ keeps track of multiplicities.  The corresponding theta symbol is $\theta(a,n,b; \alpha, \overline{\alpha^\prime})$.

The module action of ${\mathcal A}_k(SU(3))$ on ${\mathcal E}$ is def\/ined by a fusion graph that will be called $G$ in the last section: it is the Cayley graph describing the action of $\sigma = \sigma_{(1,0)}$ (the basic representation) on the simple objects, called $a,b,\ldots$,  of ${\mathcal E}$.  More generally, the action $n \times a = \sum_b  F_{na}^b  b$ of the simple objects $n$ of ${\mathcal A}_k$  is described by the so-called annular matrices $F_n$ obeying the same recurrence relation as the fusion matrices $N_n$, only the seed $G=F_{(1,0)}$, is dif\/ferent.

 For any fusion graph $G$ describing a quantum module ${\mathcal E}$ of type $SU(3)$, we shall obtain two coherence equations for the system of triangular cells on the graph $G$.
 In physicists' parlance, they are obtained from the two previously discussed Kuperberg  equations by coupling $SU(3)$ to ${\mathcal E}$-type matter, and taking a trace.

\subsubsection{Normalization and closure relations}\label{section3.4.1}

\begin{figure}[h]
\centering
\includegraphics{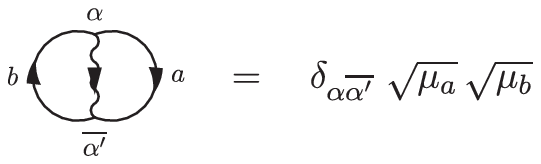}
\caption{Normalization of the fundamental intertwiner.}
\label{fig:normalizetheta}
\end{figure}

When $n$ refers to the fundamental (and basic) representation $\sigma = (1,0)$ we usually do not  put any explicit label on the corresponding oriented line of the wire model diagram, and the morphism $a \otimes \sigma \rightarrow b$  is described by  an oriented edge from $a$ to $b$ on the fusion graph of ${\mathcal E}$. This intertwiner, at the moment, is only def\/ined up to scale, so we need to f\/ix its normalization. This is done by giving a value to the corresponding theta symbol: It is convenient to set:
$\theta(a,b; \alpha,\overline{\alpha^\prime}) \equiv
 \theta(a, \sigma, b; \alpha,\overline{\alpha^\prime}) =
\delta_{\alpha \overline{\alpha^\prime}}    \sqrt{ \mu_a}    \sqrt{ \mu_b}
$, depicted by Fig.~\ref{fig:normalizetheta}.

\begin{figure}[h]
\centering
\includegraphics{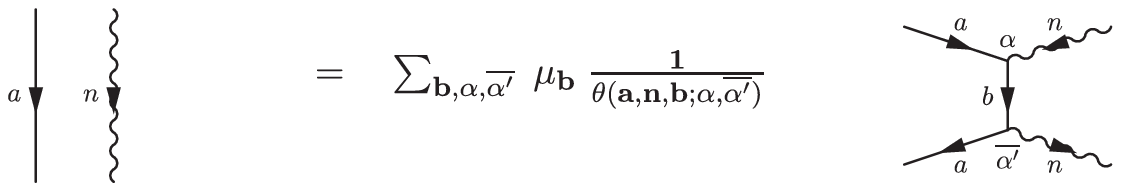}
\caption{Closure relation 1.}
\label{fig:atensorsigma}
\end{figure}

\begin{figure}[h]
\centering
\includegraphics{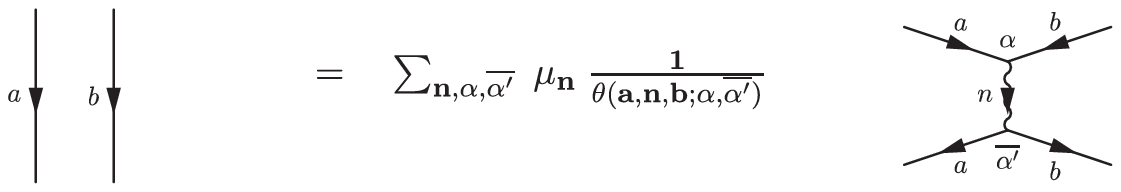}
\caption{Closure relation 2.}
\label{fig:atensorb}
\end{figure}

At some point we shall need the following closure relations:  Figs.~\ref{fig:atensorsigma} and ~\ref{fig:atensorb}. They can be seen as generalizations of those already obtained for $SU(2)$ (see
Figs.~\ref{fig:closureAA} and~\ref{fig:closureAE}). Notice that  Fig.~\ref{fig:atensorb}  implies
$\mu_a = \sum_n (\mu_n / \theta(a,n,a; \alpha, \overline{\alpha^\prime})) \mu_a \mu_a  \delta_{\alpha \overline{\alpha^\prime}}$
so that $1/\mu_a = \sum_n \mu_n /\theta(a,n,a; \alpha, \overline{\alpha})$.

\begin{figure}[h]
\centering
\includegraphics{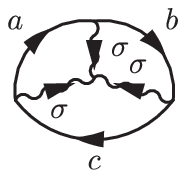}
\caption{The elementary $TET$ symbol for $SU(3)$ coupled to  ${\mathcal E}$-type  matter.}
\label{fig:TETsigma3abc}
\end{figure}

We shall not discuss the extended recoupling theory for $SU(3)$ in full generality, but, as already stressed, we have to associate a number to any elementary oriented triangle of the graph~${\mathcal E}$.
This number is a special kind of $6j$ symbol where three basic representations~$\sigma$ meet at a triple point on the wire model. It is diagrammatically denoted by the particular tetrahedral symbol
$\TET {a}{\sigma}{b}{\sigma}{c}{\sigma}$ displayed on Fig.~\ref{fig:TETsigma3abc}. One should compare it with the $SU(2)$ symbol  $\TET {a}{p}{b}{n}{c}{m}$, on the r.h.s.\ of Fig.~\ref{fig:AEEvertexwithaloop}, which is automatically $0$ when all $m$, $n$, $p$ coincide with~$\sigma$.  From now on we shall only consider tetrahedral symbols of that type, called triangular cells since their values  ${\mathcal{T}}_{abc}^{\alpha\beta\gamma}$ depend only on triangles $  \mu_{a,b,c}^{\alpha,\beta,\gamma}$ on the fusion graph of  ${\mathcal E}$, where $\alpha$, $\beta$, $\gamma$ denote oriented edges from $a$ to $b$, $b$ to $c$ and $c$ to $a$. In many cases  there is only one edge between vertices and the notation  $\tau(a,b,c)$ is suf\/f\/icient.

\subsubsection{The quadratic equation (Type~I)}\label{section3.4.2}

\begin{figure}[h]
\centering
\includegraphics{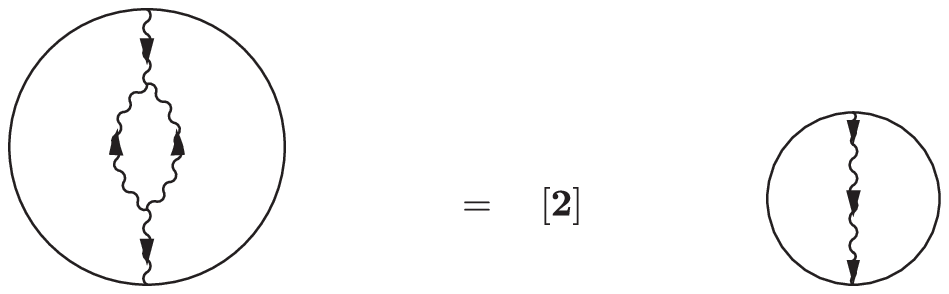}
\caption{The traced quadratic relation for $SU(3)$ coupled to ${\mathcal E}$-type matter.}
\label{fig: tracedquadraticrelation}
\end{figure}

\begin{figure}[h]
\centering
\includegraphics{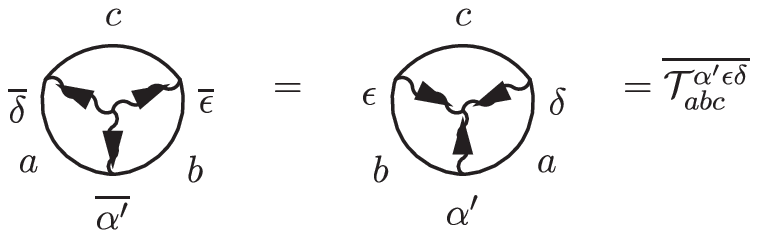}
\caption{Relation between ${\mathcal T}$ and  $\overline{\mathcal T}$.}
\label{fig:tetconjugate}
\end{figure}

We plug the f\/irst Kuperberg equation in a loop of type ${\mathcal E}$: see Fig.~\ref{fig: tracedquadraticrelation}. The right hand side (a~theta symbol) is already known. The left hand side, evaluated using f\/irst the  closure equations Fig.~\ref{fig:atensorsigma}  twice\footnote{Here $n=\sigma$, the fundamental representation.},
\begin{figure}[h!]
\centering
\includegraphics{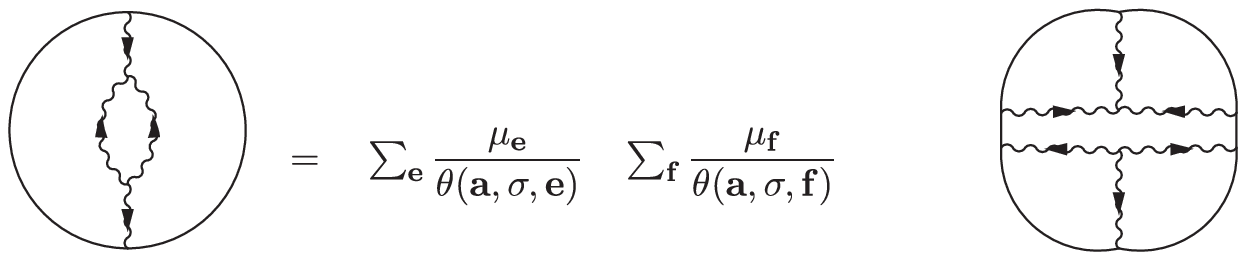}
\caption{An intermediate step.}
\label{fig: traceofatensorsigmarighthandside}
\end{figure}
\begin{figure}[h!]
\centering
\includegraphics{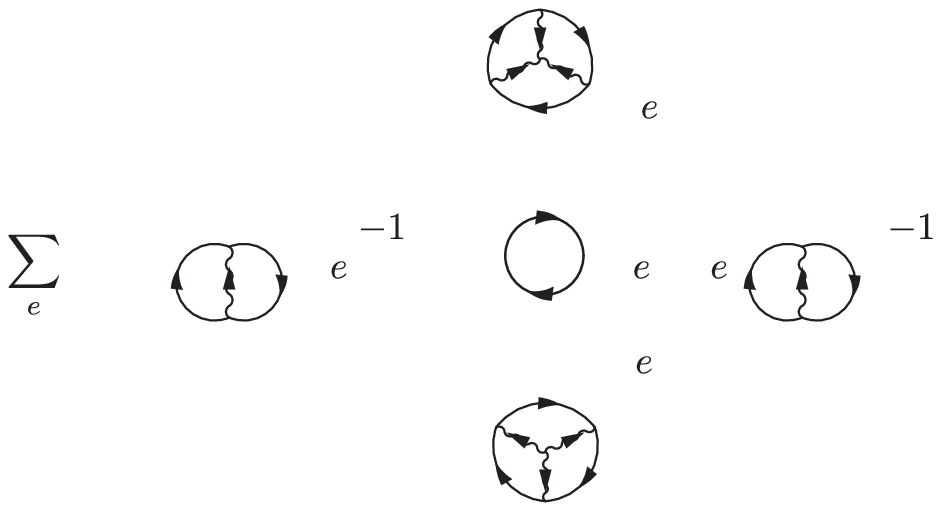}
\caption{The left hand side of the traced quadratic equation.}
\label{fig:quadraticeqlefthandside}
\end{figure}
with Fig.~\ref{fig: traceofatensorsigmarighthandside}  appearing at the end of this intermediate step, and then using the closure relation Fig.~\ref{fig:atensorb} once,  becomes a sum of products of more elementary diagrams, as described on Fig.~\ref{fig:quadraticeqlefthandside}. Explicitly, the evaluation of the l.h.s.\ reads\footnote{We use a short notation for the triangle functions $\theta$ in order not to clutter the expressions.}:
\begin{gather*}
{\rm l.h.s.} =  \sum_c \sum_d \sum_n
 \frac{\mu_c}{\theta(a,\sigma,c)}  \frac{\mu_d}{\theta(d,\sigma,b)}  \frac{\mu_n}{\theta(c,n,d)} \delta_{cd}   {\mathcal{T}}_{abc}^{\alpha\beta\gamma} \;  \overline{{\mathcal{T}}_{abd}^{\alpha^\prime \beta\gamma}} \\
\phantom{{\rm l.h.s.}}{} =  \sum_c \sum_n
\frac{\mu_c}{\theta(a,\sigma,c)}  \frac{\mu_c}{\theta(c,\sigma,b)}  \frac{\mu_n}{\theta(c,n,c)}  {\mathcal{T}}_{abc}^{\alpha\beta\gamma} \overline{{\mathcal{T}}_{abc}^{\alpha^\prime \beta\gamma}} \\
\phantom{{\rm l.h.s.}}{} =  \sum_c
\frac{\mu_c}{\theta(a,\sigma,c)}  \frac{\mu_c}{\theta(c,\sigma,b)} {\mu_c}^{-1}  {\mathcal{T}}_{abc}^{\alpha\beta\gamma} \overline{{\mathcal{T}}_{abc}^{\alpha^\prime \beta\gamma}} \\
\phantom{{\rm l.h.s.}}{} =  \sum_c
{\theta(a,\sigma,c)^{-1}}   {\mu_c}   {\theta(c,\sigma,b)^{-1}}     {\mathcal{T}}_{abc}^{\alpha\beta\gamma} \overline{{\mathcal{T}}_{abc}^{\alpha^\prime \beta\gamma}}.
\end{gather*}
On the f\/irst line, we have replaced one tetrahedral symbol by its conjugate, using the equality displayed on Fig.~\ref{fig:tetconjugate}.

The sum over $\mu_n/\theta(c,n,c)$ on the second line was replaced by ${\mu_c}^{-1}$ from the last equality obtained in Section~\ref{section3.4.1}.
We did not write explicit labels on arcs and vertices of Fig.~\ref{fig:quadraticeqlefthandside}, but a~summation should be carried out, not only on the arcs $e$ but over all ``dummy vertex labels'', i.e.\ those not already labelling the vertices of Fig.~\ref{fig: tracedquadraticrelation}. It is clear that the sequence of operations is easier to follow when it is displayed  graphically; this is what we shall do later when we come to discuss the quartic equation.

Using $\theta(a,\sigma,c)= \sqrt{a}\sqrt{c}$ in the above, the left and right hand sides give respectively
\begin{gather*}
{\rm l.h.s.}  =  \sum_c \frac{1}{\sqrt \mu_c}  \frac{1}{\sqrt \mu_b}   {\mathcal{T}}_{abc}^{\alpha\beta\gamma} \overline{{\mathcal{T}}_{abc}^{\alpha^\prime \beta\gamma}},\qquad
{\rm r.h.s.} =  [2] \sqrt{\mu_a} \sqrt{\mu_b}.
\end{gather*}
So that one obtains f\/inally the following equality:
\begin{gather*}
\sum_{c,\beta,\gamma}
{\mathcal{T}}_{a b c}^{\alpha \beta \gamma}    \overline {\mathcal{T}_{a b c}^{\alpha^\prime \beta \gamma}} = [2] \,  \delta_{\alpha \alpha^\prime}    { \mu_a}   { \mu_b}.
\end{gather*}

\subsubsection{The quartic equation (Type~II)}\label{section3.4.3}

 \begin{figure}[h]\vspace*{2mm}
\centering
\includegraphics{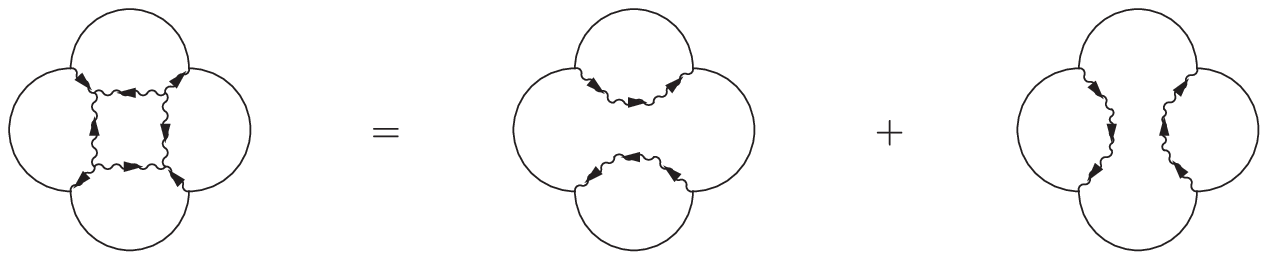}
\caption{The traced quartic relation for $SU(3)$ coupled to ${\mathcal E}$-type matter.}
\label{fig: tracedcuarticrelation}
\end{figure}
\begin{figure}[h!]
\centering
\includegraphics{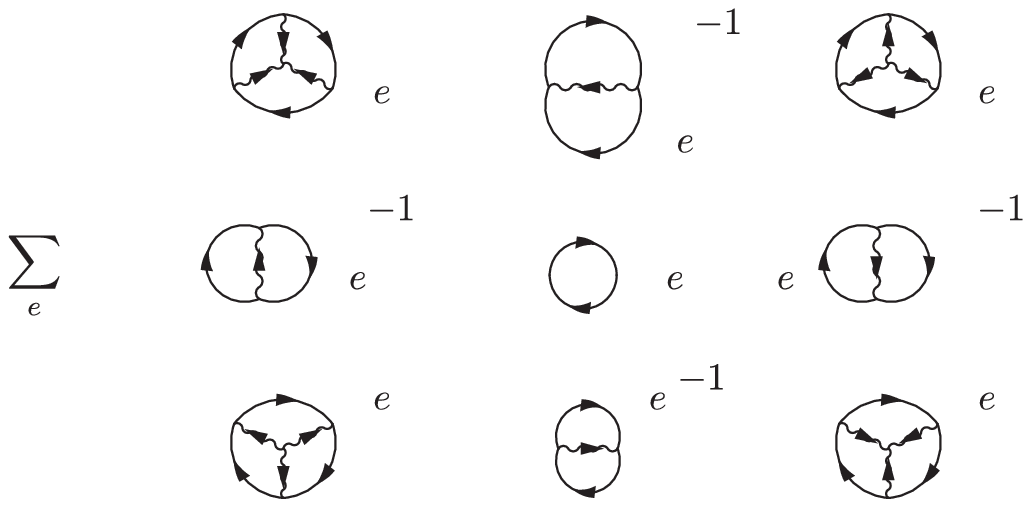}
\caption{The left hand side of the traced quartic equation.}
\label{fig:cuarticeqlefthandside}
\end{figure}

The calculation is similar for the quartic equation:  the f\/irst step is to plug the second Kuperberg equation into a loop of type ${\mathcal E}$, this is displayed on Fig.~\ref{fig: tracedcuarticrelation}.
On the left hand side,  we use the closure equations (of Fig.~\ref{fig:atensorsigma}) four times on the four sides of the square.
In this way one obtains a new graph with prefactors.
Then one uses the the other closure equations (see Fig.~\ref{fig:atensorb}) on two horizontal parallel edges, together with the fact that only $n=0$ contributes in the sum over intermediate states, so that the previous graph splits into two disconnected pieces.
Then one uses again the closure equation of Fig.~\ref{fig:atensorb}, twice, on the two pairs of parallel edges (again, only the term $n=0$ contributes).
One ends up with a product of four tetrahedral symbols, multiplied by prefactors.
Those prefactors can themselves be written in terms of theta symbols, or loops.
When written in terms of TQFT wires (graphs), the left hand side of the quartic equation therefore decomposes as a sum of products involving  four tetrahedral symbols, four inverse  ``propagators'' $\theta$ (which are proportional to Kronecker deltas) and one loop (a quantum dimension). This is displayed on Fig.~\ref{fig:cuarticeqlefthandside}.

\begin{figure}[h!]
\centering
\includegraphics{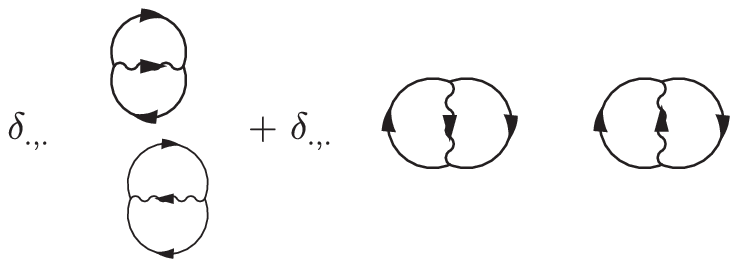}
\caption{The right hand side of the traced quartic equation.}
\label{fig:cuarticeqrighthandside}
\end{figure}

On the right hand side, the calculation is simpler.
One only has to use the closure equations of Fig.~\ref{fig:atensorb} on each of the two terms, together with the fact that only the intermediate state $n=0$ contributes. We obtain in this way a sum of two terms  (Fig.~\ref{fig:cuarticeqrighthandside}) where each term is a~product of two theta symbols. As usual, each such symbol gives a delta function and a~product of square roots of quantum dimensions.
After simplif\/ication, and using the equation of Fig.~\ref{fig:tetconjugate} that relates a tetrahedral symbol with its conjugate,  one gets\footnote{On the left hand side, one has to sum over all indices that do not appear on the right hand side.}:
\begin{gather*}
\sum \frac{1}{ \mu_e}
{\mathcal{T}}_{a e b}^{\beta_1 \beta_2 \alpha_2}
\overline {\mathcal{T}}_{a e d}^{\beta1 \beta_4 \alpha_1}
{\mathcal{T}}_{c e d}^{\beta_3 \beta_4 \alpha_4}
\overline {\mathcal{T}}_{c e b}^{\beta_3 \beta_2 \alpha_3}
=
  \delta_{\alpha_1 \alpha_2}   \delta_{\alpha_3 \alpha_4}     \mu_a     \mu_b     \mu_c +
 \delta_{\alpha_1 \alpha_4}    \delta_{\alpha_2 \alpha_3}      \mu_a     \mu_c     \mu_d.
\end{gather*}
The diagrammatic interpretation of this equation, in terms of the triangles of a fusion graph, is actually simpler than what it looks. We shall return to it later.

\subsubsection{Gauge freedom}\label{section3.4.4}

\looseness=-1
Coherence equations of Types~I and~II usually do not determine completely the cell values of the triangles of a fusion graph. More precisely, given a solution, i.e.\ a coherent assignment of complex numbers to all triangles of the graph, one can construct other solutions  by using arbitrary unitary matrices of size~$s_{ab} \times s_{ab}$, where~$s_{ab}$ denote the multiplicity of edges going from vertex~$a$ to vertex~$b$.
In particular, for graphs with single edges, a gauge choice associates phase factors to edges.
The ef\/fect of a gauge transformation is to multiply the values of triangular cells sharing a given edge (simple or multiple) by the chosen unitary factors. This will be illustrated later.

\subsubsection[$SU(3)$ braiding, Hecke and generalized JTL algebras]{$\boldsymbol{SU(3)}$ braiding, Hecke and generalized JTL algebras}\label{section3.4.5}

\begin{figure}[h]
\centering
\includegraphics{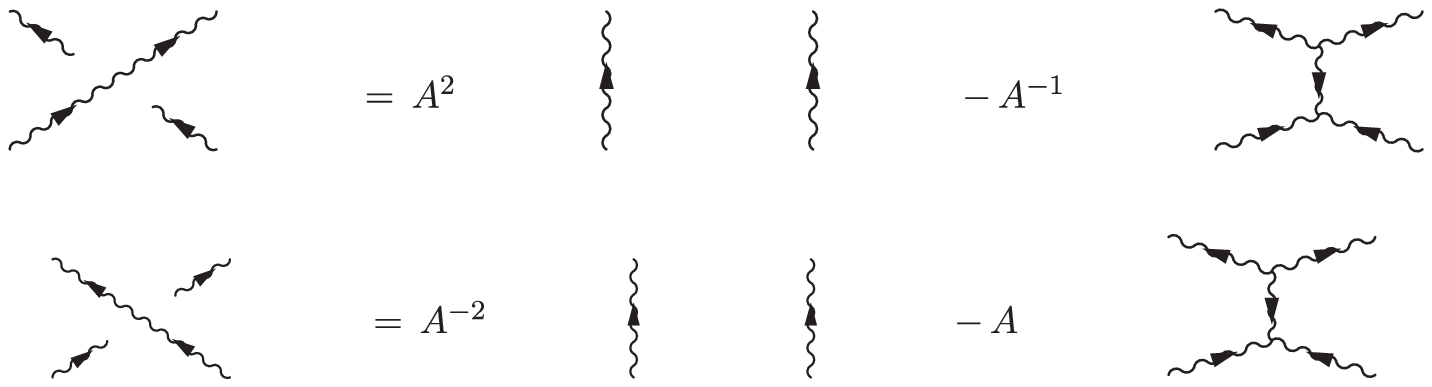}
\caption{The braid relation for $SU(3)$.}
\label{fig:braidingSU3}
\end{figure}

\begin{figure}[t]
\begin{minipage}[b]{0.45\linewidth}
\centering
\includegraphics{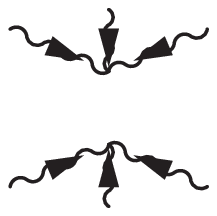}
\caption{The operator $F$.}
\label{fig: theoperatorF}\end{minipage}
\qquad
\begin{minipage}[b]{0.45\linewidth}
\centering
\includegraphics{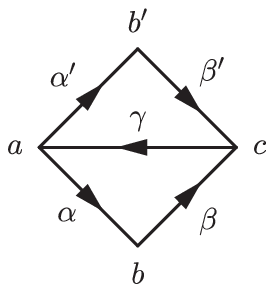}
\caption{From triangular cells to Hecke.}
\label{fig:fromcellstoHecke}
\end{minipage}
\end{figure}

In the $SU(3)$ theory, one also chooses a cubic root $A = \exp(\frac{i \pi}{3 \kappa})$ of $q$ and introduces a braiding def\/ined by  Fig.~\ref{fig:braidingSU3}.
It reads $X =  A^{2}   \one -  A^{-1}   U$ where  $U$ now denotes the double-triple vertex of $SU(3)$ that replaces the $SU(2)$ cup-cap generator.
The $SU(3)$ wire model, in its primitive version, using the graphical elements presented in Section~\ref{sec:su3graphicalelements}, is def\/ined by the three relations displayed in Figs.~\ref{fig:quadraticrelation}, \ref{fig:cuarticrelation},  \ref{fig:braidingSU3}.
The operators $U_n = \one \otimes \cdots \otimes \one \otimes U$ generate an algebra $C_\beta(SU(3))$; they still obey the def\/ining relations of the Hecke algebra with parameter $\beta = q+q^{-1}$ (see Section~\ref{section3.2.4}), in particular $U^2=\beta U$,  but with an extra relation that replaces the usual $SU(2)$ Jones--Temperley--Lieb relation. It  reads: $(U_3 U_2 U_1 - (U_1+U_3)) (U_2 U_{3} U_2  -  U_2)=0$ and can be obtained by imposing the vanishing of a quantum Young projector def\/ined in the algebra of the braid group with four strands (see  for instance \cite{DiFZuber} or  \cite{DF-meanders}). Notice that $C_\beta(SU(2)) \subset C_\beta(SU(3))$.
The operator $F = U_1 U_2 U_1 - U_1$ pictured on Fig.~\ref{fig: theoperatorF} is also of interest. Simple graphical calculations show that $F_n^2 = [2][3] F_n$,   $F_{n} F_{n\pm 1} F_{n} = [2]^2 F_{n}$ and $F_i F_j = F_j F_i$ when $\vert i - j \vert >2$, but, as noticed in \cite{EvansEtalbis},  the operators $F_{n}$ do not generate an algebra isomorphic with  $C_\beta(SU(2))$ because the last commutation relation fails for $\vert i-j \vert =2$.
As  with $SU(2)$, one can represent the $SU(3)$ intertwiners on vector spaces of paths def\/ined on the fusion graph of any module-category ${\mathcal E}$ associated with ${\mathcal A}_{k}(SU(3))$, or on the fusion graphs of ${\mathcal A}_{\infty}(SU(3))$, i.e.\ $SU(3)$ itself and its subgroups. The fusion graphs are now oriented:  one graph is associated with the action of $\sigma$, and another one, with opposite orientation, describes the action of $\overline \sigma$. From now on, for def\/initeness, we choose the f\/irst. These graphs still obey the rigidity constraint that, in the $SU(2)$ case, implies the condition  of being simply laced, but now, there may be more than one (oriented) edge $\xi_{ab}$ between vertices $a$ and $b$,  and for this reason it is usually denoted $\xi_{ab}^\alpha$.
An elementary path could be def\/ined as a succession of matching edges $\xi_{ab}^\alpha$, corresponding to the action of the $\sigma$ generator, but  it is convenient to consider more general elementary paths containing also edges of the type $\overline{\xi}_{ab}^\alpha$ corresponding to the action of the $\overline \sigma$ generator.
The $SU(2)$ ``round trip'' of length $2$, that was associated with the existence of an intertwiner from $\sigma^2$ to $\one$ is now replaced by an elementary triangle in the fusion graph, associated with the existence of an $SU(3)$ intertwiner $C$ from $\sigma^2$ to $\overline \sigma$. The action of $C: {\rm  Path}^{p} \mapsto {\rm Path}^{p+1}$ can only be def\/ined on paths containing both types of edges, as it is clear from its def\/inition,  reads:
\begin{gather*}
C_n(\xi(1)\cdots \xi(n-1) \xi_{ab}^\alpha \xi_{bc}^\beta \xi({n+2}) \cdots \xi(p))\\
\qquad{}
=
 \sum_{\gamma} \frac{
{ {\mathcal{T}}_{a b c}^{\alpha \beta \gamma}}}{\sqrt{{ \mu_a    \mu_c}}}
 (\xi(1)\cdots \xi(n-1) \overline{\xi}_{ac}^{\gamma}\xi({n+2}) \cdots \xi(p)).
\end{gather*}
 Formally, the prefactor, that was $\sqrt{{  \mu_b}/{ \mu_a}}$ for a single $SU(2)$ elementary round trip $\xi_{ab} \xi_{ba}$ is now replaced, for an $SU(3)$ elementary triangle $\xi_{ab}^\alpha \xi_{bc}^\beta \xi_{ca}^\gamma$, by the prefactor ${ {\mathcal{T}}_{a b c}^{\alpha \beta \gamma}}/\sqrt{{ \mu_a    \mu_c}}$.
Compo\-sing~$C$ with its adjoint gives the operator $U$: for every integer $n >0$,  $U_n$, as an endomorphism of  ${\rm Path}^p$, acts as follows on elementary paths with arbitrary length,  origin and extremity:
\begin{gather*}
U_n(\xi(1)\cdots \xi(n-1) \xi_{ab}^\alpha \xi_{bc}^\beta \xi({n+2}) \cdots \xi(p))\\
\qquad{} =
 \sum_{\alpha^\prime \beta^\prime}
 {\mathcal U}_{\alpha \beta}^{\alpha^\prime \beta^\prime }
 (\xi(1)\cdots \xi(n-1) \xi_{ab^\prime}^{\alpha^\prime} \xi_{b^\prime c}^{\beta^\prime} \xi({n+2}) \cdots \xi(p)),
\end{gather*}
where
\begin{gather*}
 {\mathcal U}_{\alpha \beta}^{\alpha^\prime \beta^\prime } = \sum_{\gamma} \frac{1}{ \mu_{a} \mu_{c}}
 {\mathcal{T}}_{a b c}^{\alpha \beta \gamma}
\overline{{\mathcal{T}}_{a b^\prime c}^{\alpha^{\prime}  \beta^{\prime}  \gamma }}.
\end{gather*}
The last summation runs over all edges $\gamma$ such that $(\xi_{ab}^\alpha \xi_{bc}^\beta \xi_{ca}^\gamma)$ and $(\xi_{ab^\prime}^{\alpha^\prime} \xi_{b^\prime c}^{\beta^\prime} \xi_{ca}^\gamma)$ make a rhombus with diagonal $\gamma$, i.e.\ a pair of elementary triangles sharing the edge $\xi_{ca}^\gamma$ from~$c$ to~$a$.

A matching pair  of edges $(\alpha \beta)$ determines three vertices $a$, $b$, $c$ of a triangle, so that when the triangular cells $ {\mathcal{T}}_{a b c}^{\alpha \beta \gamma}$ are known, the above equation associates to every pair $(a,c)$ of neighboring vertices of the fusion graph an explicit square matrix whose lines and columns are labelled by triangles sharing the given two vertices. Its matrix elements associate $\CC$-numbers to the various, possibly degenerated, rhombi with given diagonal (only its endpoints are given).
This representation of the Hecke algebra can sometimes be obtained by other means: there are general formula for representations associated with ${\mathcal A}_k(SU(N))$ itself,  given in~\cite{Wenzl:An}, and all Hecke representations associated with modules over ${\mathcal A}_k(SU(3))$ were obtained by ``computer aided f\/lair''  in~\cite{DiFZuber}. The point to be stressed is that the obtained representation is now {\it deduced} from the values of triangular cells. This was discussed in \cite{EvansSU3Cells}.
For illustration we shall only give one example of this construction (see the end of Section~\ref{section4.3.2}).

\section{Triangular cells on fusion graphs}\label{section4}

A motivation for introducing a good part of the material described so far was to explain and justify the coherence equations of Types~I and~II.
In the coming section we shall make no explicit use of the TQFT wire models but we shall see on specif\/ic examples, i.e.\ on specif\/ic fusion graphs,  how coherence equations can be solved.
The remaining part of this article can therefore be read almost independently of the previous sections.
We shall consider the following $\mathcal{A}_k (SU(3))$ modules:  $\mathcal{A}_k$ itself, $\mathcal{E}_5$, $\mathcal{E}_9$, $\mathcal{E}_{21}$.
We shall also study the graph $\mathcal{Z}_9$ that turns out to be rejected as describing an $SU(3)$ module.

\subsection{Self-connection on graphs}\label{section4.1}

An $SU(3)$ module $\mathcal{E}$ at level $k$ is specif\/ied by its fusion graph, also called  $\mathcal{E}$.
The adjacency matrix, called $G$,  encodes the action of the fundamental generator $\sigma = \sigma_{(1,0)}$ of $\mathcal{A}_k (SU(3))$ on the simple objects $a,b,c,\dots$ of $\mathcal{E}$, represented as vertices.
There exist a Perron--Frobenius positive measure on the set of vertices, it gives  their quantum dimension denoted\footnote{Warning: in the previous sections, this was the notation for $q$-numbers but there should be no confusion.}  $[a]=\mu_a$ for the vertex $a$.

We consider the following oriented triangles of the fusion graph:
\[
\Tribis{a}{b}{c}{\alpha}{\beta}{\gamma}{}
 = \Delta_{abc}^{\alpha\beta\gamma},
\]
where $\alpha$, $\beta$, $\gamma$ denote oriented edges.
A cell-system $\mathcal{T}$ on a fusion graph $G$ is a map from the set of oriented triangles of $G$ to the complex numbers, that we denote as follows:
\begin{equation*}
\mathcal{T} \left(\,  \Tribis{a}{b}{c}{\alpha}{\beta}{\gamma}{} \right) = \Tribis{a}{b}{c}{\alpha}{\beta}{\gamma}{\tau} = \mathcal{T}(\Delta_{abc}^{\alpha \beta \gamma}) = \mathcal{T}_{abc}^{\alpha \beta \gamma}  .
\end{equation*}
The complex numbers $\mathcal{T}_{abc}^{\alpha \beta \gamma}$ are the triangular cells of $G$ or simply ``cells'' in the sequel, since they will all be of that type.
As discussed previously, in order to guarantee the existence of a~module-category $\mathcal{E}$ described by the fusion graph $G$, the set of triangular cells associated with this graph has to obey the following  coherence equations.

\vspace{-3mm}

\paragraph{Type  I equations.}
For each frame\footnote{It may be degenerated, it is then a single edge.} $\frameI{a}{b}{\alpha}{\alpha'}$ of the graph $G$, {ie} a double edge $\alpha$, $\alpha'$ from $a$ to $b$,
we have a quadratic equation (or ``small pocket equation''):
\begin{equation*}
\mathcal{T} \left( \Tridouble{a}{b}{c}{\alpha}{\beta}{\gamma}{\alpha'} \right)=
\sum_{c,\beta,\gamma} \  \Tribis{a}{b}{c}{\alpha}{\beta}{\gamma}{\tau}    \overline{ \Tribis{a}{b}{c}{\alpha'}{\beta}{\gamma}{\tau}\;\;} = 2_q [a] [b].
\end{equation*}

\paragraph{Type  II  equations.}
For each frame\footnote{It may also be degenerated: see the discussion later.} \basisBig{a_1}{a_2}{a_3}{a_4}{\alpha_1}{\alpha_4}{\alpha_3}{\alpha_2} of the graph $G$,
we have a quartic equation (or ``big pocket equation''):
\begin{gather}
\sum_{c,\beta_1,\beta_2,\beta_3,\beta_4} \quad \frac{1}{[c]} \mathcal{T} \left( \frameIIIeq{a_1}{a_2}{a_3}{a_4}{c}{\alpha_1}{\alpha_4}{\alpha_3}{\alpha_2} \right)
\nonumber\\
\nonumber \qquad{} =  \sum_{c,\beta_1,\beta_2,\beta_3,\beta_4} \quad \frac{1}{[c]} \Tribis{a_1}{a_2}{c}{\alpha_1}{\beta_2}{\beta_1}{\tau} \qquad \overline{\Tribis{a_3}{a_2}{c}{\alpha_2}{\beta_2}{\beta_3}{\tau} } \qquad
\Tribis{a_3}{a_4}{c}{\alpha_3}{\beta_4}{\beta_3}{\tau} \qquad  \overline{\Tribis{a_1}{a_4}{c}{\alpha_4}{\beta_4}{\beta_1}{\tau}} \\
\qquad{} =  [a_1][a_2][a_3]\delta_{\alpha_1,\alpha_4}\delta_{\alpha_2,\alpha_3} +  [a_1][a_2][a_4]\delta_{\alpha_1,\alpha_2}\delta_{\alpha_3,\alpha_4}.
\label{typeIIeq}
\end{gather}

\paragraph{Gauge equivalence.}
Let $\mathcal{T}$ be a cell system for the graph ${G}$. For any pair of vertices $a$, $b$ of~${G}$, let the integer $s_{ab}$ denote the edge multiplicity from $a$ to $b$ and let $U^{ab}$ be a unitary matrix of size $s_{ab} \times s_{ab}$. We def\/ine the complex numbers:
\begin{equation}
\mathcal{T}' \left(\, \Tribis{a}{b}{c}{\alpha'}{\beta'}{\gamma'}{} \right) = \sum_{\alpha,\beta,\gamma} \mathcal{T} \left( \, \Tribis{a}{b}{c}{\alpha}{\beta}{\gamma}{}\right) (U^{ab})_{\alpha'\alpha} \, (U^{bc})_{\beta'\beta} \,(U^{ca})_{\gamma'\gamma}.
\label{gauge}
\end{equation}
These complex numbers satisfy Type~I  and Type~II  equations,  therefore $\mathcal{T}'$ is also a cell system for the graph ${G}$. The two cell systems are called gauge equivalent. An equivalence class of cell systems is called a self-connection on ${G}$. A quantity is an invariant of the cell-system if its value does not depend on the gauge choices.

\subsection{Solving the cell system}\label{section4.2}

 The number of edges from a given vertex $a$ to a given vertex $b$ is equal to the matrix-ele\-ment~$({G})_{ab}.$ The number of cells of ${G}$ is $\frac{1}{3} \sum_a ({G}{G}{G})_{aa} = \frac{1}{3}\textrm{Tr}({G}^3)$. These complex numbers have to satisfy Type~I  and Type~II  equations. There is one Type~I  equation for each pair of vertices linked by one (or more) edge(s). So the number of Type~I  equations is given by the number of non-zero matrix elements of ${G}$. For Type~II  equations, the number of frames is given by $\textrm{Tr}({G}{G}^T{G}{G}^T)$ where ${G}^T$ is the transpose matrix of~${G}$. But the following frames (related by the two diagonal symmetries):
\[
\basisBig{a_1}{a_2}{a_3}{a_4}{\alpha_1}{\alpha_4}{\alpha_3}{\alpha_2}
\qquad \qquad
\basisBig{a_3}{a_2}{a_1}{a_4}{\alpha_2}{\alpha_3}{\alpha_4}{\alpha_1}
\qquad \qquad
\basisBig{a_1}{a_4}{a_3}{a_2}{\alpha_4}{\alpha_1}{\alpha_2}{\alpha_3}
\]
lead to the same equations, so the number of Type~II  equations is reduced\footnote{Even after considering diagonal symmetries, two dif\/ferent frames can lead to the same equation. In some cases, there is no equation associated with a frame of this type: this happens when the fusion graph is such that there is no vertex~$c$ on the top of the pyramid in equation~(\ref{typeIIeq}).}.

The number of Type~I  and Type~II  equations is much larger than the number of cells that one has to compute. Nevertheless, these equations do not  f\/ix completely their values and one can make use of gauge freedom for that. In the examples treated, we try to choose the most convenient solution (either by imposing real values for the  cells, if possible, or by making a~gauge choice exhibiting the symmetry of the graph). The resolution technique for a cell system depends largely on the chosen example.

\subsubsection{Graphs with single edges}\label{section4.2.1}

Most fusion graphs of $SU(3)$ have only single edges (exceptions are  the $\mathcal{A}^*$ and $\mathcal{D}_{3s}$ series, the special twisted conjugate $\mathcal{D}_{9}^{tc}$ and the exceptional $\mathcal{E}_9$). Since there is no multiplicity, we can drop the edge labeling symbols and denote the cells simply by $\mathcal{T}_{abc}$. The cell values, on a given triangle, for two gauge equivalent cell systems $\mathcal{T}$ and $\mathcal{T}'$  will only dif\/fer by phase factors since the unitary matrices in (\ref{gauge}) are one-dimensional in this case.

\paragraph{(i) Type~I  equations.} Since $\alpha = \alpha'$, these equations reduce to:
\[
\sum_c \left| \mathcal{T}_{abc} \right|^2 = 2_q [a] [b] ,
\]
and only give quadratic constraints on the modulus of the cells.
\paragraph{(ii) Type~II  equations.} These equations read:
\begin{equation}
\sum_c \frac{1}{[c]} \mathcal{T}_{a_1a_2c}   \overline{\mathcal{T}_{a_3a_2c}}   \mathcal{T}_{a_3a_4c}   \overline{\mathcal{T}_{a_1a_4c}} = [a_1][a_2][a_3]\delta_{a_2,a_4} +  [a_1][a_2][a_4]\delta_{a_1,a_3}.
\label{TypeIIfulldeg}
\end{equation}

\paragraph{Doubly degenerated case ($\boldsymbol{a_1=a_3}$ and $\boldsymbol{a_2=a_4}$).} Here both terms on the right hand side of~(\ref{TypeIIfulldeg}) contribute. The equations read:
\[
\sum_c \frac{1}{[c]} \left| \mathcal{T}_{a_1a_2c} \right|^4 = [a_1]^2[a_2] + [a_1][a_2]^2
\]
and give quartic constraints for the modulus of the cells.

\paragraph{Singly degenerated cases ($\boldsymbol{a_1=a_3}$ or $\boldsymbol{a_2=a_4}$).} Here only the f\/irst or the second term of~(\ref{TypeIIfulldeg}) remains, and the equations read:
\begin{gather*}
  \sum_c \frac{1}{[c]} \left| \mathcal{T}_{a_1a_2c} \right|^2   \left| \mathcal{T}_{a_1a_4c} \right|^2 = [a_1][a_2][a_4] \qquad   (a_1 = a_3), \\
  \sum_c \frac{1}{[c]} \left| \mathcal{T}_{a_1a_2c} \right|^2   \left| \mathcal{T}_{a_3a_2c} \right|^2 = [a_1][a_2][a_3] \qquad  (a_2 = a_4).
\end{gather*}
Here also the equations only give constraints for the modulus of the cells.

\paragraph{Non degenerate case.}  Here the right hand side of~(\ref{TypeIIfulldeg}) vanishes and the Type~II  equations give therefore constraints on the phase factors of the cells, once the values of their modulus have been f\/ixed by using the previous equations.

\subsubsection{Graphs with double edges}\label{section4.2.2}
For these graphs,  because of the possible multiplicities, the summations to be performed in Type~I  and~II equations can be quite involved. We refer to the~$\mathcal{E}_9$ case treated below for an example of the discussion.

\subsection[Classical and quantum $SU(3)$:  cells for  ${\cal A}_k$ graphs]{Classical and quantum $\boldsymbol{SU(3)}$:  cells for  $\boldsymbol{{\cal A}_k}$ graphs}\label{section4.3}

\subsubsection[Cells for  the ${\cal A}_k$ graphs]{Cells for  the $\boldsymbol{{\cal A}_k}$ graphs}\label{section4.3.1}

The fusion graph of ${\cal A}_k(SU(3))$ (Weyl alcove) appears as the truncation at level $k$ of its classical analogue, the  ${\cal A}_\infty$ graph. Altitude is $\kappa = k+3$, i.e.\ $q^\kappa=-1$.
Using triangular coordinates (distances to the walls of the Weyl alcove), the quantum dimension\footnote{In this section $q$-numbers are just denoted $[n]$ rather than~$n_q$.} of an irreducible representation $(k,l)$ reads simply $\mu_{k,l}=[k+1][l+1][k+l+2]/[2]$.
Consider now the rhombus with vertices $(k,l)$, $(k+1,l-1)$, $(k+1,l)$, $(k,l+1)$ made of two opposite triangles $\tau^\uparrow$ and $\tau^\downarrow$ sharing the (horizontal) diagonal $(k,l)$, $(k+1,l)$.
The Type~I  equations read
$\vert \tau^\uparrow \vert^2 + \vert \tau^\downarrow \vert^2 = [2]   \mu_{k,l}   \mu_{k+1,l}$.
The Type~II equations (degenerated subcases) applied to the same two triangles read
$\vert \tau^\uparrow \vert^4 / \mu_{k,l+1} + \vert \tau^\downarrow \vert^4  / \mu_{k+1,l-1} = \mu_{k,l}    \mu_{k+1,l}  (\mu_{k,l} +  / \mu_{k+1,l} )$
and
$\vert \tau^\uparrow \vert^2   \vert \tau^\downarrow \vert^2 / \mu_{k,l} =  \mu_{k,l+1}    \mu_{k+1,l}    \mu_{k+1,l-1}$.
Dividing the second equation by the third gives a second degree equation $ \mu_{k+1,l-1}   \lambda +   \mu_{k,l+1}   1/\lambda  = \mu_{k,l} +  \mu_{k+1,l}$, with  $\lambda = \vert \tau^\uparrow \vert^2 / \vert \tau^\downarrow \vert^2$, leading to the\footnote{The other root $\lambda = \frac{[k + 1] [3 + k + l]}{[k + 2] [2 + k + l]}$ leads to a solution that is not compatible with some of the Type~II equations (the doubly degenerate ones) and should be rejected.}
root $\lambda = [l+2]/[l]$. Def\/ine $f(k,l)=\vert \tau^\uparrow \vert^2 [2]/ [l+2]=\vert \tau^\downarrow \vert^2 [2]/ [l]$.
The l.h.s.\ of the f\/irst (quadratic) equation, using $[2] [l+1] = [l]+[l+2]$,
 gives $\vert \tau^\uparrow \vert^2 + \vert \tau^\downarrow \vert^2 =  f(k,l) ([l]+[l+2])/[2] = f(k,l) [l+1]$ and the r.h.s. is $[2] \, \mu_{k,l} \, \mu_{k+1,l}.$
 This determines immediately $f(k,l)$, hence also $\vert \tau^\uparrow \vert $ and  $\vert  \tau^\downarrow \vert $.
For  a triangle $(k,l) \mapsto (k+1,l) \mapsto (k,l+1)$, i.e.\ a ``pointing up  triangle'', one f\/inds:
\[
\vert \tau^\uparrow \vert^2 =  \vert \tau_{k,l}^\uparrow \vert^2  =  \frac{[k+1][k+2][l+1][l+2][k+l+2][k+l+3]}{[2]^2}.
\]

For  a triangle $(k,l)  \mapsto (k+1,l-1)  \mapsto (k+1,l)$, i.e.\  a ``pointing down  triangle''
\[
\vert \tau^ \downarrow \vert^2 = \vert \tau_{k,l}^\downarrow \vert^2 = \frac{[k+1][k+2][l][l+1][k+l+2][k+l+3]}{[2]^2}.
\]
Using obvious gauge choices it is easy to make the cells real.
Another proof of the same general formulae was given in~\cite{EvansSU3Cells}.

Another way to compute the cell values is to use a simple recurrence, starting from the f\/irst triangle $(0,0) \mapsto (1,0) \mapsto (0,1)$ and using the symmetries of the graph.
For illustration we consider below the ``classical case'', i.e.\  the  ${\cal A}_\infty$ graph corresponding to the Lie group $SU(3)$ itself.
It is enough to use the Type~I  coherence equations. The calculation goes as follows:
\begin{itemize}\itemsep=-1pt
\item Edge $(0,0)\mapsto (1,0)$. Applying $I$ to $\tau = (0,0) \mapsto (1,0) \mapsto (0,1)$ gives $\vert \tau \vert^2 = 2 \times 1 \times 3 = 6$.
\item Edge $ (1,0) \mapsto (0,1)$. Applying $I$ to $ (0,0) \mapsto (1,0) \mapsto (0,1)$ and $\tau = (1,1) \mapsto (1,0) \mapsto (0,1)$ sharing the common edge $(1,0) \mapsto (0,1)$, gives $6 + \vert \tau \vert^2  = 2\times 3\times 3$ hence $\vert \tau \vert^2  = 12$.
\item Edge $(0,1)\mapsto (1,1)$. Applying $I$ to  $(1,1) \mapsto (1,0) \mapsto (0,1) $ and $\tau = (1,0) \mapsto (2,0) \mapsto (1,1)$ sharing the common edge $(1,1) \mapsto (1,0)$, gives $ 12 + \vert \tau \vert^2  = 2 \times 3 \times 8 $ hence $\vert \tau \vert^2  = 36$.
\item Edge $(2,0) \mapsto (1,1)$. Applying $I$ to  $(1,0) \mapsto (2,0) \mapsto (1,1)$ and $\tau = (2,0) \mapsto (1,1) \mapsto (2,1)$ sharing the common edge $(2,0) \mapsto (1,1)$, gives $ 36 + \vert \tau \vert^2  = 2 \times 6 \times 8 $ hence $\vert \tau \vert^2  =  60$.  Etc.
\end{itemize}

\begin{figure}[htbp]
\centering
\scalebox{0.50}{\includegraphics{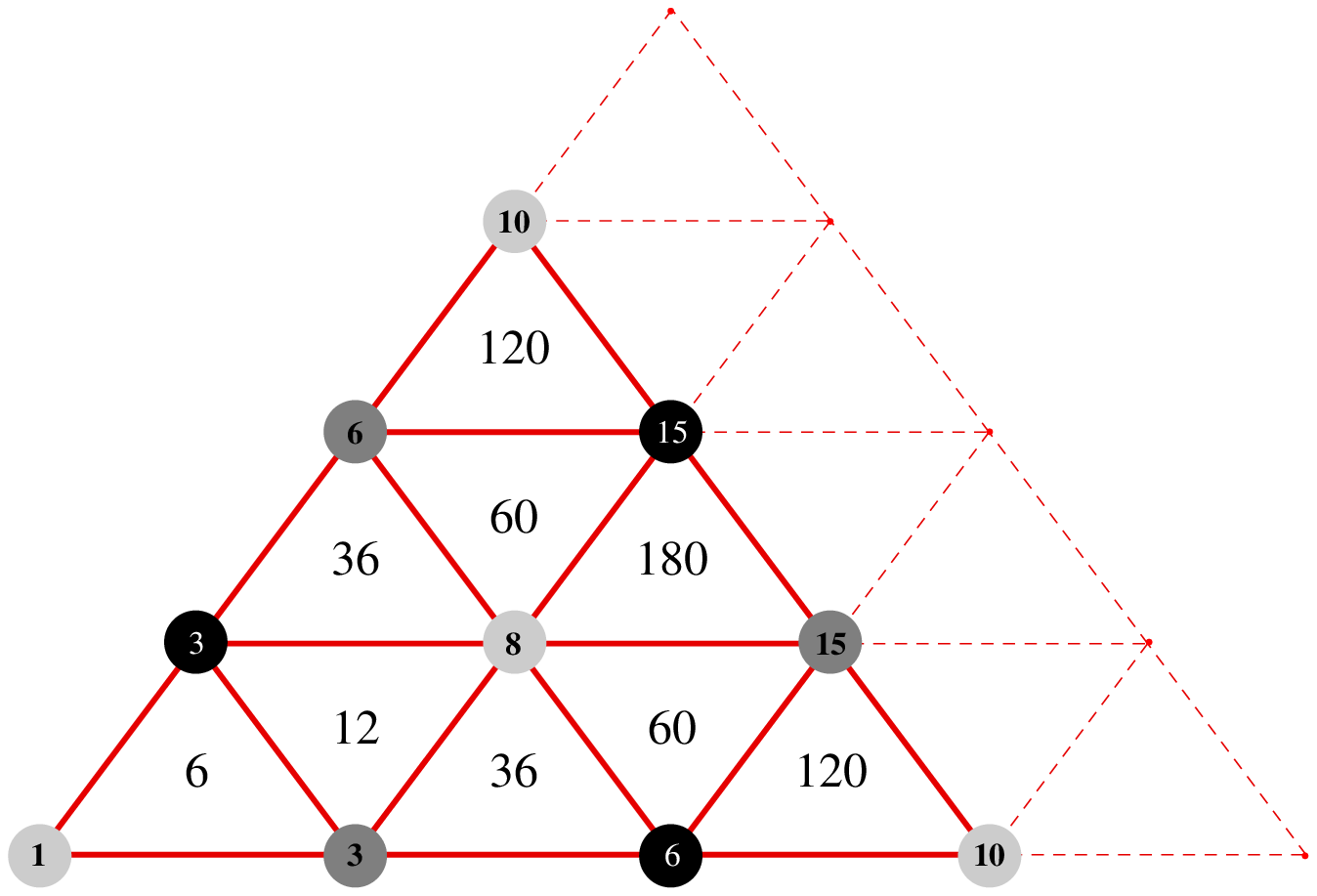}}
\caption{Square moduli of cells for classical $SU(3)$.}
\label{fig: cellsforclassicalsu3}
\vspace{-1mm}
\end{figure}

\looseness=-1
We can recover these values from the previous general formulae by taking $q=1$. Indeed the coherence equations~I  and~II    hold not only for the quantum $SU(3)$ graphs, but also for the classical ones, i.e.\ for the fusion graph of  ${\cal A}_\infty$ and for fusion graphs (that may be called $SU(3)$ McKay graphs) associated with all the f\/inite subgroups of $SU(3)$.
The case of $SU(3)$ itself is displayed on Fig.~\ref{fig: cellsforclassicalsu3}.
The numbers inside the triangles give the squared absolute values $\vert \tau \vert^2$ of the cells.

\subsubsection{From triangular cells to Hecke: an example}\label{section4.3.2}

\begin{figure}[t]
\centerline{\scalebox{0.65}{\includegraphics{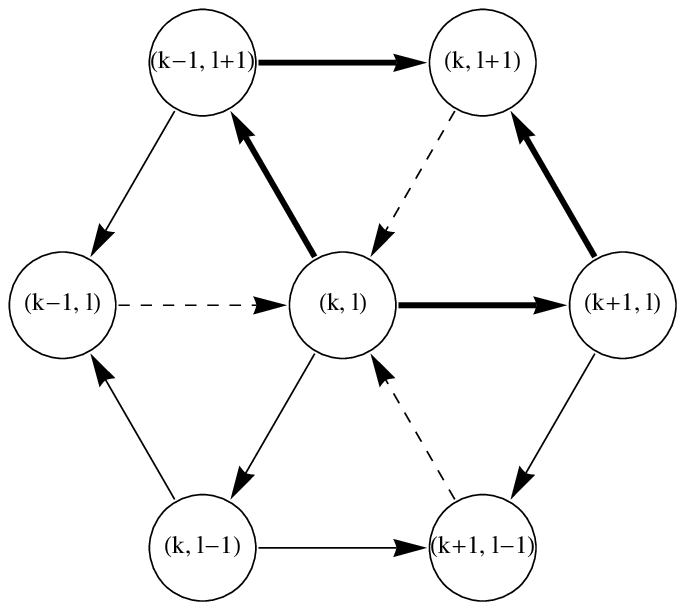}}}
\caption{${\mathcal A}_k(SU(3))$: three rhombi.}
\label{fig:threerhombus}
\end{figure}

As announced in Section~\ref{section3.4.5}, we show on this example how to deduce matrices representing the Hecke algebra from the obtained self-connection.
Here all edges are simple, so that choosing a~pair of neighboring vertices amounts to choose an edge of the fusion graph;  for the same reason, there is no summation over $\gamma$ in the formula for ${\mathcal U}$ given at the end of  Section~\ref{section3.4.5}. Consider the edge $(k,l) \gets (k,l+1)$. Naively there is only one rhombus with this edge as a diagonal, but there are actually four of them, i.e.\ four matrix elements for a $2\times 2$ matrix ${\mathcal U}$,  since one should consider degenerated cases and also take order into account. The f\/irst degenerated rhombus (set $\alpha=\alpha^\prime= (k,l) \to (k+1,l)$, $\beta=\beta^\prime=(k+1,l)\to (k,l+1)$ in the expression of ${\mathcal U}$) gives ${\mathcal U}_{11}=\frac{1}{\mu_{k,l}\mu_{k,l+1}} \tau^\uparrow_{k,l} \tau^\uparrow_{k,l}= \frac{[k+2]}{[k+1]}$. The next degenerated rhombus (set $\alpha=\alpha^\prime= (k,l) \to (k-1,l+1)$, $\beta=\beta^\prime=(k-1,l+1)\to (k,l+1)$) gives ${\mathcal U}_{22}=\frac{1}{\mu_{k,l}\mu_{k,l+1}} \tau^\downarrow_{k-1,l+1} \tau^\downarrow_{k-1,l+1}= \frac{[k]}{[k+1]}$. The two non diagonal terms (they are equal since we choose the phases to make ${\mathcal T}$ real) correspond to the non-degenerated case $((\alpha \beta ),(\alpha^\prime, \beta^\prime))$ displayed on Fig.~\ref{fig:threerhombus}. Their value is ${\mathcal U}_{12}={\mathcal U}_{21}= \frac{1}{\mu_{k,l}\mu_{k,l+1}} \tau^\downarrow_{k,l} \tau^\downarrow_{k-1,l+1}=\frac{\sqrt{[k][l+2]}}{[l+1]}$. To the edge $(k,l) \gets (k,l+1)$ is therefore associated a $2\times 2$ matrix ${\mathcal U} =
\left(\begin{array}{c c}
{\mathcal U}_{11} & {\mathcal U}_{12}\\  {\mathcal U}_{21} & {\mathcal U}_{22}  \end{array}\right)
$.
In the same way one can build two other $2\times 2$ matrices, i.e.\ a~total of three, see Fig.~\ref{fig:threerhombus}, for each of the three edges $(k,l) \gets (k,l+1)$,   $(k,l) \gets (k-1,l)$, $(k,l) \gets (k+1,l-1)$ radiating towards the vertex $(k,l)$.
Other explicit expressions for  ${\mathcal U}$ matrices can be found in \cite{DiFZuber} and \cite{EvansSU3Cells}.


\subsection[The $\mathcal{E}_5$ graph]{The $\boldsymbol{\mathcal{E}_5}$ graph}\label{section4.4}

\begin{figure}[H]
\centerline{\scalebox{0.55}{\includegraphics{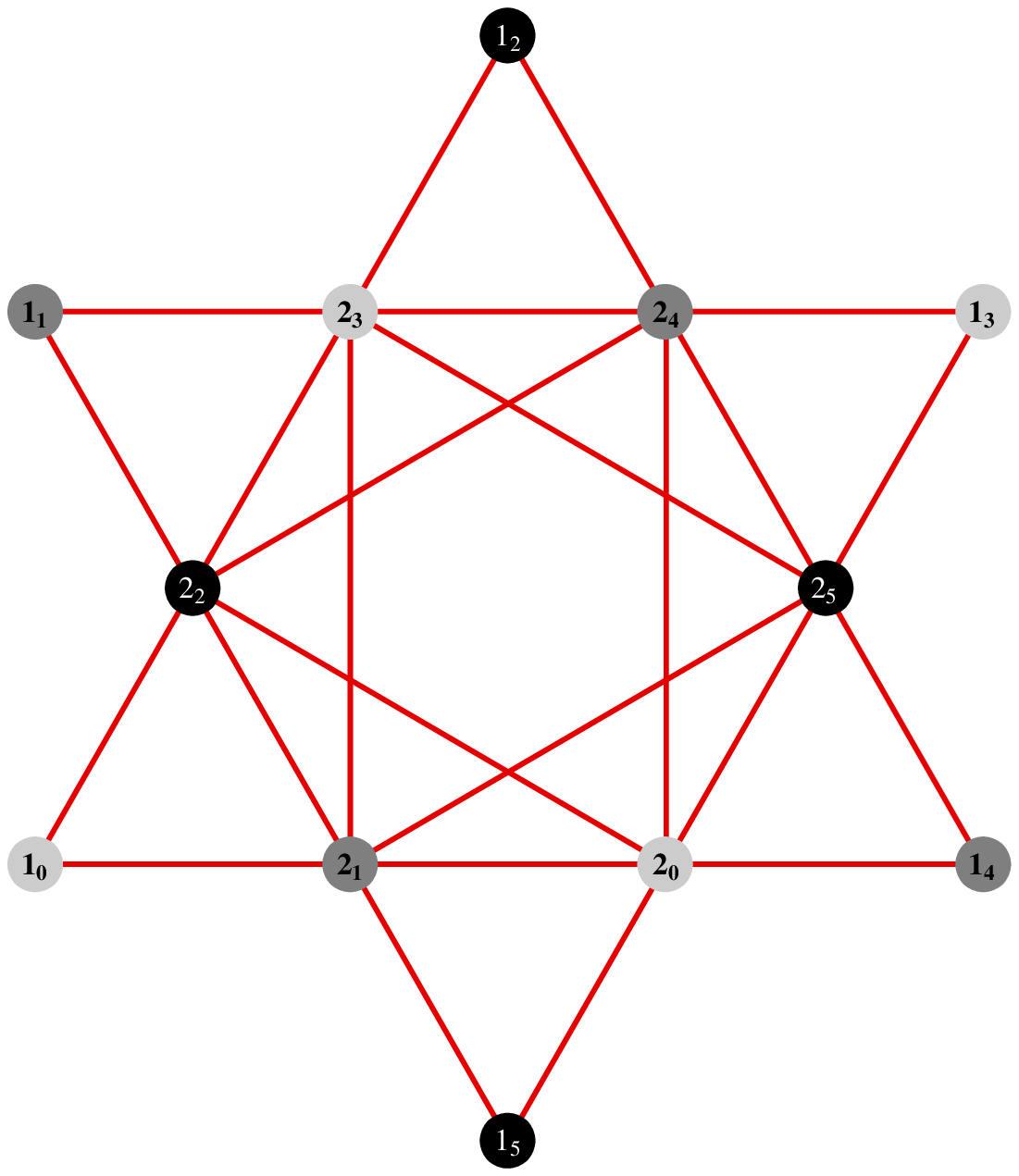}}}
\caption{The $\mathcal{E}_5$ fusion graph.}
\label{E5-graph}
\end{figure}

$\mathcal{E}_5$ is the graph of level $k=5$ and altitude $\kappa=k+3=8$, i.e.\ $q^8=-1$, displayed in Fig.~\ref{E5-graph}. It has 12 vertices that we denote $1_i$ and $2_i$, with $i$ from 0 to 5. The quantum dimensions of the vertices are $\left[1_i\right]=1$, $\left[2_i\right]=3_q= 1 + \sqrt{2}$.

The number of oriented triangles (and therefore the number of cells) of the graph is $\textrm{Tr}({G})/3 = 14$. There are 6 external triangles with value denoted $\tau_i$, 6 internal triangles with value denoted~$\mu_i$ and~2 larger internal triangles with value denoted~$\nu_0$ and~$\nu_1$:
$$
\tau_i = \Tri{1_i}{2_{i+1}}{2_{i+2}}{\tau}  \qquad \qquad   \mu_i=\Tri{2_i}{2_{i+1}}{2_{i+2}}{\tau} \qquad\qquad \nu_0=\Tri{2_0}{2_2}{2_4}{\tau} \qquad \qquad \nu_1=\Tri{2_1}{2_3}{2_5}{\tau}
$$
The symmetry of the graph suggests (but it can also be deduced from the equations) that the six cells $\tau_i$ have equal modulus. This is also the case for the six cells $\mu_i$, and for the two last ones. We set: $|\tau_i|=|\tau|$, $|\mu_i|=|\mu|$ for $i=0$ to $5$ and $|\nu_0|=|\nu_1|=|\nu|$. These numbers are invariants of the cell-system.

The number of Type~I  frames is $\sum_{a,b} ({G})_{ab}=24$. With the above symmetry considerations, the number of independent equations associated with Type~I  frames is reduced to 3. The Type~I  frames and the corresponding equations are:
\begin{alignat*}{3}
 & \frameIsingle{1_i}{2_{i+1}} \qquad \textrm{and} \qquad \frameIsingle{2_i}{1_{i-2}} \qquad   & &  |\tau|^2  =  2_q 3_q&  \\
& \frameIsingle{2_i}{2_{i+1}} \qquad \qquad                                         & &  |\tau|^2 + 2   |\mu|^2  =  2_q 3_q 3_q  & \\
& \frameIsingle{2_i}{2_{i+4}} \qquad \qquad                                         & &  |\mu|^2 + |\nu|^2  = 2_q 3_q 3_q&
\end{alignat*}
We have $2_q = \sqrt{2+\sqrt{2}}$ and $3_q = 1+\sqrt{2}$, so we obtain:
\begin{equation}
|\tau|^2 = \sqrt{10+7\sqrt{2}}, \qquad |\mu|^2 = \sqrt{5+\frac{7}{2}\sqrt{2}}, \qquad |\nu|^2 = \sqrt{29+\frac{41}{2}\sqrt{2}}.
\label{E5squarecells}
\end{equation}

The number of Type~II  frames is $\textrm{Tr}({G}{G}^{T}{G}{G}^{T})=108$. There are 24 double degenerated frames, but considering symmetries, we end up with 3 independent equations:
\begin{alignat*}{3}
& \basisBig{1_i}{2_{i+1}}{1_i}{2_{i+1}}{}{}{}{} \qquad \textrm{and} \qquad \basisBig{2_{i}}{1_{i-2}}{2_{i}}{1_{i-2}}{}{}{}{} \quad \qquad && \frac{1}{3_q} |\tau|^4  =  3_q + 3_q 3_q & \\[0.9ex]
& \basisBig{2_i}{2_{i+1}}{2_i}{2_{i+1}}{}{}{}{}  \qquad && |\tau|^4 + \frac{1}{3_q} \big(|\mu|^4 + |\mu|^4\big)  =  2(3_q3_q3_q) & \\[0.9ex]
& \basisBig{2_i}{2_{i+4}}{2_i}{2_{i+4}}{}{}{}{}  \qquad&& \frac{1}{3_q}\big(|\mu|^4+|\nu|^4\big)    =   2(3_q3_q3_q)&
\end{alignat*}
There are 72 singly degenerated frames, 36 of them with the same vertices on the second diagonal, they lead to the following two independent equations:
\begin{alignat*}{3}
& \basisBig{1_i}{2_{i+1}}{2_i}{2_{i+1}}{}{}{}{} \qquad \textrm{and} \qquad \basisBig{1_i}{2_{i+1}}{2_{i+3}}{2_{i+1}}{}{}{}{} \quad \qquad&& \frac{1}{3_q} |\tau|^2|\mu|^2  = 3_q 3_q &\\[0.9ex]
& \basisBig{2_i}{2_{i+1}}{2_{i+3}}{2_{i+1}}{}{}{}{}  \qquad&& \frac{1}{3_q} \big(|\mu|^4 + |\mu|^2|\nu|^2\big)  = 3_q3_q3_q&
\end{alignat*}
The other singly degenerated frames, with the same vertices in the f\/irst diagonal, are:
\[
\basisBig{2_i}{1_{i-2}}{2_i}{2_{i+1}}{}{}{}{} \qquad \qquad \qquad
\basisBig{2_i}{1_{i-2}}{2_i}{2_{i+4}}{}{}{}{} \qquad \qquad \qquad
\basisBig{2_i}{2_{i+1}}{2_i}{2_{i+4}}{}{}{}{}
\]
and lead to the same equations. We can check that the value for the modulus of the cells obtained in~(\ref{E5squarecells}) satisfy these 5 equations.
We consider now the non-degenerated case. There are 12 non-degenerated Type~II  frames, consisting of the following three frames:\vspace{0.4mm}
\[
\basisBig{2_0}{2_4}{2_3}{2_1}{}{}{}{}  \qquad \qquad \qquad \basisBig{2_1}{2_5}{2_4}{2_2}{}{}{}{} \qquad \qquad \qquad \basisBig{2_2}{2_0}{2_5}{2_3}{}{}{}{}
\]
and those obtained by diagonal symmetry. They lead to the same equation:
\[
\frac{1}{3_q} |\mu|^2\mu^*(\nu_0+\nu_1) = 0.
\]
We can choose a solution by setting $\nu_1 = -\nu_0$.
Using gauge freedom (see our discussion in the next paragraph), we can make all the cells real. A solution to the cell system of $\mathcal{E}_5$ is therefore given by:
\begin{alignat*}{3}
& \tau_i =  \big(10+7\sqrt{2}\big)^{\frac{1}{4}},  \qquad&& \mu_i  = \left(5+\frac{7}{2}\sqrt{2}\right)^{\frac{1}{4}},  &\\
& \nu_0 =  \left(29+\frac{41}{2}\sqrt{2}\right)^{\frac{1}{4}},\qquad & &    \nu_1  = -\left(29+\frac{41}{2}\sqrt{2}\right)^{\frac{1}{4}}.&
\end{alignat*}

\begin{figure}[h!]
\centering
\scalebox{0.65}{\includegraphics[height=3.0in]{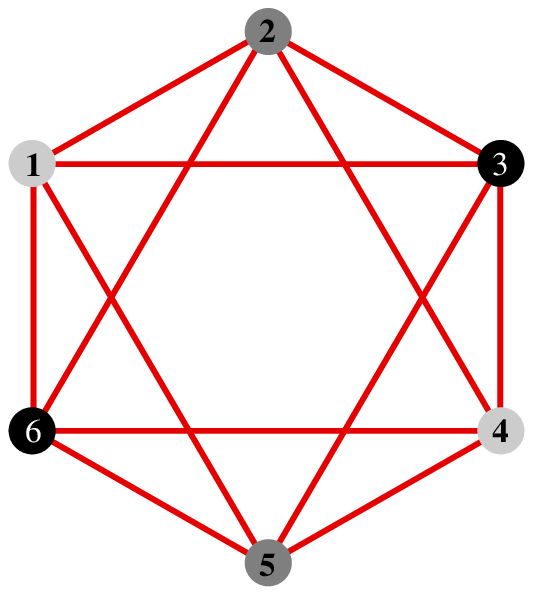}}
\setlength{\abovecaptionskip}{0.pt}
\caption{The inner octahedron of the ${\mathcal E}_5$ graph.}
\label{fig:octahedronE5}
\end{figure}

\paragraph{Fixing the phases using gauge freedom.}
The inner part of the fusion graphs, forgetting the six ``external triangles'' that have equal absolute cell values $\vert \tau \vert$, looks like a  plane projection of an octahedron, Fig.~\ref{fig:octahedronE5},
 seen from a direction orthogonal to one of  its triangular faces, so that we have\footnote{For the sake of this discussion, vertices are re-labelled as on Fig.~\ref{fig:octahedronE5}.} two large triangles (264 and 135), with the same absolute cell values $\vert \nu \vert$),  lying on top of each other and  six other triangles, with same absolute cell values $\vert \mu \vert$, on the sides. The following sequence of gauge choices (choice of a phase for an edge) make $7$ triangles of the octahedron real and positive, only one  (associated with triangle $246$) is still free at the end of the process.
The successive choices are given below, writing pairs for a chosen edge and the corresponding triangle made real and positive by that choice:
$(13, 135)$, $(12, 123)$, $(34, 345)$,  $(56, 561)$, $(24, 234)$, $(46, 456)$, $(62, 612)$. As we have seen, equations of Type~II  tell us (three times)  that the sum of the two pyramids whose union is the octahedron vanishes, so that
$\nu_0+\nu_1=0$, but gauge f\/ixing makes $\nu_0$ real positive, so $\nu_1$ is both real and negative.
Finally, we return to the six external triangles. For each of them, the edge common with the octahedron was already gauged f\/ixed.  This leaves us with two edges, still un-gauged, for each of them, so we can make all these cells $\tau$ real and positive. Remark:  the product of cells corresponding to the $8$ triangles of the octahedron is independent of the succession of gauge choices. In other words, the number $C= \mu^6 \nu_0  \nu_1= -  \frac{239 + 169 \sqrt{2}}{2}$ is an invariant.
An analogous discussion can be carried out for all examples discussed in this paper.

\subsection[The $\mathcal{E}_9$ case]{The $\boldsymbol{\mathcal{E}_9}$ case}\label{section4.5}

\begin{figure}[th]
\centerline{\scalebox{0.55}{\includegraphics{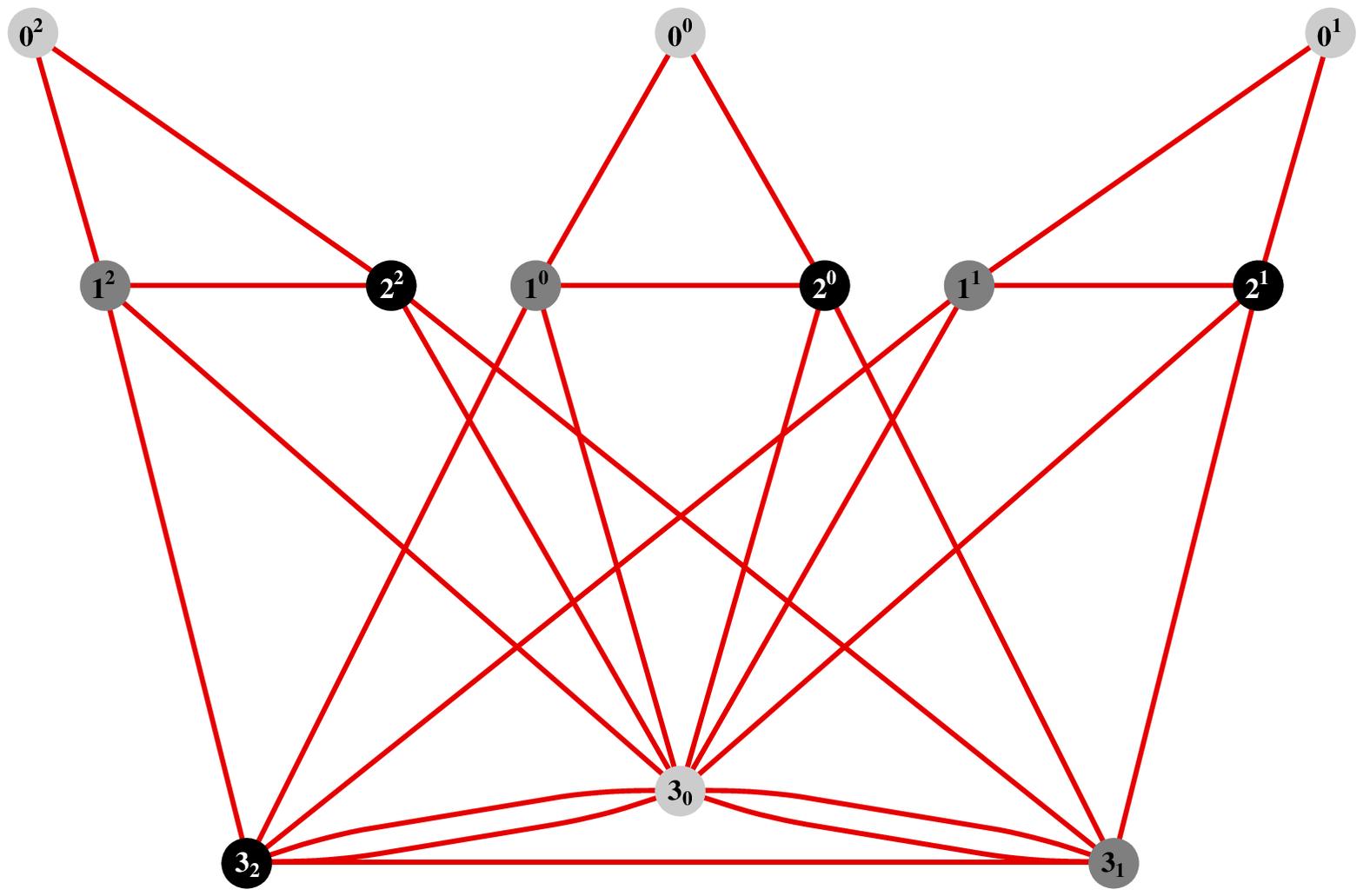}}}
\caption{The $\mathcal{E}_9$ fusion graph.}
\label{E9-graph}
\end{figure}

$\mathcal{E}_9$ is the graph of level $k=9$ and altitude $\kappa=k+3=12$, i.e.\ $q^{12}=-1$, displayed in Fig.~\ref{E9-graph}. It has 12 vertices that we label $0^j$, $1^j$, $2^j$, $3_0$, $3_1$, $3_2$, with $j=0,1,2$. The superscript $j$ ref\/lects the $\mathbb{Z}_3$ symmetry of the three wings of the graph. The quantum dimensions of the vertices are $[0^j]=1$, $[1^j]=[2^j]= 3_q = 1 + \sqrt{3}$ for $j=0,1,2$, $[3_0] = 3+2\sqrt{3}$ and $[3_1]=[3_2] = 3+\sqrt{3}$.

The number of oriented triangles (and therefore the number of cells) of the graph is ${\rm Tr}({G}_{\mathcal{E}_9})/3 $ $= 22$.
Six cells have only single edges between their vertices, they are denoted:
\[
a^j  = \Tribis{1^j}{2^j}{0^j}{}{}{}{\tau}   \qquad \qquad  b^j =\Tribis{1^j}{2^j}{3_0}{}{}{}{\tau}
\]
for $j=0,1,2$. The other cells involve the frames $\frameI{3_0}{3_1}{\alpha_1}{\alpha_2}$ and $\frameI{3_2}{3_0}{\beta_1}{\beta_2}$ with double edges, so that we have to introduce distinct edge labels. We denote as follows the 16 remaining cells:{\samepage
\begin{equation}
c^j_k  = \Tribis{3_0}{3_1}{2^j}{\alpha_k}{}{}{\tau}   \qquad \qquad
d^j_{\ell} = \Tribis{3_2}{3_0}{1^j}{\beta_{\ell}}{}{}{\tau}   \qquad \qquad
e_{k\ell} = \Tribis{3_1}{3_2}{3_0}{}{\beta_{\ell}}{\alpha_k}{\tau}
\label{16cellsE9}
\end{equation}
for $j = 0,1,2$ and $k,\ell = 1,2$.}

\paragraph{Fixing the modulus of the cells.}
We f\/irst consider Type~I  and Type~II  equations that only involve the modulus of the cells.
There are 24 Type~I  frames, but two of them involve double edges and the corresponding equations involve the phase of the cells. Among the~22 other Type~I  frames, only 13 of the corresponding equations are linearly independent. The three Type~I  frames
$\frameIsingle{0^j}{1^j}$ (and $\frameIsingle{2^j}{0^j}$) lead to $|a^j|^2 = 2_q[0^j][1^j] = \sqrt{2}(2+\sqrt{3})$, so that the three cells $a^j$ have equal modulus. From the three Type~I  frames $\frameIsingle{1^j}{2^j}$ we get $|a^j|^2 + |b^j|^2 = 2_q [1^j][2^j] = 2^{-1/2} { (1+ \sqrt{3})}^3$ so that the three cells $b^j$ have the same modulus, with $|b^j|^2 = 2^{-1/2} 3^{1/2} {(1+ \sqrt{3})}^2$.
The following 7 Type~I  frames give equations that involve the modulus of the 16 cells def\/ined in~(\ref{16cellsE9}):
\begin{alignat*}{3}
& \frameIsingle{1^j}{3_2}    \qquad&& |d^j_{1}|^2 + |d^j_{2}|^2  = 2_q[1^j][3_2] = 2^{-1/2} 3^{1/2} {(1+ \sqrt{3})}^3,&\\
& \frameIsingle{3_1}{2^j}     \qquad&& |c^j_{1}|^2 + |c^j_2|^2  = 2_q[3_1][2^j] =  2^{-1/2} 3^{1/2} {(1+ \sqrt{3})}^3, &\\
& \frameIsingle{3_1}{3_2}     \qquad \qquad&& |e_{11}|^2+|e_{12}|^2+|e_{21}|^2+|e_{22}|^2  = 2_q[3_1][3_2] = 2^{-1/2} 3 {(1+ \sqrt{3})}^3.&
\end{alignat*}
All other Type~I  frames (excepting the two involving double edges) lead to equivalent equations and give no further information.
Among the Type~II  frames leading to equations involving the modulus of the cells, only the 6 following ones lead to new equations:
\begin{equation*}
\begin{array}{@{}ccccc}
\basisBig{3_0}{3_1}{3_0}{3_1}{\alpha_1}{\alpha_1}{\alpha_1}{\alpha_1} &\qquad \qquad \qquad&
\basisBig{3_0}{3_1}{3_0}{3_1}{\alpha_2}{\alpha_2}{\alpha_2}{\alpha_2} &\qquad \qquad \qquad&
\basisBig{3_0}{3_1}{3_0}{3_1}{\alpha_1}{\alpha_1}{\alpha_2}{\alpha_2} \\
{ } & { } & { } & { } & { } \\
\basisBig{3_2}{3_0}{3_2}{3_0}{\beta_1}{\beta_1}{\beta_1}{\beta_1} &\qquad \qquad \qquad&
\basisBig{3_2}{3_0}{3_2}{3_0}{\beta_2}{\beta_2}{\beta_2}{\beta_2} &\qquad \qquad \qquad&
\basisBig{3_2}{3_0}{3_2}{3_0}{\beta_1}{\beta_2}{\beta_2}{\beta_1}
\end{array}
\end{equation*}
At this point we have $7+6=13$ linearly independent equations involving the 16 undetermined parameters that are the modulus of the cells def\/ined in~(\ref{16cellsE9}).

\paragraph{Looking for a symmetric solution.}
The three cells $a^j$ have the same modulus, as well as the three cells $b^j$; this ref\/lects the $\mathbb{Z}_3$ symmetry of the graph. We look for a solution of the cell system of $\mathcal{E}_9$ ref\/lecting this symmetry, and we def\/ine the numbers $|\tau|^2$ and $|\mu|^2$ by setting:
\begin{equation}
|c^j_1|^2 = |d^j_1|^2 = |\tau|^2, \qquad |c^j_2|^2 = |d^j_2|^2 = |\mu|^2   \qquad \textrm{for} \quad j=0,1,2.
\label{symsolE9}
\end{equation}
With this choice, the f\/irst two  equations of Type~I  become equivalent, and we end with only two linearly independent  equations of Type~I.  The equations corresponding to the Type~II  frames def\/ined previously read:
\begin{gather*}
\frac{3}{[2^j]}|\tau|^4 + \frac{1}{[3_2]} (|e_{11}|^4 + 2|e_{11}|^2|e_{21}|^2 + |e_{21}|^4)  =  [3_0][3_0][3_1]+[3_0][3_1][3_1], \\
\frac{3}{[2^j]}|\mu|^4 + \frac{1}{[3_2]} (|e_{12}|^4 + 2|e_{12}|^2|e_{22}|^2 + |e_{22}|^4)  =  [3_0][3_0][3_1]+[3_0][3_1][3_1], \\
\frac{3}{[2^j]}|\tau|^2|\mu|^2 + \frac{1}{[3_2]} \left(|e_{11}|^2 |e_{12}|^2 + |e_{12}|^2 |e_{21}|^2 + |e_{11}|^2|e_{22}|^2 + |e_{21}|^2|e_{22}|^2\right)  =  [3_0][3_0][3_1], \\
\frac{3}{[1^j]}|\tau|^4 + \frac{1}{[3_1]} (|e_{11}|^4 + 2|e_{11}|^2|e_{12}|^2 + |e_{12}|^4)  =  [3_0][3_0][3_1]+[3_0][3_2][3_2], \\
\frac{3}{[1^j]}|\mu|^4 + \frac{1}{[3_1]} (|e_{21}|^4 + 2|e_{21}|^2|e_{22}|^2 + |e_{22}|^4)  =  [3_0][3_0][3_1]+[3_0][3_2][3_2] ,\\
\frac{3}{[1^j]}|\tau|^2|\mu|^2 + \frac{1}{[3_1]} \left(|e_{11}|^2 |e_{21}|^2 + |e_{12}|^2 |e_{21}|^2 + |e_{11}|^2|e_{22}|^2 + |e_{12}|^2|e_{22}|^2\right)  = [3_0][3_0][3_2].
\end{gather*}
These equations impose the symmetry condition $|e_{12}|=|e_{21}|$.  The last three equations of Type~II   are then equivalent to the f\/irst three. Also the third  equation of Type~II can be written in terms of the other two and of the two  equations of Type~I. At this point, with the symmetric condition~(\ref{symsolE9}), we have 5 parameters, $|\tau|$, $|\mu|$, $|e_{11}|$, $|e_{12}|$ and $|e_{22}|$ but only four linearly independent equations, so we still have one unknown parameter, say $|e_{22}|$. Solving this set of (quadratic and quartic) equations for the modulus of the cells gives two solutions depending on this parameter:
\begin{gather}
\left\lbrace \begin{array}{l}
|\tau|^2  =   r_+, \\
|\mu|^2   =  r_-, \\
|e_{11}|^2  =  \rho_- - \rho_+ + |e_{22}|^2, \\
|e_{12}|^2 = |e_{21}|^2  =  \rho_+ - |e_{22}|^2,
\end{array} \right.
\qquad
\left\lbrace \begin{array}{l}
|\tau|^2  =  r_-, \\
|\mu|^2  =  r_+, \\
|e_{11}|^2  =  \rho_+ - \rho_- + |e_{22}|^2, \\
|e_{12}|^2 = |e_{21}|^2  =  \rho_- - |e_{22}|^2,
\end{array} \right.
\label{solmid}
\end{gather}
where $ r_{\pm} = \sqrt{78+45\sqrt{3}} \pm \sqrt{12+7\sqrt{3}} =  2^{-3/2} 3^{1/2} (1+ \sqrt{3})^3 \pm    2^{-1} 3^{1/4} (1+ \sqrt{3})^2 $ and $\rho_{\pm} =  3 \big((2+\sqrt{3})^{\frac{3}{2}} \pm  \sqrt{12+7\sqrt{3}} \big) =  2^{-3/2} 3 (1+ \sqrt{3})^3 \pm    2^{-1} 3^{5/4} (1+ \sqrt{3})^2.$
In order to determine the last parameter $|e_{22}|$, we consider the following frames:
\begin{equation}
\basisBig{3_0}{3_1}{3_0}{3_1}{\alpha_1}{\alpha_2}{\alpha_2}{\alpha_1}   \qquad \qquad \qquad
\basisBig{3_1}{3_2}{3_1}{3_2}{}{}{}{}
\label{e22fix}
\end{equation}
As usual, the equations associated with these frames are obtained by summing over all vertices and edges building the pyramid of equation~(\ref{typeIIeq}). For the frame on the r.h.s.\ of~(\ref{e22fix}), we can only sum over the vertex $3_0$ but all edges connecting the frame to $3_0$ are double, so that the summation f\/inally consists of $2^4 = 16$ terms! The Type~II  equations corresponding to these two frames read:
 \begin{gather*}
\frac{3}{[2^j]}|\tau|^2|\mu|^2 + \frac{1}{[3_2]}\big(|e_{11}|^2|e_{12}|^2+|e_{21}|^2|e_{22}|^2+
e_{11}\overline{e_{12}}e_{22}\overline{e_{21}}+
\overline{e_{11}}e_{12}\overline{e_{22}}e_{21}\big)  =  [3_0][3_1][3_1], \\
\frac{1}{[3_0]} \big( |e_{11}|^4 + |e_{12}|^4 + |e_{21}|^4 + |e_{22}|^4 + 2 \big(|e_{11}|^2+|e_{22}|^2\big)\big(|e_{12}|^2+|e_{21}|^2\big)\\
\qquad{} + 2 (e_{11}\overline{e_{12}}e_{22}\overline{e_{21}}+
\overline{e_{11}}e_{12}\overline{e_{22}}e_{21})\big)  =  [3_1][3_1][3_2]+[3_1][3_2][3_2].
\end{gather*}
We can combine them in order to get rid of the factors involving phases, and we get a new equation involving only the modulus of the cells. This allows to f\/ix the last parameter $|e_{22}|$. For each of the two solutions in~(\ref{solmid}), we obtain two possible values for $|e_{22}|$, namely:
\[
\left\lbrace \begin{array}{l}
|e_{22}|^2  =  0, \\
|e_{22}|^2  =  \rho_+,
\end{array} \right.
 \qquad
\left\lbrace \begin{array}{l}
|e_{22}|^2 = 0, \\
|e_{22}|^2  =  \rho_-.
\end{array} \right.
\]
So one would obtain four solutions for the modulus of the $e_{ij}$, with the given symmetry conditions, but actually only three remain: one of them is rejected because it leads to a negative value for  $|e_{11}|$. The two solutions with $|e_{22}|^2 \neq 0$ lead to $e_{12} = e_{21} = 0$, a condition making the equations involving phases more simple, we make this choice in the following.

\paragraph{Fixing the phases of the cells.}
There are no equations involving the phases of  $a_j$ and $b_j$ and we can make their values real by a gauge choice. We denote the other cells as follows:
\[
c_k^j= |c_k^j|e^{i\theta_{jk}}, \qquad d^j_{\ell}= |d^j_{\ell}|e^{i\psi_{j{\ell}}}, \qquad e_{k \ell} = |e_{k\ell}|e^{i\phi_{k\ell}}
\]
for $j = 0,1,2$ and $k,\ell = 1,2$. The Type~I  frames involving double edges and the corresponding equations are:
\begin{alignat*}{3}
& \frameI{3_0}{3_1}{\alpha_1}{\alpha_2} \qquad \qquad & &  c_1^1\overline{c_2^1}+c_1^2\overline{c_2^2}+c_1^3\overline{c_2^3} +
e_{11}\overline{e_{12}}+e_{21}\overline{e_{22}} =  0, & \\[1ex]
& \frameI{3_2}{3_0}{\beta_1}{\beta_2}  \qquad \qquad & &  d_1^1\overline{d_2^1}+d_1^2\overline{d_2^2}+d_1^3\overline{d_2^3} +
e_{11}\overline{e_{21}}+e_{12}\overline{e_{22}}=  0. &
\end{alignat*}
Looking for a solution with the symmetric condition (\ref{symsolE9}) together with the choice $|e_{12}|{=}|e_{21}|{=}0$,
the above two equations lead to the following constraints on phases:
\begin{gather}
\label{phase1} e^{i(\theta_{11} - \theta_{12})} + e^{i(\theta_{21} - \theta_{22})} + e^{i(\theta_{31} - \theta_{32})}  =  0, \\
\label{phase2} e^{i(\psi_{11} - \psi_{12})} + e^{i(\psi_{21} - \psi_{22})} + e^{i(\psi_{31} - \psi_{32})}  =  0.
\end{gather}
We now consider the following Type~II  frames and the corresponding equations:
\begin{gather*}
\basisBig{0_2^j}{3_0}{3_2}{3_0}{}{}{\beta_2}{\beta_1} \\  \frac{1}{[0_1^j]}(|b^j|^2 d_1^j\overline{d_2^j}) + \frac{1}{[3_1]}\big(
|c_1^j|^2 e_{11}\overline{e_{21}} + |c_2^j|^2 e_{12}\overline{e_{22}} + c_1^j\overline{c_2^j}e_{12}\overline{e_{21}} + \overline{c_1^j}c_2^je_{11}\overline{e_{22}} \big)=0.
\end{gather*}
With the chosen solution $|e_{12}|=|e_{21}|=0$, and from the fact that $\frac{1}{[0_1^j]} |b^j|^2 = \frac{1}{[3_1]}|e_{11}||e_{22}|$, the equations read:
\begin{gather*}
e^{i(\psi_{j1}-\psi_{j2})}+e^{i(\theta_{j2}-\theta_{j1})}e^{i(\phi_{11}-\phi_{22})} =0 \qquad  \textrm{for} \quad j=0,1,2.
\end{gather*}
Setting $\phi_{11}-\phi_{22}=\pi \textrm{ (mod } 2\pi)$, we get $e^{i(\psi_{j1}-\psi_{j2})} = e^{i(\theta_{j2}-\theta_{j1})}$ and the two equations~(\ref{phase1}) and~(\ref{phase2}) become equivalent. A solution of~(\ref{phase1}) is given by $\theta_{j1}=\frac{2 \pi j}{3}$ and $\theta_{j2} = -\frac{2 \pi j}{3}$.

\paragraph{A solution to the cell-system of $\boldsymbol{\mathcal{E}_9}$.}
With  $|\tau|^2 = r_+$, $|\mu|^2 = r_-$, $|e_{11}|^2 = \rho_-$ and $|e_{22}|^2 = \rho_+$, a solution is obtained by setting $\phi_{11} = 0$, $\phi_{22} = \pi$, $\psi_{j1} = \theta_{j2}$ and $\psi_{j2} = \theta_{j1}$. One can check that these cells satisfy all Type~I  and Type~II  equations and therefore give a solution to the cell-system of~$\mathcal{E}_9$. Summarizing results, our solution\footnote{Another gauge nonequivalent solution is discussed in the next paragraph.} is given by:
\begin{alignat*}{3}
& a^j  =  2^{-1/4} (1 + \sqrt 3),  \qquad   & & b^j  =  2^{-1/4}    3^{1/4}     (1 + \sqrt 3),&  \\
& c_1^j  =  \sqrt{r_+}   e^{\frac{2 i \pi j}{3}}, \qquad  &&  c_2^j = \sqrt{r_-}  e^{\frac{ - 2 i \pi j}{3}}, & \\
& d_1^j  =  \sqrt{r_+}   e^{\frac{ - 2 i \pi j}{3}},    \qquad && d_2^j  =  \sqrt{r_-}   e^{\frac{2 i \pi j}{3}}, &  \\
& e_{11}  =  \sqrt{\rho_-},   \qquad && e_{22}  =  - \sqrt{\rho_+}, & \\
& e_{12}  =   0,    \qquad && e_{21}  =  0&
\end{alignat*}
for $j=0,1,2$, with
\[
r_{\pm} = \sqrt{78+45\sqrt{3}} \pm \sqrt{12+7\sqrt{3}} \qquad \textrm{and} \qquad
\rho_{\pm} = 3(2+\sqrt{3})^{\frac{3}{2}} \pm 3 \sqrt{12+7\sqrt{3}}
\]
that can be written
\begin{gather*}
r_{\pm} =  2^{-3/2} 3^{1/2} (1+ \sqrt{3})^3 \pm    2^{-1} 3^{1/4} (1+ \sqrt{3})^2  \qquad \textrm{and } \\
\rho_{\pm} = 2^{-3/2} 3 (1+ \sqrt{3})^3 \pm    2^{-1} 3^{5/4} (1+ \sqrt{3})^2.
\end{gather*}

\paragraph{Gauge freedom and invariants for $\boldsymbol{\mathcal{E}_9}$.}
Because of gauge freedom, the given solution for the moduli and phases of $c_k^j$, $d_{\ell}^j$ and $e_{k\ell}$ is not unique (the moduli of cells $a^j$ and $b^j$ are obviously invariant).
Moreover the existence of two double lines in the graph makes the discussion of gauge freedom slightly more involved than in the other examples discussed in this paper.

Def\/ine the matrix $M = \lbrace \lbrace e_{11}, e_{12} \rbrace, \lbrace e_{21}, e_{22} \rbrace \rbrace$.
Gauge freedom allows one to modify the values $e_{ij}  \mapsto {e_{ij}}^\prime$ of these four cells by introducing two arbitrary $2 \times 2$ unitary matrices $U$, $V$ (for the two double lines)  together with a phase $\exp(i \alpha)$, and setting
$M^\prime =  \exp(i \alpha)  U . M . V$.  This implies  ${M^\prime}^\dag . M^\prime = V^\dag . M^\dag . U^\dag . U . M . V = V^\dag . M^\dag . M . V  =  V^{-1} . M^\dag . M . V$
The trace and determinant of the matrix  $M^\dag . M$ are therefore gauge invariant quantities. They read
\begin{gather*}
{\rm Det} (M^\dag . M)  =  {\rm Det}(M) \overline { {\rm Det}(M)} =  (-e_{12} e_{21} + e_{11} e_{22}) \overline{[-e_{12}e_{21} + e_{11} e_{22}]} = 2^{-1} 3^2  (1 +  \sqrt 3 )^4, \\
{\rm Tr} (M^\dag . M)  =  e_{11} \overline{e_{11}} + e_{12} \overline{e_{12}} + e_{21} \overline{e_{21}} + e_{22} \overline{e_{22}} =  2^{-1/2}   3   (1+ \sqrt 3)^3.
\end{gather*}

Def\/ine $c^0=(c_1^0, c_2^0)$,  $c^1=(c_1^1, c_2^1)$, $c^2=(c_1^2, c_2^2)$.
Gauge freedom gives $c^j \mapsto {c^{j}}^\prime = U.c^j   e^{i \phi_j}   e^{i \psi_j}$ where $U$ is a $2 \times 2$ unitary matrix (the same for all~$j$ since these pairs of cells share the same double-edge) and where $\phi_j$, $\psi_j$ are phases.
One sees immediately that the following quantities are gauge invariant:
\begin{gather*}
{c^0}^\dag . c^0  =  {c^1}^\dag . c^1={c^2}^\dag . c^2 = r_{-} + r_{+}= 2^{-1/2} 3^{1/2}  (1 +  \sqrt 3 )^3,  \\
({c^0}^\dag . c^1)({c^1}^\dag . c^0) =  ({c^1}^\dag . c^2)({c^2}^\dag . c^1)= ({c^2}^\dag . c^0)({c^0}^\dag . c^2)=r_{-}^2 +r_{+}^2 - r_{-}  r_+ = {2}^{-1}{3} \big(1+\sqrt{3}\big)^5,\\
({c^0}^\dag . c^1)({c^1}^\dag . c^2) ({c^2}^\dag . c^0)  =  (r_{-} -  \exp(i \pi/3) r_{+})^3 = {2}^{-3/2} {3}^{3/2}  \big(1+\sqrt{3}\big)^{15/2}  e^{i \arctan \big(\frac{3^{3/4}  (1+\sqrt{3} )}{\sqrt{2}}\big)}.
\end{gather*}
The last invariant is complex and goes to its complex conjugate (also an invariant) under an odd permutations of the branches  0, 1, 2. This shows that we obtain actually two gauge nonequivalent solutions for the cell system,  their
existence ref\/lects the symmetry of the given graph (keep $j=0$ but permute the lower indices $j=1$, $j=2$ in the solution for $c$ and $d$ given previously).

We can also introduce vectors $d^0=(d_1^0, d_2^0)$,  $d^1=(d_1^1, d_2^1)$, $d^2=(d_1^2, d_2^2)$. The discussion is entirely similar for the coef\/f\/icients $d_k^j$ (introduce a $2\times 2$ unitary matrix $V$ and $U(1)$ phases),  and leads to analogous gauge invariant quantities with the same values as above.

Together with the modules of the cells $a^j$ and $b^j$, for $j=0,1,2$,  the above invariants provide a~tool that enable one to check compatibility of seemingly dif\/ferent solutions.
In order to illustrate this, we mention:

 $\bullet$ A solution\footnote{The other (gauge nonequivalent) solution is obtained by $\lambda_{\pm}  \leftrightarrow \lambda_{\mp} $ and $e_{21}  \leftrightarrow  -e_{12}$.} given by~\cite{EvansSU3Cells} (the authors did not directly calculated the cells for ${\mathcal E}_9$, instead they determined those of its exceptional module ${\mathcal M}_9$, an easier task, and then use the fact that the former is a $Z_3$ orbifold of the later):
\begin{alignat*}{3}
& a^j =  2^{-1/4} (1 + \sqrt 3),  \qquad && b^j  =  2^{-1/4}    3^{1/4}    (1 + \sqrt 3), & \\
& c_1^j =  \lambda_{+}    e^{\frac{2 i \pi (j-1)}{3}},  \qquad  &&  c_2^j =\lambda_{-}   e^{\frac{ - 2 i \pi (j-1)}{3}}, & \\
& d_1^j = \lambda_{-}    e^{\frac{ - 2 i \pi (j-1)}{3}}, \qquad && d_2^j =\lambda_{+}  e^{\frac{2 i \pi (j-1)}{3}}, & \\
& e_{11}  =  0, \qquad  && e_{22}  =   0, &  \\
& e_{12}  =    - \frac{[4]_q}{ \sqrt {[2]_q}}    \sqrt{[2]_q    [2]_q + \sqrt{ [2]_q  [4]_q}}, \qquad  &&
e_{21}  =    \frac{[4]_q}{ \sqrt {[2]_q}}    \sqrt{[2]_q    [2]_q - \sqrt{ [2]_q \,[4]_q}}, &
\end{alignat*}
with  $ \lambda_{\pm} = \sqrt{  [2]_q} \sqrt{  [2]_q [4]_q \pm \sqrt{ [2]_q [4]_q}}$.

$\bullet$ A solution given (without proof) by \cite{Oc:private}:
\begin{alignat*}{3}
& a^j = {2^{-1/4}} {\big(1+\sqrt{3}\big)}, \qquad  && b^j = {2^{-1/4} 3^{1/4}}\big(1+\sqrt{3}\big), &\\
& c_2^0 =  0,  \qquad  &&  c_1^0 =  {2^{-1/4} 3^{1/4}} \big(1+\sqrt{3}\big)^{3/2}, &\\
& c_2^1 =  {2^{-1/4}} {3^{1/2}} { \big(1+\sqrt{3}\big)}, \qquad  &&
c_1^1  =  {2^{-1/4} 3^{1/4}} \big(1+\sqrt{3}\big), &  \\
& c_2^2 =   {2^{-1/4} 3^{1/2}} \big(1+\sqrt{3}\big)  \left(\frac{1}{2} \big(  1-\sqrt{3}\big) - i  {2^{-1/2} 3^{1/4}} \right),
 \qquad  &&
c_1^2  =  {2^{-1/4} 3^{1/4}}  \big(1+\sqrt{3}\big), &\\
& d_2^0 =  0, \qquad   && d_1^0  =  c_1^0, & \\
& d_2^1 =  c_2^1, \qquad   && d_1^1  =  c_1^1, &   \\
&  d_2^2 =  c_2^2, \qquad  && d_1^2  = \overline {c_1^2},&  \\
& e_{12} = 2^{-5/4}  \big(3+\sqrt{3}\big)
\big(
 \big(1-\sqrt{3}\big) -
i {2^{1/2}{3}^{1/4}}
\big),
   \qquad && e_{21} ={2^{-1/4} {{3^{1/2}} \big(1+\sqrt{3}\big)} }, &
\\
& e_{22} =    {2^{-5/4} } {3^{1/2} \big(1+\sqrt{3}\big)^{1/2}} \Big(- \sqrt{3-\sqrt{3}}+ i
   \sqrt{1+\sqrt{3}}\Big), \qquad &&  e_{11}  = 0. &
\end{alignat*}
These solutions look very dif\/ferent and obey all equations of Types~I and~II, but one can check that the gauge invariants are equal\footnote{One can choose unitaries $U$ and $V$  making equal the scalars ${d^j}^\dag . c^j$, for $j=0,1,2$. This equality  holds for the three solutions given previously.}.

\subsection[The $\mathcal{E}_{21}$ case]{The $\boldsymbol{\mathcal{E}_{21}}$ case}\label{section4.6}

\begin{figure}[h!]\label{E21}
\centerline{\includegraphics[scale=0.50]{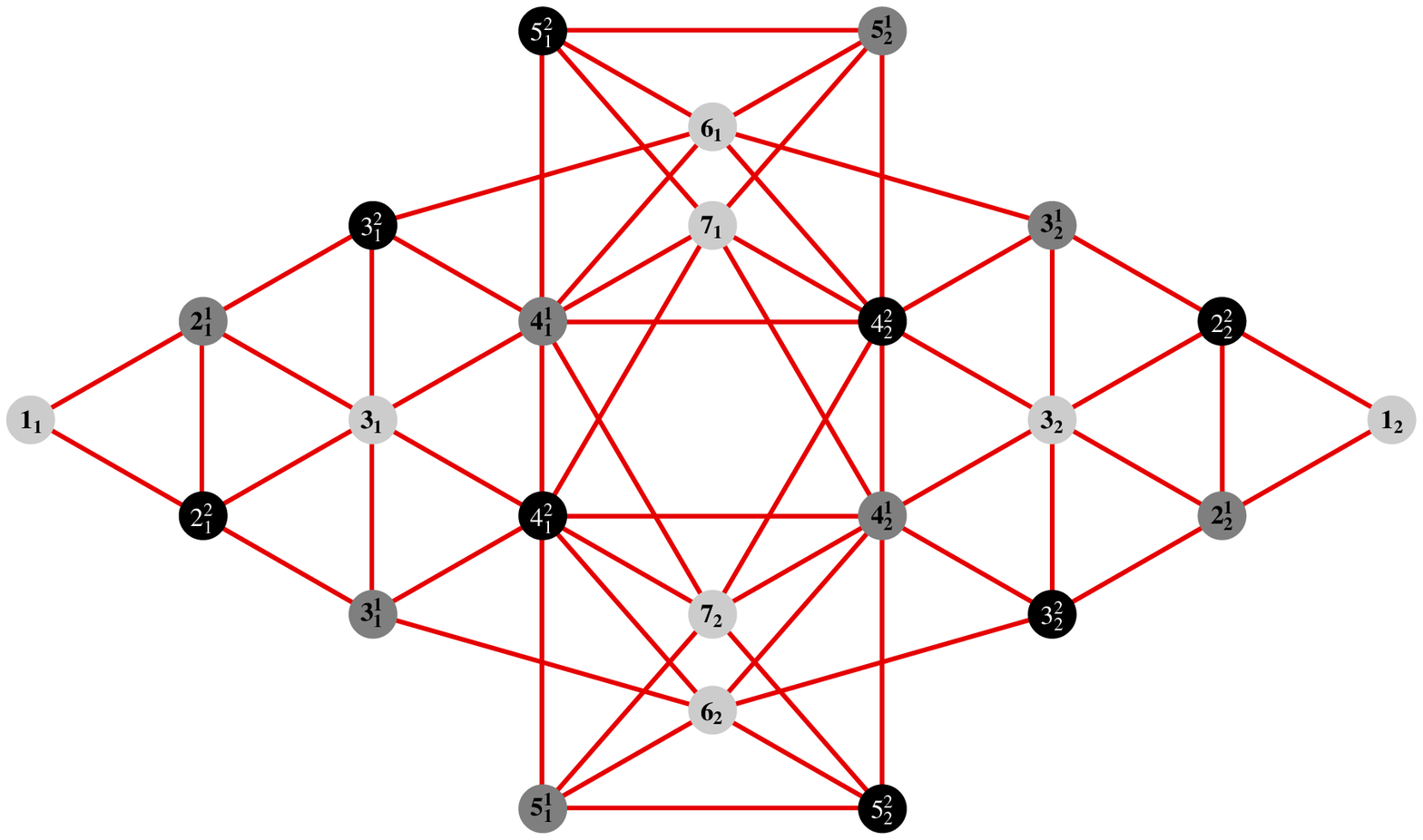}}
\caption{The $\mathcal{E}_{21}$ fusion graph.}
\label{graphE21}
\end{figure}

The $\mathcal{E}_{21}$ graph is the graph of level $k=21$ and altitude $\kappa=k+3=24$, i.e.\ $q^{24}=-1$, displayed in Fig.~\ref{graphE21}. It is a graph with 24 vertices that we divide into $7$ slices: $1_k$, $2_k^i$, $3_k$, $3_k^i$, $4^i_{k}$, $5^i_k$, $6_k$ and $7_k$, where $i=1,2$ and $k=1,2$. For a given vertex $a^j_k$, the superscript  $j=\{0,1,2\}$ gives the triality of the vertex, with the convention that triality 0 is omitted. The subscript $k=\{1,2\}$ refers to the $\mathbb{Z}_2$ symmetry of the graph with respect to the central point of the graph (it is the composition of the left-right and up-down symmetries). Notice that this $\mathbb{Z}_2$ symmetry preserves the triality of the vertices: $a^j_1 \leftrightarrow a^j_2$.
We introduce the function $g[a,b,c,d]=(a+b\sqrt 2 + c \sqrt 3 + d \sqrt 6)/2$.
The quantum dimensions of vertices are: $[1_k]= 1_q = 1=g[2,0,0,0]$, $[2^i_k] = 3_q=g[2,1,0,1]$, $[3_k]  = 2_q 4_q=g[4,2,2,2]$, $[3^i_k]  = \frac{3_q4_q}{2_q}=g[4,1,2,1]$, $[4^i_k]  =3_q5_q=g[6,5,4,3]$, $[5^i_k]  = 9_q=g[4,3,2,1]$, $[6_k]  = \frac{4_q5_q}{2_q}=g[4,4,2,2]$ and $[7_k]  =\frac{4_q7_q}{2_q}=g[6,4,4,2]$, for $i=1,2$ and $k=1,2$.

The number of oriented triangles of the graph is ${\rm Tr}(G)^3/3=40$. Due to symmetry considerations (left-right and up-down), there are actually only 14 dif\/ferent square moduli for the cells. The types of triangles are labelled as follows:
\begin{alignat*}{3}
& \text{two  of type  } \{ 1 \rightarrow 2 \rightarrow 2 \}  \text{ label } \alpha_1,   \qquad &&
\text{two of type  } \{5 \rightarrow 5 \rightarrow 7 \}    \text{ label } \nu_2, &  \\
& \text{four  of type  } \{2 \rightarrow 3 \rightarrow 3\}   \text{ label } \alpha_2,\qquad   &   &
\text{four  of type  } \{4 \rightarrow 5 \rightarrow 6 \}  \text{ label } \mu_1, &  \\
& \text{two  of type  } \{3 \rightarrow 4 \rightarrow 4\}    \text{ label } \alpha_3 , \qquad &   &
\text{four of type  } \{4 \rightarrow 5 \rightarrow 7 \}   \text{ label } \mu_2, &  \\
& \text{two  of type  } \{2 \rightarrow 2 \rightarrow 3 \}   \text{ label } \gamma_2,\qquad  &   &
\text{two  of type  } \{4_k \rightarrow 4_{k+1} \rightarrow 6 \}  \text{ label } \sigma_1, &  \\
& \text{four  of type  } \{ 3 \rightarrow 3 \rightarrow 4\}  \text{ label } \gamma_3 \qquad &   &
\text{two of type  } \{4_k \rightarrow 4_{k+1} \rightarrow 7 \}   \text{ label } \sigma_2 \text{ (small)}, &  \\
& \text{four  of type  } \{3 \rightarrow  4 \rightarrow 6 \}    \text{ label } \beta \qquad &   &
\text{two  of type  } \{4_k \rightarrow 4_{k+1} \rightarrow 7 \}  \text{ label } \rho  \text{ (large)}, & \\
& \text{two  of type  } \{5 \rightarrow 5 \rightarrow 6 \}      \text{ label } \nu_1,\qquad   &   &
\text{four  of type  } \{4_k \rightarrow 4_{k} \rightarrow 7 \} \text{ label } \lambda. &
\end{alignat*}
Their square moduli are displayed below, where we omit indices when no confusion is possible,
\begin{alignat*}{4}
& |\alpha_1|^2  =  \Trimod{2}{2}{1}  \qquad&&  |\alpha_2|^2  =  \Trimod{3}{3}{2}  \qquad&&  |\alpha_3|^2  = \Trimod{4}{4}{3}& \\
&&& |\gamma_2|^2  =  \Trimod{2}{2}{3}  \qquad& & |\gamma_3|^2  = \Trimod{3}{3}{4} & \\
& |\nu_1|^2  =  \Trimod{5}{5}{6}  \qquad&& |\mu_1|^2  =  \Trimod{4}{5}{6}  \qquad&& |\sigma_1|^2  =  \Trimod{4_k}{4_{k+1}}{6} &\\
& |\nu_2|^2  =  \Trimod{5}{5}{7}  \qquad&& |\mu_2|^2  =  \Trimod{4}{5}{7}  \qquad&& |\sigma_2|^2  =  \Trimod{4_k}{4_{k+1}}{7} (\textrm{small})& \\
& |\beta|^2  =  \Trimod{3}{4}{6}  \qquad&& |\lambda|^2  =  \Trimod{4_k\;}{4_k}{7} \qquad&& |\rho|^2  =  \Trimod{4_k\;}{4_{k+1}}{7} (\textrm{big})&
\end{alignat*}

The number of Type~I equations is $\sum_{a,b}(G)_{ab}=60$. Symmetries of the graph greatly reduce the number of independent equations. For example, the coherence equations for the two Type~I  frames $\frameIsingle{1_k}{2^1_k}$ show that the square modulus of the two cells $\{1_1\rightarrow 2_1^1 \rightarrow 2_1^2\}$ and $\{1_2 \rightarrow 2_2^1 \rightarrow 2_2^2\}$ are the same, and equal to $|\alpha_1|^2 = 2_q[1_k][2_k^i] = 2_q3_q$. This justif\/ies the symmetry conditions introduced previously. The same argument applies for the other frames. In the following, in order to ease the notation, we shall omit superscripts and subscripts when no confusion is possible.
We write below the equations for Type~I frames in terms of quantum dimensions and $q$-numbers. A subset of Type~I  equations allows us to f\/ix step-by-step the values of the modulus square of the cells:
\begin{alignat}{3}
& \frameIno{1}{2}   \qquad && |\alpha_1|^2 = 2_q[1_k][2^i_k] =  2_q3_q, & \nonumber\\
& \frameIno{2}{3^i}   \qquad && |\alpha_2|^2  = 2_q[2^i_k][3^i_k] =  3_q^24_q, & \nonumber\\
& \frameIno{3^i}{6} \qquad && |\beta|^2 = 2_q[3^i_k][6_k] =  \frac{3_q4_q^25_q}{2_q}, & \nonumber\\
& \frameIno{2}{2}  \qquad && |\gamma_2|^2+|\alpha_1|^2 = 2_q[2^i_k][2^i_k] =  2_q3_q^2, & \nonumber\\
& \frameIno{3^i}{3}  \qquad && |\alpha_2|^2+|\gamma_3|^2 = 2_q[3_k^i][3_k] =  2_q3_q4_q^2,& \nonumber\\
& \frameIno{3}{4}  \qquad && |\alpha_3|^2+|\gamma_3|^2 =  2_q[3_k][4^i_k] =  2_q^23_q4_q5_q, & \nonumber\\
& \frameIno{4_k}{4_k}  \qquad && |\alpha_3|^2+ 2 |\lambda|^2 = 2_q[4^i_k][4^i_k] =  2_q3_q^25_q^2. & \label{E21typeI-1}
\end{alignat}

The following Type I frames give relations between the square moduli of the cells:
\begin{alignat}{3}
& \frameIno{4}{5}    \qquad  &&  |\mu_1|^2+|\mu_2|^2 = 2_q [4^i_k][5^i_k] = 2_q 3_q 5_q 9_q, &  \nonumber\\
& \frameIno{4_k}{4_{k+1}}    \qquad &&  |\sigma_1|^2+|\sigma_2|^2+|\rho|^2 = 2_q[4^i_k][4^i_k] = 2_q 3_q^25_q^2, &  \nonumber\\
& \frameIno{5}{6}    \qquad  && |\mu_1|^2+|\nu_1|^2 = 2_q[5^i_k][6_k] =  4_q5_q9_q , & \nonumber\\
& \frameIno{5}{7}   \qquad  && |\mu_2|^2+|\nu_2|^2 = 2_q[5^i_k][7_k] =  4_q7_q9_q. & \label{E21-typeI-2}
\end{alignat}

In order to f\/ix the remaining values, we have to consider also Type~II  equations. We use three Type~II  degenerated frames and the corresponding equations
\begin{alignat}{3}
& \basisBigLarge{3_1^1}{4_1^2}{5_1^1}{4_1^2}{}{}{}{} \qquad \qquad && |\beta|^2 |\mu_1|^2  =  [6_2][3^1_1] [5^1_1] [4_1^2] =  \frac{3_q^2 4_q^2 5_q^2 9_q}{2_q^2} , & \nonumber\\[1ex]
& \basisBigLarge{3_1^2}{4_1^1}{4_2^2}{4_1^1}{}{}{}{} \qquad \qquad && |\beta|^2 |\sigma_1|^2 = [6_1][3^2_1][4_2^2][4^1_1] =
\displaystyle \frac{3_q^3 4_q^2 5_q^3}{2_q^2}, &  \label{TypeIIE21}\\[1ex]
& \basisBigLarge{6_2}{4_2^1}{7_1}{4_2^1}{}{}{}{} \qquad \qquad && |\sigma_1|^2 |\rho|^2 =  [4^2_1][6_2][7_1][4^1_2] =   \frac{3_q^2 4_q^2 5_q^37_q}{2_q^2}.
& \nonumber
\end{alignat}

\paragraph{The solution for the square moduli.}
The values of $\alpha_1$, $\alpha_2$ and $\beta$ are read immediately from equations (\ref{E21typeI-1}).
From them we calculate $\gamma_2$, $\gamma_3$, $\alpha_3$ and $\lambda$ as follows:
\begin{alignat*}{3}
& |\gamma_2|^2  =  2_q3_q^2 - |\alpha_1|^2 = 3_q4_q ,   \qquad &&  |\gamma_3|^2  =  2_q3_q4_q^2 - |\alpha_2|^2 = 3_q4_q5_q, & \\
& |\alpha_3|^2  =  2_q^23_q4_q5_q - |\gamma_3|^2 = 3_q^24_q5_q , \qquad && | \lambda|^2  =  \frac{1}{2}(2_q3_q^25_q^2-|\alpha_2|^2)=\frac{1}{2}3_q^25_q6_q. &
\end{alignat*}
Here we have used the following q-numbers identities to simplify the expressions (remember that $q^{24}=-1$): $2_q 2_q = 1+3_q$,  $2_q3_q=2_q+4_q$, $2_q4_q=3_q+5_q$ and $2_q5_q=4_q+6_q$.
We then solve equations (\ref{TypeIIE21}) to determine the square moduli of $\mu_1$, $\sigma_1$  and $\rho$:
\begin{gather*}
|\mu_1|^2  =  \frac{1}{|\beta|^2}\frac{3_q^24_q^25_q^29_q}{2_q^2}=\frac{3_q5_q9_q}{2_q},   \qquad
|\sigma_1|^2  =  \frac{1}{|\beta|^2}\frac{3_q^34_q^25_q^3}{2_q^2}=\frac{3_q^25_q^2}{2_q}, \\
|\rho|^2  =  \frac{1}{|\sigma_1|^2}\frac{3_q^24_q^25_q^37_q}{2_q^2}=\frac{4_q^25_q7_q}{2_q}.
\end{gather*}
Finally, the square moduli of the remaining cells are obtained from equations (\ref{E21-typeI-2}):
\begin{alignat*}{3}
& |\mu_2|^2 =  2_q3_q5_q9_q - |\mu_1|^2 =   \frac{3_q^25_q9_q}{2_q},  \qquad&&
|\nu_1|^2  =  4_q5_q9_q - |\mu_1|^2 =   \frac{5_q^29_q}{2_q}, & \\
& |\nu_2|^2  =  4_q7_q9_q - |\mu_2|^2 =   \frac{7_q9_q}{2_q} , \qquad&&
  |\sigma_2|^2  =  2_q3_q^25_q^2 - |\sigma_1|^2 - |\rho|^2 =   \frac{1}{2}\left(\frac{3_q^25_q^26_q}{2_q4_q} \right). &
\end{alignat*}
As before we have used $q$-numbers identities in order to simplify the expressions. We can also write these expressions in the following way.
The {\it square} of the modulus of a cell  can always be expressed as $f[a,b,c,d]  =  \sqrt{2\,  g[a,b,c,d]}= \sqrt{a+b\sqrt{2} +c \sqrt{3} +d \sqrt{6}}$,  for appropriate rational arguments:
 \begin{alignat*}{3}
&     | \alpha_1 |^2 = f[10, 5, 4, 4], \qquad &&
      | \alpha_2 |^2 = f[272, 191, 156, 111],  &\\
&        | \alpha_3 |^2 = f[5896, 4169, 3404, 2407], \qquad &&
        | \gamma_2 |^2 = f[32, 22, 18, 13], &\\
&| \gamma_3 |^2 = f[686, 485, 396, 280],\qquad &&
          | \nu_1 |^2 = f[1508, 1066, 870, 615], &\\
& | \mu_1 |^2 = f[596, 421, 344, 243],\qquad &&
            | \sigma_1 |^2 = f[2224, 1571, 1284, 907], & \\
& | \nu_2 |^2 = f[118, 83, 68, 48],\qquad & & |              \mu_2 |^2 = f[5120, 3620, 2956, 2090], &\\
& |\sigma_2 |^2 = f[1112, 1571/2, 642, 907/2], \qquad &&   | \beta |^2 = f[2560, 1810, 1478, 1045], & \\
&| \lambda |^2 = f[2948, 4169/2, 1702, 2407/2],\qquad & &
                  | \rho |^2 = f[11002, 15559/2, 6352, 8993/2]. &
\end{alignat*}

\paragraph{Fixing the phases using gauge freedom.}
We f\/irst observe that the ${\mathcal E}_{21}$ fusion graph does not contain double edges, so that the discussion is somewhat easier than with ${\mathcal E}_{9}$.
We consider the graph as built from an inner part (the big rectangle delimited by vertices $5_1^1$, $5_1^2$, $ 5_2^1$, $5_2^2$), and two external ``wings'', each one made of $9$ triangles, each one attached to the rectangular inner part by three triangles sharing edges of type $(6,4)$ and $(4_1^1, 4_1^2)$, or $(4_2^1, 4_2^2)$, with the inner rectangle.
The large inner rectangle contains the (projection) of a central octahedron with vertices $4_1^1$, $7_1$, $4_2^2$, $4_2^1$, $7_2$, $4_1^2$, like it was in the ${\mathcal E}_{5}$ case. This octahedron can be built in three dif\/ferent ways as a~bi-pyramid with a rectangular base, the base being one of the three Type~II  frames displayed on the f\/irst line of~(\ref{rectanglesE21}). The inner rectangle contains also the projection of an upper octahedron with vertices $4_1^1$, $5_1^2$, $5_2^1$, $4_2^2$, $6_1, 7_1$, with its three rectangular sections displayed as Type~II  frames on the second line of~(\ref{rectanglesE21}). Finally it contains its symmetric, the lower octahedron $4_1^2$, $5_1^1$, $5_2^2$, $4_2^1$, $6_2$, $7_2$ with its three rectangular sections appearing on the third line of~(\ref{rectanglesE21}). We end up with $3 \times 3 = 9$  coherence equations for the non-degenerated frames of Type~II  (these equations involve the cell themselves, not their modulus). Each equation will equate to zero a sum of two terms since we deal with bi-pyramids, each term being itself a product of four cells.

All the cells of the same type have the same modulus but they are not necessarily equal. For this reason we need to be specif\/ic about the actual triangles that we consider.
Nevertheless, and as we did when we discussed the ${\mathcal E}_{5}$ case,  we can make real and positive almost all cells of the three octahedra.  At the end we are left with a few possible ambiguities:  for the sake of the discussion it will be enough
to distinguish between  $\rho^{\prime} = {\mathcal T}(7_2, 4_1^1, 4_2^2)$ and $\rho^{\prime\prime} = {\mathcal T}(7_1, 4_2^1, 4_1^2)$, both with modulus  $| \rho |$, and also between $\sigma_2^{\prime} = {\mathcal T}(7_2, 4_2^1, 4_1^2)$ and $\sigma_2^{\prime\prime} = {\mathcal T}(7_1, 4_1^1, 4_2^2)$, both with modulus  $| \sigma |$. All the other cells of the large central rectangle can be made real and positive by simple gauge choices;  in the same way, one can also make $\rho^{\prime}$, $\rho^{\prime\prime}$, $\sigma_2^{\prime}$ and  $\sigma_2^{\prime\prime}$ real, but an ambiguity about their signs still remains.

Any of the three equations of Type~II  associated with the central octahedron leads to the constraint $\sigma_2^{\prime} \rho^{\prime} + \sigma_2^{\prime\prime} \rho^{\prime\prime} = 0$. For instance the one associated with the frame $(4_1^1, 4_2^2, 4_2^1, 4_1^2)$ reads $\frac{1}{[7_1]}\lambda \lambda  \rho^{\prime\prime} \sigma_2^{\prime\prime}  + \frac{1}{[7_2]}\lambda \lambda \rho^{\prime} \sigma_2^{\prime} = 0$, hence the result since $[7_1]=[7_2]$.
This already implies that either one or three among these four cells should be negative.

The three equations associated with the upper octahedron involve only the unknown cell $\sigma_2^{\prime\prime}$, all others being already gauge f\/ixed to a real and positive quantity. For instance, the equation associated with the frame $(4_1^1, 5_1^2, 5_2^1, 4_2^2 )$ reads:  $\frac{1}{[6_1]} \mu_1 \mu_1  \nu_1 \sigma_1 + \frac{1}{[7_1]}  \mu_2 \mu_2  \nu_2  \sigma_2^{\prime\prime} = 0$. From the known values of the moduli, one checks that this equation holds automatically, provides one takes $\sigma_2^{\prime\prime} = - | \sigma_2|$. The other two equations lead to the same conclusion.

The three equations associated with the lower octahedron involve only the unknown cell~$\sigma_2^{\prime}$ and the discussion is similar (the equations, in terms of q-numbers, are exactly the same, by symmetry). The conclusion, here, is that $\sigma_2^{\prime}$ has to be negative as well:  $\sigma_2^{\prime} = - | \sigma_2|$.

We reach the conclusion that all the cells of the large inner rectangle can be made real and positive, with the exception of three of them:  $\sigma_2^{\prime}$, $\sigma_2^{\prime\prime}$ and, say,  $\rho^{\prime\prime}$, that should be negative (of course $\rho^{\prime} = - \rho^{\prime\prime}$)
\begin{equation}
\begin{array}{@{}ccccc}
\basisBigLarge{4_2^1}{4_1^2}{4_1^1}{4_2^2}{}{}{}{} & \qquad \qquad &
 \basisBigLarge{7_1}{4_2^1}{7_2}{4_1^1}{}{}{}{} & \qquad \qquad&
\basisBigLarge{4_1^2}{7_1}{4_2^2}{7_2}{}{}{}{}  \\
{} & \\
\basisBigLarge{4_1^1}{4_2^2}{5_2^1}{5_1^2}{}{}{}{} & \qquad &
 \basisBigLarge{4_2^2}{7_1}{5_1^2}{6_1}{}{}{}{} & \qquad &
\basisBigLarge{6_1}{4_1^1}{7_1}{5_2^1}{}{}{}{} \\
{} & \\
 \basisBigLarge{4_2^1}{4_1^2}{5_1^1}{5_2^2}{}{}{}{} & \qquad &
\basisBigLarge{7_2}{4_2^1}{6_2}{5_1^1}{}{}{}{} & \qquad &
\basisBigLarge{4_1^2}{6_2}{5_2^2}{7_2}{}{}{}{}
\end{array}
\label{rectanglesE21}
\end{equation}

{\samepage Finally, for each of the two wings, we are left  with $9$ triangles on each side. Although three edges in common with the inner rectangle have already been gauged f\/ixed on each side, there is enough gauge freedom (choice of phases)  on the remaining edges to make all of the corresponding cells real and positive; for the f\/irst wing one can for instance start from~$1_1$,~$2_1^1$,~$2_1^2$,  make it real positive by the choice of a phase on the edge~$1_1$,~$2_1^1$ and proceed towards the center while determining an appropriate sequence of independent gauge choices.

}

We already gave the absolute values for the cells for all possible types. Our solution for the system is obtained by taking all the cells real-valued and positive, with the exception of three cells, that should be real and negative:
one of the two cells of type $\rho$, say the one that we called~$\rho^{\prime\prime}$, and the two cells of type $\sigma_2$, namely those that we called $\sigma_2^{\prime}$ and $\sigma_2^{\prime\prime}$.

\subsection[A counter example: the $\mathcal{Z}_9$ graph]{A counter example: the $\boldsymbol{\mathcal{Z}_9}$ graph} \label{section4.7}

\begin{figure}[h]\centering
 \scalebox{0.4}{\includegraphics{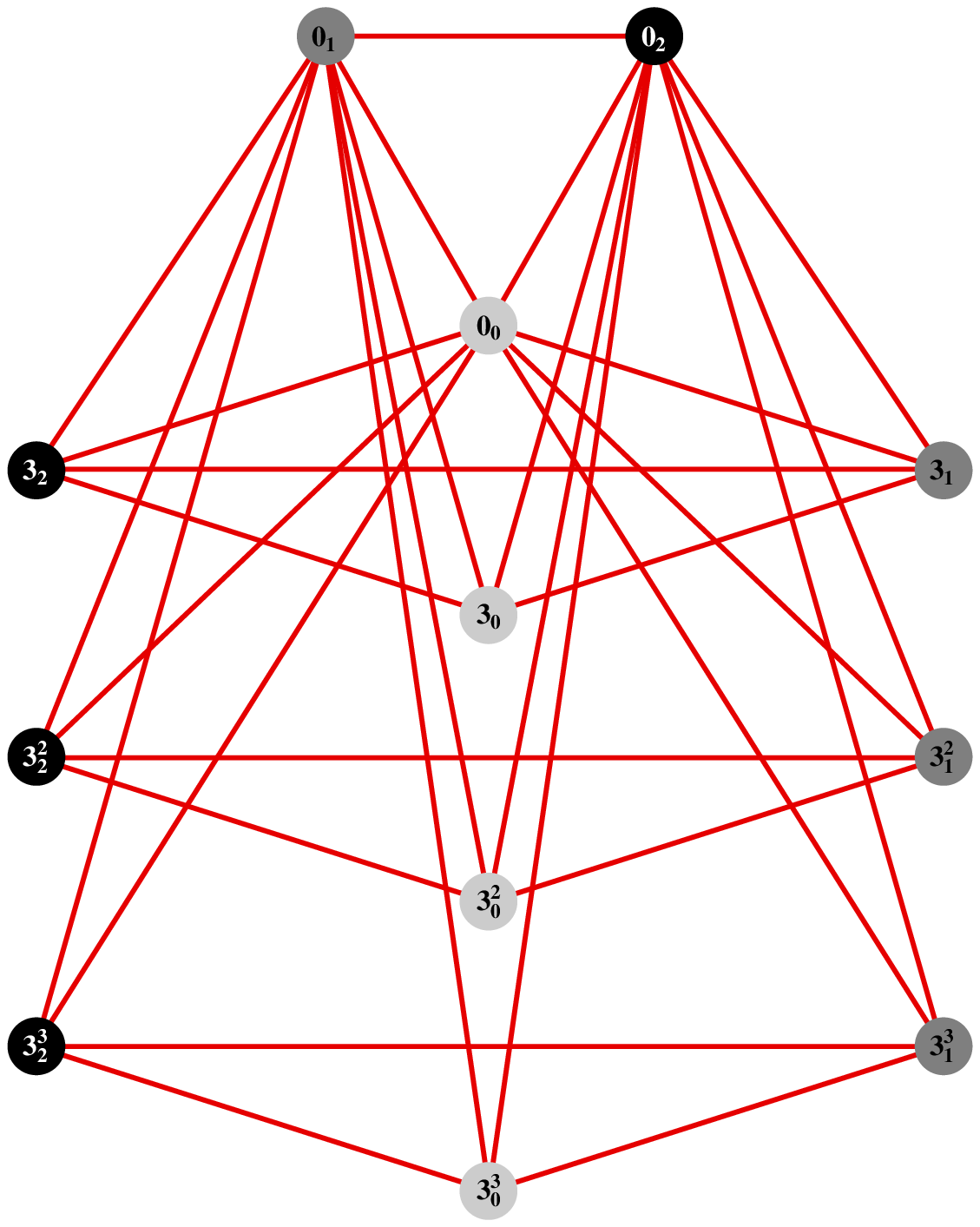}\qquad\qquad \qquad\includegraphics{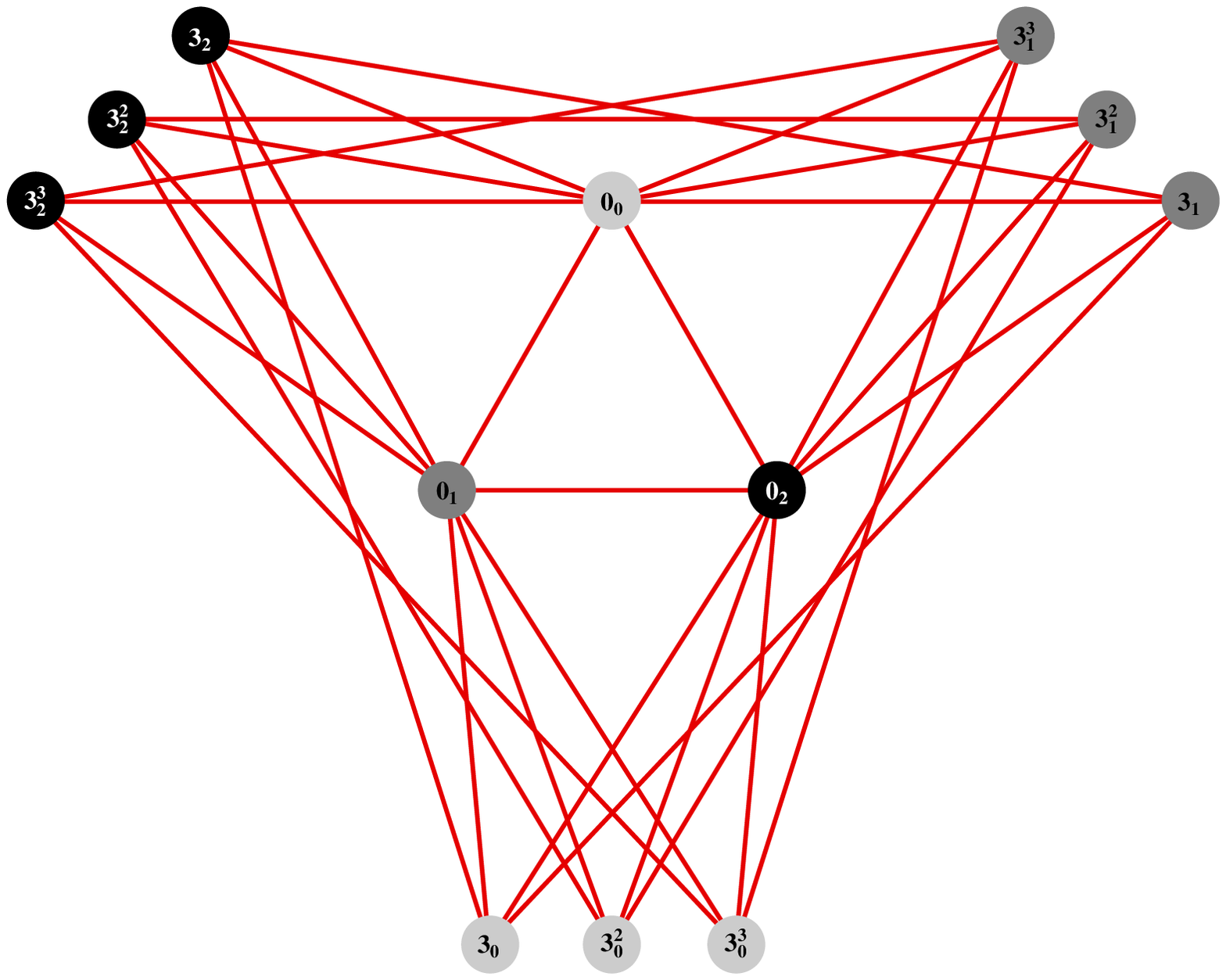}}
\caption{Two equivalent graphical representation of the $\mathcal{Z}_9$ fusion graph, rejected as a fusion graph for a module-category of type $SU(3)$.}
\label{Z9-graph}
\end{figure}

The $\mathcal{Z}_9$ is a graph that could appear at level $9$, i.e.\ at an altitude $\kappa = 12$: $q^{12}= - 1$.
We shall see that coherence equations for triangular cells cannot be satisf\/ied in this case.
The graph is displayed on Fig.~\ref{Z9-graph}.

Vertices are denoted
by $0_i$ and $3_{i}^{j}$ for $i,j = 1,2,3$. The subscript $i$ refers to the vertex triality. Quantum dimension of vertices are $[0_i] = 1$ and $[3_i^j] = 1/\sqrt{3}$, for $i,j=1,2,3$. There are 22 cells for the $\mathcal{Z}_9$ graph. The $\mathbb{Z}_3$ symmetry of the graph suggests that the cell values should also exhibit this $\mathbb{Z}_3$ symmetry. In order to prove the incompatibility of the system, we do not impose, at f\/irst, any symmetry and search for the most general solution to the cell system of $\mathcal{Z}_9$, if any.

There are 30 Type~I  equations for $\mathcal{Z}_9$, involving the square values of the 22 cells. These equations lead to relations between square values of the cells, and we are left with 7 linear equations of Type~I  involving 10 unknown square values. Among the Type~II  equations, 30  correspond to fully degenerated cases. Taking into account the relations obtained previously, this set reduces to~7 second degree equations involving the same 10 unknown square values. One can then check that this system of~7 linear equations and~7 second degree equations involving~10 square values of cells does not have any solution.

Let us be more specif\/ic and show explicitly that one reaches a contradiction, considering a~subset of the system of equations described above involving only 7 square values of cells. We introduce the following notation for the square values of cells:
\begin{gather*}
a = |\mathcal{T}_{0_00_10_2}|^2 ,\qquad   b_j = |\mathcal{T}_{0_13_{2}^{j}0_0}|^2,   \qquad c_j = |\mathcal{T}_{0_13_2^j3_0^j}|^2,  \qquad \textrm{for} \quad j=1,2,3.
\end{gather*}
We consider the following Type~I  frames and the fully degenerated Type~II  frames associated with the same edges:
\begin{gather*}
\frameIsingle{0_1}{3_2^j} \qquad \qquad \qquad \basisBig{0_1}{3_2^j}{0_1}{3_2^j}{}{}{}{}
\end{gather*}
The corresponding Type~I  and Type~II  equations are:
\begin{alignat}{4}
& c_j+b_j = 2_q [0_1][3_2^j] \qquad && \Longleftrightarrow \qquad&&   c_j+b_j   =  \sqrt{\frac{1}{3}(2+\sqrt{3})},&\nonumber\\
& b_j^2 +\frac{1}{[3_2^j]}c_j^2  =  [0_1][0_1][3_2^j]+[0_1][3_2^j][3_2^j] \qquad && \Longleftrightarrow \qquad && b_j^2 + \sqrt{3} c_j^2  =  \frac{1}{3}(1+\sqrt{3}).&
\end{alignat}
The solution $(b_j, c_j)$ of this quadratic equation is
\begin{gather}
b_j =  \displaystyle \frac{1}{6}\big(3\sqrt{2} \pm (3-\sqrt{3})\big), \qquad
c_j = \displaystyle \frac{1}{6}\big(\sqrt{6} \mp (3-\sqrt{3})\big).
\label{z9sol}
\end{gather}
We now consider the following Type~I  frame and the fully degenerated Type~II  frame associated with the same edge:
\begin{gather*}
\frameIsingle{0_0}{0_1} \qquad \qquad \qquad \basisBig{0_0}{0_1}{0_0}{0_1}{}{}{}{}
\end{gather*}
The corresponding Type~I  and Type~II  equations are:
\begin{gather*}
a+ b_1+b_2+b_3  =   2_q = \sqrt{2+\sqrt{3}}, \qquad
a^2 + \sqrt{3} (b_1^2 + b_2^2 + b_3^2 )  =  2.
\end{gather*}
We have therefore two dif\/ferent ways of calculating the same coef\/f\/icient $a$.
Using the solutions for~$b_j$ obtained in (\ref{z9sol}) it is easy to check that the above two equations are incompatible.
As claimed in~\cite{Ocneanu:Bariloche}, the graph $\mathcal{Z}_9$, that appears in the list \cite{DiFZuber}, and leads to a {\it bona fide} collection of annular matrices with non-negative integer coef\/f\/icients providing a representation of the fusion algebra of ${\mathcal A}_9(SU(3))$, has nevertheless to be rejected as fusion graph describing a module-category of type $SU(3)$ since the coherence equations fail.


 \section{Comments}\label{section5}

\paragraph{Coherence equations for others Lie groups.}
 For the Lie groups $SU(n)$, we have $n-1$ fundamental representations $\sigma_1,  \sigma_2, \ldots,\sigma_{n-1}$. All of them already exist at level $k=1$.  Call $\sigma_0$ the trivial representation.
In the case of $SU(3)$,  the main ingredient of the coherence equations is that $\sigma_1^{\otimes 3}$ contains $\sigma_0$. This is replaced for,  $SU(n)$, $n >2$,  or for other Lie groups,
by the fact that,  for appropriate fundamental representations $a$, $b$, $c$, the product $\sigma_a \otimes \sigma_b \otimes \sigma_c$ contains~$\sigma_0$.
For $SU(n)$, one can also use the fact that the conjugate of $\sigma_p$ is $\overline \sigma_p = \sigma_{n-p}$.
 What matters is that a triple point will also appear in the wire diagrams. The dif\/ference with the $SU(3)$ case is that one has to keep track of dif\/ferent fundamental representations (one may introduce dif\/ferent kinds of colors for the wires).
 For the same reason the coherence equations will again involve triangular cells, but the edges of triangles, possibly of dif\/ferent colors, will now refer to distinct representations.

 \subsection*{Acknowledgments}

Part of this paper was written while one of us (R.C.) was a visiting professor at CERN, whose hospitality and support are gratefully acknowledged.

\pdfbookmark[1]{References}{ref}
\LastPageEnding

\end{document}